\documentclass{article}
\usepackage[english]{babel}
\usepackage{booktabs}
\usepackage{graphicx}
\usepackage{sidecap}
\usepackage{hyperref}
\usepackage{amssymb, amsmath}
\usepackage{dsfont}
\usepackage{diagbox,multirow}
\usepackage{hhline}
\usepackage{ulem}
\usepackage[dvipsnames]{xcolor}
\usepackage{array}
\usepackage{colortbl}
\usepackage{subfig}
\usepackage{enumerate}
\usepackage[shortlabels]{enumitem}
\usepackage{makecell}
\usepackage{float}
\usepackage{framed}
\usepackage{varwidth}
\usepackage[export]{adjustbox}
\usepackage{xparse}
\usepackage {tikz}
\usetikzlibrary{arrows,decorations.pathmorphing,backgrounds,positioning,fit,matrix}
\usetikzlibrary{decorations.markings}

\usepackage{amsbsy}

\pgfdeclarelayer{myback}
\pgfsetlayers{myback,background,main}
\definecolor {processblue}{cmyk}{0.96,0,0,0}
\definecolor {darkred}{RGB}{190,11,0}

\tikzset{mycolor/.style = {line width=1bp,color=#1}}%
\tikzset{myfillcolor/.style = {draw,fill=#1}}%

\newlength\mylen
\tikzset{
	bicolor/.style 2 args={
		dashed,dash pattern=on 20pt off 20pt,-,#1,
		postaction={draw,dashed,dash pattern=on 20pt off 20pt,-,#2,dash phase=20pt}
	},
}

\NewDocumentCommand{\fhighlight}{O{blue!40} m m}{%
	\draw[myfillcolor=#1] (#2.north west)rectangle (#3.south east);
}

\usetikzlibrary{fit}
\tikzset{%
	highlight/.style={rectangle,rounded corners,fill=red!15,draw,
		fill opacity=0.5,thick,inner sep=0pt}
}

\newtheorem{thm}{Theorem}[section]
\newtheorem{lem}[thm]{Lemma}
\newtheorem{prop}[thm]{Proposition}
\newtheorem{cor}[thm]{Corollary}
\newtheorem{defn}{Definition}[section]

\newtheorem{rem}{Remark}

\newcommand{\N}{\mathbb{N}}
\newcommand{\R}{\mathbb{R}}

\newcommand{\diag}{\mathop{\mathrm{diag}}}\newcommand{\e}{\textup{e}}\renewcommand{\i}{\textup{i}}\newcommand{\bfzero}{\mathbf 0}\newcommand{\bfuno}{\mathbf 1}

\newcommand{\bfnn}{{\boldsymbol n}}\newcommand{\bfff}{{\boldsymbol f}}\newcommand{\bfgg}{{\boldsymbol g}}\newcommand{\bftt}{{\boldsymbol t}}\newcommand{\bfii}{{\boldsymbol i}}\newcommand{\bfjj}{{\boldsymbol j}}\newcommand{\bfmm}{{\boldsymbol m}}\newcommand{\bfpp}{{\boldsymbol p}}\newcommand{\bfkk}{{\boldsymbol k}}\newcommand{\bfhh}{{\boldsymbol h}}\newcommand{\bfww}{{\boldsymbol w}}\newcommand{\bfxx}{{\boldsymbol x}}
\newcommand{\bftheta}{{\boldsymbol\theta}}\newcommand{\bfnnu}{{\boldsymbol\nu}}\newcommand{\bfss}{{\boldsymbol s}}

\newcommand*\quot[2]{{^{\textstyle #1}\big/_{\textstyle #2}}}
\DeclareMathSymbol{\shortminus}{\mathbin}{AMSa}{"39}

\begin{document}
\title{Asymptotic spectra of large (grid) graphs with a uniform local structure\footnote{This is a preprint.}}
\author{
Andrea Adriani $^{(a)}$, Davide Bianchi $^{(b)}$, Stefano Serra-Capizzano $^{(c,d)}$}

\date{}

\maketitle

\section*{Abstract}
We are concerned with sequences of graphs having a grid geometry,
with a uniform local structure in a bounded domain $\Omega\subset {\mathbb
R}^d$, $d\ge 1$. We assume $\Omega$ to be Lebesgue measurable with regular boundary
and contained, for convenience, in the cube $[0,1]^d$. When
$\Omega=[0,1]$, such graphs include the standard Toeplitz graphs
and, for $\Omega=[0,1]^d$,  the considered class includes
$d$-level Toeplitz graphs. In the general case, the underlying
sequence of adjacency matrices has a canonical eigenvalue
distribution, in the Weyl sense, and we show that we can associate
to it a symbol $f$. The knowledge of the symbol and of its basic
analytical features provide many informations on the eigenvalue
structure, of localization, spectral gap, clustering, and distribution type. Few generalizations are also considered in
connection with the notion of generalized locally Toeplitz
sequences and applications are discussed, stemming e.g. from the
approximation of differential operators via numerical schemes.\\ 
\ \\
\noindent
(a) \, Department of Science and High Technology,   University of Insubria,  Via Valleggio 11, 22100 Como, Italy\\
(aadriani@uninsubria.it)\\
(b)\, Department of Science and High Technology,   University of Insubria,  Via Valleggio 11, 22100 Como, Italy\\
(d.bianchi9@uninsubria.it)\\
(c)\, Department of Science and High Technology,   University of Insubria,  Via Valleggio 11, 22100 Como, Italy\\
 (stefano.serrac@uninsubria.it)\\
(d)\, Department of Information Technology, Uppsala University, Uppsala, Sweeden \\ (stefano.serra@it.uu.se)\\


\section{Introduction}
In this work we  are interested in defining and studying a large
class of graphs enjoying few structural properties:
\begin{description} 
\item[a)] when we
look at them from ``far away'', they should reconstruct
approximately a given domain $\Omega\subset [0,1]^d$, $d\ge 1$,
i.e., the larger is the number of the nodes the more accurate is
the reconstruction of $\Omega$; 
\item[b)] when we look at them
``locally'', that is from a generic internal node, we want that
the structure is uniform, i.e., we should to be unable to understand
where we are in the graphs, except possibly when the considered
node is close enough to the boundaries of $\Omega$.
\end{description} 

 Technically, we are not concerned with a single graph, but with a whole
sequence of graphs, where $\Omega$ and the internal structure are
fixed, independently of the index (or multi-index) of the graph
uniquely related to the cardinality of nodes: thus the resulting
sequence of graphs has a grid geometry, with a uniform local
structure, in a bounded domain $\Omega\subset {\mathbb R}^d$, $d\ge 1$.   We assume the domain $\Omega$ to be Lebesgue measurable with regular boundary, which is for us a boundary $\partial \Omega$ of zero Lebesgue measure,
and contained for convenience in the cube $[0,1]^d$. We will call {\itshape regular} such a domain. When $\Omega=[0,1]$, it is worth observing that such graphs include the standard
Toeplitz graphs (see \cite{Ghorban2012} and Definition \ref{def:toeplitz-graph}) and for $\Omega=[0,1]^d$  the
considered class includes $d$-level Toeplitz graphs (see Definition \ref{def:d-toeplitz-graph}). 

The main result is the following: given a sequence of graphs
having  a grid geometry with a uniform local structure in a domain $\Omega$, the underlying sequence
of adjacency matrices has a canonical eigenvalue distribution, in
the Weyl sense (see \cite{GS,BS} and references therein), and we
show that we can associate to it a symbol function $f$. More precisely,
when $f$ is smooth enough, if $N$ denotes the size of the
adjacency matrix (i.e. the number of the nodes of the graph), then
the eigenvalues of the adjacency matrix are approximately values
of a uniform sampling of $f$ in its definition domain, which
depends on $\Omega$ (see Definition \ref{def:eig-distribution} for
the formal definition of eigenvalue distribution in the Weyl sense
and the results on Section \ref{sec:main} for the precise
characterization of $f$ and of its definition domain).

 The knowledge of the symbol and of some of its basic analytical features provide a
 lot of information on the eigenvalue structure, of localization, spectral gap, clustering, and distribution type. 

 The mathematical tools are taken from the field of Toeplitz (see
 the rich book by B\"ottcher and Silbermann \cite{BS} and \cite{GS,Tillinota,tyrtL1}) and
 Generalized Locally Toeplitz (GLT) matrix sequences (see  \cite{Tilliloc,glt,glt-bis}): for a recent account on the GLT theory,
which is indeed quite related to the present topic, we refer to the following books and reviews \cite{glt-book-1,glt-book-2,glt-book-3,GMS18}. Interestingly enough, as discussed
 at the end of this paper, many numerical schemes (see e.g. \cite{FE-book,IgA-book,FD}) for
 approximating partial differential equations (PDEs) and operators lead
 to sequences of structured matrices which can be written as linear combination  of adjacency 
matrices, associated with the graph  sequences described here. More specifically, if the physical domain of the differential
operator is $[0,1]^d$ (or any $d$-dimensional rectangle) and the coefficients are constant, then
we encounter $d$-level (weighted) Toeplitz graphs, when approximating the underlying PDE by using e.g. equispaced Finite Differences or uniform Isogeometric Analysis (IgA). On the other hand, under the same assumptions on the underlying operator, quadrangular  and triangular Finite Elements lead to block $d$-level Toeplitz structures, where the size of the blocks is related to the degree of the polynomial space of approximation and to the dimensionality $d$ (see \cite{FEM-paper}). Finally, in more generality the GLT case is encountered by using any of the above numerical techniques, also with non-equispaced nodes/triangulations, when dealing either with a general domain $\Omega$ or when the coefficients of the differential operator are not constant.

The paper is organized as follows. In Section \ref{sec:intro} we
collect all the machinery we need for our derivations: we will
first review basic definitions and notations from graph theory,
from the field  of Toeplitz and $d$-level Toeplitz matrices, and
then we provide the definitions of canonical spectral
distribution, spectral clustering etc. In Section \ref{sec:main}
we give the formal definitions of sequences of graphs  having  a
grid geometry, with a uniform local structure, in regular domains $\Omega\subset  [0,1]^d$, $d\ge 1$, and we
prove the main results, by identifying the related symbols.
Section \ref{sec:appl} contains specific applications, including
the analysis of the spectral gaps and the connections with the
numerical approximation of differential
operators by local methods, such Finite Differences, Finite Elements, Isogeometric
Analysis etc. Finally,  Section \ref{sec:final} is devoted to draw
conclusions and to present open problems.

\section{Background and definitions}\label{sec:intro}

In this section we present some definitions, notations, and
(spectral) properties associated with graphs (see \cite{sp-gr} and
references therein) and, in particular, with Toeplitz graphs
\cite{Ghorban2012}.

\begin{defn}[Standard set]
Given $n\in \N$, we call $[n]:=\{1,\ldots,n\}$ the \textnormal{standard set} of cardinality $n$.
\end{defn}

\begin{defn}[Graph]
We will call a (finite) \textnormal{graph} the quadruple $G=(V,E,w, \kappa)$, defined by a
set of nodes 
$$
V=\left\{ v_1, v_2,\ldots, v_n \right\},
$$
a weight function $w: V\times V \to \R$, a set of edges 
$$
E=\{(v_i,v_j)|\, v_i,v_j\in V, \; w(v_i,v_j)\neq0\}
$$
between the nodes and a \textnormal{potential term} $\kappa :V \to \R$. The non-zero values $w(v_i,v_j)$ of the weight function $w$ are called \textnormal{weights} associated to the edge $(v_i,v_j)$. Given an edge $e=(v_i,v_j)\in E$, the nodes $v_i,v_j$ are called \textnormal{end-nodes} for the edge $e$. An edge $e\in E$ is said to be \textnormal{incident} to a node $v_i\in V$ if there exists a node $v_j\neq v_i$ such that
either $e=(v_i,v_j)$ or $e=(v_j,v_i)$. A \textnormal{walk} of length $k$ in $G$ is a set of nodes $v_{i_1}, v_{i_2},\ldots,v_{i_k}, v_{i_{k+1}}$ such
that for all $1\leq r\leq k$, $(v_{i_r},v_{i_{r+1}})\in E$. A \textnormal{closed walk} is a walk for which $v_{i_1}=v_{i_{k+1}}$. A \textnormal{path}
is a walk with no repeated nodes. A graph is \textnormal{connected}
if there is a walk connecting every pair of nodes. 

A graph is said to be \textnormal{unweighted} if $w(v_i,v_j)\in \{0,1\}$ for every $v_i,v_j\in V$. In that case the weight function $w$ is uniquely determined by the edges which belong to $E$. 

A graph is said to be \textnormal{undirected} if the weight function $w$ is symmetric, i.e.,  for every edge $(v_i,v_j)\in E$ then $(v_j,v_i)\in E$ and $w(v_i,v_j)=w(v_j,v_i)$. In this case the edges $(v_i,v_j)$ and $(v_j,v_i)$ are considered equivalent and the edges are formed by unordered pairs of vertices. Two nodes $v_i,v_j$ of an undirected graph are said to be \textnormal{neighbors} if $(v_i,v_j)\in E$ and we will write $v_i\sim v_j$. On the contrary, if $(v_i,v_j)\notin E$ then we will write $v_i\nsim v_j$.

An undirected graph with unweighted edges and no self-loops (edges from a node to itself) is said to be \textnormal{simple}. 

Throughout this work, we will always consider undirected
graphs without self-loops. 

Since the potential term $\kappa$ will come into play only in Section \ref{sec:appl}, we will often omit it in the next subsections and in Section \ref{sec:main}. Moreover, when dealing with simple graphs we will use the simplified notation $G=(V,E)$.
\end{defn}

\begin{defn}[Sub-graph, interior and boundary nodes]\label{def:sub-graph}
Given a graph $\bar{G}=(\bar{V},\bar{E},\bar{w},\bar{\kappa})$ we say that a graph $G=(V,E,w,\kappa)$ is a (proper) \textnormal{sub-graph} of $\bar{G}$, and we write $G\subset \bar{G}$, if
\begin{enumerate}[(i)]
	\item $V\subset \bar{V}$;
	\item $E=\{(v_i,v_j) \in \bar{E} \, | \, v_i,v_j \in V \}\subset \bar{E}$;
	\item $w=\bar{w}_{|E}$;
	\item $\kappa(v_i) = \bar{\kappa}(v_i)$ for every node $v_i$ such that there not exists $\bar{v}_j \in \bar{V}\setminus V$, $\bar{v}_j \sim v_i$. \label{item:kappa=kappa}
\end{enumerate} 
Sometimes we will call $\bar{G}$ the \textnormal{mother graph}. The set of nodes 
$$
\mathring{V}:=\left\{ v_i \in V \, | \, v_i \nsim \bar{v}_j \, \forall \, \bar{v}_j \in \bar{V}\setminus V \right\}
$$
is called \textnormal{interior} of $V$ and its elements are called \textnormal{interior nodes}. Vice-versa, the set of nodes
$$
\partial V:=\left\{ v_i \in V \, | \, v_i \sim \bar{v}_j \, \mbox{for some } \, \bar{v}_j \in \bar{V}\setminus V \right\}
$$
is called \textnormal{boundary} of $V$ and its elements are called \textnormal{boundary nodes}. Therefore, condition \eqref{item:kappa=kappa} can be restated saying that $\kappa = \bar{\kappa}_{|\mathring{V}}$. Observe that we do not request that $\kappa =\bar{\kappa}$ on the boundary of $V$. 
\end{defn}

\begin{defn}[Degree of a node]
	In an undirected graph, the {\textnormal{degree}} of a node $v_i\in V$, denoted by  $\deg(v_i)$, is the sum of weights associated to the edges incident to $v_i$, that is,
	\begin{equation*}
\deg(v_i):= \sum_{v_j\sim v_i}w(v_i,v_j).
	\end{equation*} 
\end{defn}

\begin{defn}[Adjacency matrix]
Every graph $G=(V,E,w,\kappa)$ with $\kappa\equiv 0$ can be represented as a matrix $W=\left(w_{i,j}\right)_{i,j=1}^n\in\mathbb{R}^{n\times n}$, called the \textnormal{adjacency matrix} of the graph. In particular, there is a bijection between the set of weight functions $w : V\times V \to \R$ and the set of a adjacency matrices $W \in \R^{|V|\times |V|}$.

The entries of the adjacency matrix $W$ are
\begin{equation}
w_{i,j} =w(v_i,v_j), \qquad \forall \, v_i,v_j\in V.
\end{equation}
In short, the adjacency matrix tells which nodes are connected and the \textquoteleft weight\textquoteright of the connection. If the graph does not admit self-loops, then the diagonal elements of the
adjacency matrix are all equal to zero. In the particular case of
an undirected graph, the associated adjacency matrix is symmetric,
and thus its eigenvalues are real \cite{bhatia}.
\end{defn}

We will always label the eigenvalues in non-decreasing order: $\lambda_1\leq\lambda_2\leq\ldots\leq\lambda_n$.

	\begin{defn}[Isomorphism between graphs]
		Given two graphs $G=(V,E,w,\kappa), G'=(V',E',w',\kappa')$ with
		$$
		V=\{v_1,\ldots, v_n\}, \qquad V'=\{ v_1',\ldots, v_m'\},
		$$ 
		we say that $G$ is isomorph to $G'$, and we write $G \simeq G'$, if 
		\begin{enumerate}[i)]
			\item $n=m$, i.e., $|V|=|V'|$ where $|\cdot|$ is the cardinality of a set;
			\item there exists a permutation $P$ over the standard set $[n]$ such that $w(v_i,v_j)=w'(v'_{P(i)}, v'_{P(j)}), \kappa(v_i)=\kappa'(v'_{p(i)})$.
		\end{enumerate}
		In short, two graphs are isomorphic if they contain the same number of vertices connected in the same way. Notice that an isomorphism between graphs is characterized by the permutation matrix $P$.
\end{defn}

As an immediate consequence of the previous definition, it holds that $G\simeq G'$ if and only if there exists a permutation matrix $P$ such that $W=PW'P^{-1}=PW'P^{T}$, where $W, W'$ are the adjacency matrices of $G$ and $G'$, respectively.

\begin{defn}[Linking-graph operator]\label{def:linking-graph}
Given $\nu\in \N$, we will call {\itshape linking-graph operator} for the reference node set $[\nu]$ any non-zero $\R^{\nu\times \nu}$ matrix, and we will indicate it with $L$. Namely, a linking-graph operator is the adjacency matrix for a (possibly not undirected) graph $G=([\nu],E,l)$, with $l$ a weight function. Trivially, note that $L$ may have nonzero elements on the main diagonal, so it admits loops. When the entries of $L$ are just in $\{0,1\}$ then we call it a simple linking-graph operator.
\end{defn}

In Section \ref{sec:main}, we will use the linking-graph operator to connect a (infinite) sequence of graphs $G_1\simeq G_2\simeq\ldots\simeq G_n$, and that will grant a uniform local structure on the graph $G:= \bigcup_{n=1}^\infty G_n$.

The set of real functions on $V$ will be denoted as $C(V)$. Trivially, $C(V)$ is isomorph to $\R^n$.  Of great importance for Section \ref{sec:appl} will be the operator $\Delta_{G} : C(V) \to C(V)$ defined below.

\begin{defn}[Graph-Laplacian]\label{def:graph-Laplacian}
Given an undirected graph with no loops $G=(V,E,w,\kappa)$, the \textnormal{graph-Laplacian} is the symmetric matrix $\Delta_G : C(V) \to C(V)$ defined as
\begin{equation*}
\Delta_G := D+K - W,
\end{equation*} 
where $D$ is the \textnormal{degree matrix} and $K$ is the \textnormal{potential term matrix}, that is,
$$
D:=\diag\left\{\deg(v_1),\deg(v_2),\ldots,\deg(v_n)\right\}, \qquad K:=\diag\left\{\kappa(v_1),\kappa(v_2),\ldots,\kappa(v_n)\right\},
$$ and $W$ is the adjacency matrix of the graph $G$, that is,
$$
W=\left(
\begin{array}{cccc}
0      &   w(v_1,v_2)  & \cdots & w(v_1,v_n)\\
w(v_1,v_2)    &    0   & \ddots & \vdots \\
\vdots & \ddots & \ddots & w(v_{n-1},v_n) \\
w(v_1,v_n) & \cdots & w(v_{n-1},v_n)     & 0\\
\end{array}
\right).
$$
Namely,
\begin{equation*}
\Delta_G= \left(
\begin{array}{cccc}
\deg(v_1)  +\kappa(v_1)     &   -w(v_1,v_2)  & \cdots & -w(v_1,v_n)\\
-w(v_1,v_2)    &    \deg(v_2) + \kappa(v_2)   & \ddots & \vdots \\
\vdots & \ddots & \ddots & -w(v_{n-1},v_n) \\
-w(v_1,v_n) & \cdots & -w(v_{n-1},v_n)     & \deg(v_n)+\kappa(v_n)\\
\end{array}
\right).
\end{equation*} 
\end{defn}

\subsection{Toeplitz matrices, $d$-level Toeplitz matrices, and symbol}\label{sec:toeplitz_matrices}

Toeplitz matrices $A_n$ are characterized by the fact that all their diagonals parallel to the main diagonal have constant values
$a_{ij}=a_{i-j}$, where  $i,j=1,\ldots,n$, for given coefficients $a_k$, $k=1-n,\ldots,n-1$:
$$A_n=\left(
\begin{array}{cccc}
a_0      &   a_{-1}  & \cdots & a_{1-n}\\
a_1    &    a_0   & \ddots & \vdots \\
\vdots & \ddots & \ddots & a_{-1} \\
a_{n-1} & \cdots & a_1     & a_0\\
\end{array}
\right).$$
When every term $a_k$ is a matrix of fixed size $\nu$, i.e., $a_k\in \mathbb C^{\nu\times \nu}$, the matrix $A_n$ is of block Toeplitz type. Owing to its intrinsic recursive nature, the definition of $d$-level (block) Toeplitz matrices is definitely more involved. More precisely, a $d$-level Toeplitz matrix is a Toeplitz matrix where each coefficient $a_k$ denotes a $(d-1)$-level Toeplitz matrix and so on in a recursive manner.
In a more formal detailed way, using a standard multi-index notation (see \cite{Tyrty} and Remark \ref{rem:multi-index} at the end of this section), a $d$-level Toeplitz matrix is of the form
$$A_\bfnn=\left(a_{\bfii-\bfjj}\right)_{\bfii,\bfjj=\bfuno}^{\bfnn}\in\mathbb C^{(n_1\cdots n_d)\times(n_1\cdots n_d)},$$
with the multi-index $\bfnn$ such that $\bfzero<\bfnn=(n_1,\ldots,n_d)$ and $a_\bfkk\in\mathbb C$,  $-(\bfnn-\bfuno)\trianglelefteq \bfkk\trianglelefteq \bfnn-\bfuno$.
If the basic elements $a_\bfkk$ denote blocks of a fixed size $\nu\ge 2$, i.e. 
 $a_\bfkk\in\mathbb C^{\nu\times \nu}$, then $A_{\bfnn,\nu}$ is a $d$-level block Toeplitz matrix,
 
 $$A_{\bfnn,\nu}=\left(a_{\bfii-\bfjj}\right)_{\bfii,\bfjj=\bfuno}^{\bfnn}\in\mathbb C^{(n_1\cdots n_d\nu)\times(n_1\cdots n_d\nu)}, \qquad a_\bfkk\in\mathbb{C}^{\nu\times \nu}.
 $$
For the sake of simplicity, we write down an example explicitly with $d=2$ and $\nu=3$:
\[
A_{\bfnn,3}=\left(
\begin{array}{cccc}
A_0      &   A_{\shortminus 1}  & \cdots & A_{1 \shortminus n_1}\\
A_1    &    A_0   & \ddots & \vdots \\
\vdots & \ddots & \ddots & A_{\shortminus 1} \\
A_{n_1\shortminus 1} & \cdots & A_1     & A_0\\
\end{array}
\right),\ \ \ 
A_{k_1}= \left(
\begin{array}{cccc}
a_{k_1,0}     &   a_{k_1,\shortminus 1}  & \cdots & a_{k_1,1\shortminus n_2}\\
a_{k_1,1}   &    a_{k_1,0} & \ddots & \vdots \\
\vdots & \ddots & \ddots & a_{k_1,\shortminus 1} \\
a_{k_1,n_2\shortminus 1} & \cdots & a_{k_1,1}      & a_{k_1,0}\
\end{array}
\right),\ \ 
\]
\[ 
a_{k_1,k_2}\in \mathbb C^{3},\ \ 
k_1\in \{1-n_1,\ldots,n_1-1\}, \ \ k_2\in \{1-n_2,\ldots,n_2-1\}.
\]
Observe that each block $A_{k_1}$ has a (block) Toeplitz structure. When $\nu=1$ then we will just write $A_{\bfnn,\nu}= A_{\bfnn}$.

Here we are interested in asymptotic results and thus it is important to a have a meaningful way for defining sequences of Toeplitz
matrices, enjoying global common properties.  A classical and successful possibility  is given by the use of a fixed function, called {\itshape the generating function}, and by taking its Fourier coefficients as entries of all the matrices in the sequence.

More specifically, given a function $\bfff:[-\pi,\pi]^d\to\mathbb C^{\nu\times \nu}$ belonging to $\textnormal{L}^1([-\pi,\pi]^d)$, we denote its Fourier coefficients by
\begin{equation}\label{fourier}
\hat{\bfff}_\bfkk=\frac1{(2\pi)^d}\int_{[-\pi,\pi]^d}\bfff(\bftheta)\e^{-\i\,\bfkk\cdot\bftheta}d\bftheta\in\mathbb C^{\nu\times \nu},\qquad\bfkk\in\mathbb Z^d, \quad \bfkk\cdot\bftheta=\sum_{r=1}^d k_r \theta_r,
\end{equation}
(the integrals are done component-wise), and we associate to $\bfff$ the family of $d$-level block Toeplitz matrices
\begin{equation}\label{toeplitz-symbol}
T_{\bfnn,\nu}(\bfff):=\left(\hat{\bfff}_{\bfii-\bfjj}\right)_{\bfii,\bfjj=\bfuno}^{\bfnn},\qquad\bfnn\in\mathbb N^d.
\end{equation}
We call $\{T_{\bfnn,\nu}(\bfff)\}_\bfnn$ the family of multilevel block
Toeplitz matrices associated with the function $\bfff$,  which is
called the generating function of $\{T_{\bfnn,\nu}(\bfff)\}_\bfnn$. If $\bfff$ is Hermitian-valued, i.e. $\bfff(\bftheta)$ is Hermitian for almost every $\bftheta$, then it is plain to see that all the matrices $T_{\bfnn,\nu}(\bfff)$ are Hermitian, simply because the Hermitian character of the generating function
and relations (\ref{fourier}) imply that $\hat{\bfff}_{-\bfkk}=\hat{\bfff}_\bfkk^H$ for all $\bfkk\in\mathbb Z^d$.
 If, in addition, $\bfff(\bftheta)=\bfff(|\bftheta|)$ for every $\bftheta$, then all the matrices
$T_\bfnn(\bfff)$ are real symmetric with real symmetric blocks $\hat{\bfff}_\bfkk$, $\bfkk\in\mathbb Z^d$.

\begin{rem}[Multi-index notation]\label{rem:multi-index}
Given an integer $d\geq 1$, a $d$-index $\bfkk$ is an element of $\mathbb{Z}^d$, that is, $\bfkk=\left(k_1, \ldots, k_d\right)$ with $k_r \in \mathbb{Z}$ for every $r=1,\ldots, d$. Through this paper we will intend $\mathbb{Z}$ equipped with the lexicographic ordering, that is, given two $d$-indices $\bfii=(i_1,\ldots,i_d)$, $\bfjj=(j_1,\ldots,j_d)$, then we write $\bfii \vartriangleleft \bfjj$ if $i_r<j_r$ for the first $r=1,2,\ldots,d$ such that $i_r\neq j_r$. The relations $\trianglelefteq,\vartriangleright,\trianglerighteq$ are defined accordingly.
	
Given two $d$-indices $\bfii,\bfjj$, we write $\bfii < \bfjj$ if $i_r< j_r$ for every $r=1,\ldots,d$. The relations $\leq,>,\geq$ are defined accordingly.

We indicate with $\bf{0},\bf{1}$ the $d$-dimensional constant vectors $\left(0,0,\ldots,0\right)$ and $\left(1,1,\ldots,1\right)$, respectively. We write $|\bfii|$ for the vector $\left(|i_1|, \ldots, |i_d|\right)$. Finally, given a $d$-index $\bfnn$ we write $\bfnn \to \infty$ meaning that $\min_{r=1,\ldots,d}\{n_r\}\to \infty$.
\end{rem}

\subsection{Weyl eigenvalue distribution and clustering}\label{ssec:weyl_distribution}

We say that a matrix-valued function
$\bfgg:D\to\mathbb C^{\nu\times \nu}$, $\nu\geq 1$, defined on a measurable set
$D\subseteq\mathbb R^m$, is measurable (resp. continuous, in
$\textnormal{L}^p(D)$) if its components $g_{ij}:D\to\mathbb C,\ i,j=1,\ldots,
\nu,$ are measurable (resp. continuous, in $\textnormal{L}^p(D)$). Let
$\mu_m$ be the Lebesgue measure on $\mathbb R^m$ and
let $C_c(\mathbb C)$ be the set of continuous functions with
bounded support defined over $\mathbb C$. Setting $d_\bfnn$ the dimension of the square matrix $X_{\bfnn,\nu}$, for $F\in C_c(\mathbb{C})$ we define $\Sigma_\lambda(F,X_{\bfnn,\nu}):=\frac{1}{d_\bfnn}\sum_{k=1}^{d_\bfnn}F(\lambda_k(X_{\bfnn,\nu}))$.

Hereafter, with the symbol $\{X_{\bfnn,\nu}\}$ we will indicate a sequence of square matrices with increasing dimensions, i.e., such that $d_\bfnn \to \infty$ as $\bfnn \to \infty$, with $\nu$ fixed and independent of $\bfnn$.

\begin{defn}[Eigenvalue distribution of a sequence of matrices]\label{def:eig-distribution}
	Let $\{X_{\bfnn,\nu}\}_\bfnn$ be a sequence of matrices and let $\bfgg:D\to\mathbb C^{\nu\times \nu}$ be a measurable Hermitian matrix-valued function defined on the measurable set
	$D\subset\mathbb R^m$, with $0<\mu_m(D)<\infty$.
	
	 We	say that $\{X_{\bfnn,\nu}\}_\bfnn$ is distributed like $\bfgg$ in the sense of the
	eigenvalues, in symbols $\{X_{\bfnn,\nu}\}_\bfnn \sim_\lambda \bfgg$, if
	\begin{equation}\label{distribution:sv-eig}
	\lim_{\bfnn\rightarrow \infty}\Sigma_{\lambda}(F,X_{\bfnn,\nu})=\frac1{\mu_m(D)}\int_D \sum_{k=1}^\nu  F\left(\lambda_k(\bfgg(\boldsymbol{y}))\right)d\mu_m(\boldsymbol{y}),\qquad\forall F\in C_c(\mathbb C),
	\end{equation}
	where $\lambda_1(\bfgg(\boldsymbol{y})), \ldots,\lambda_\nu(\bfgg(\boldsymbol{y}))$ are the eigenvalues of $\bfgg(\boldsymbol{y})$. Let us notice that in the case $\nu=1$, then identity \eqref{distribution:sv-eig} reduces to
	$$
	\lim_{\bfnn\rightarrow \infty}\Sigma_{\lambda}(F,X_{\bfnn,\nu})=\frac1{\mu_m(D)}\int_D   F\left(g(\boldsymbol{y})\right)d\mu_m(\boldsymbol{y}),\qquad\forall F\in C_c(\mathbb C).
	$$
	We will call $\bfgg$ the \textnormal{(spectral) symbol} of $\{X_{\bfnn,\nu}\}_\bfnn$.
\end{defn}

The following result on Toeplitz matrix sequences linking the definition of symbol function and generating function is due to P. Tilli.
\begin{thm}[\cite{Tillinota}]\label{thm:symbol_d-block-toeplitz}
Given a function $\bfff:[-\pi,\pi]^d\to\mathbb C^{\nu\times \nu}$ belonging to $\textnormal{L}^1([-\pi,\pi]^d)$, then
$$
\left\{T_{\bfnn,\nu}(\bfff)  \right\}_\bfnn \sim_\lambda \bfff,
$$
that is the generating function of $\left\{T_{\bfnn,\nu}(\bfff)  \right\}_\bfnn$ coincides with its symbol according to Definition \ref{def:eig-distribution}.
\end{thm}

Since in this paper we work only with undirected (i.e., symmetric) graphs, and in virtue of the applications in Section \ref{sec:appl}, we will deal with symbol functions $\bfgg$ such that $\lambda_k\left(\bfgg(\boldsymbol{y})\right)$ are real valued for every $\boldsymbol{y}\in D$, and for every $k=1,\ldots, \nu$. See for example Propositions \ref{prop:d-level-distribution},\ref{prop:diamond-distribution} and Theorem \ref{main-theorem}. 

The knowledge of the symbol function $\bfgg$ can give valuable insights on the distribution of the eigenvalues of a sequence of matrices. We refer to the Appendix \ref{appendix} where a collection of theoretical results are presented, and to Section \ref{sec:appl} where numerical experiments are provided.

\section{Diamond Toeplitz graphs}\label{sec:Diamond_Toeplitz}
In this section we are going to present the main (local) graph-structure which will be used to build more general graphs as union of sequences of sub-graphs, i.e., diamond Toeplitz graphs. The resulting graphs will be then immersed in bounded regular domains of $\R^d$ in Section \ref{sec:main}. We will proceed step by step, gradually increasing the complexity of the graph structure. 

As a matter of reference, we have the following scheme of inclusions, with the related variable coefficient versions:

\begin{figure}[H]
	\centering
\begin{tikzpicture}[baseline=-\the\dimexpr\fontdimen22\textfont2\relax]
\matrix (m)[matrix of math nodes, every node/.style={anchor=base}]
{       adjacency matrix of  Toeplitz graph (Def. \ref{def:toeplitz-graph}) & $\subset$  & Toeplitz matrix \\
		\rotatebox{270}{$\subset$}	&           &       \rotatebox{270}{$\subset$}    \\
	adjacency matrix of $d$-level Toeplitz graph (Def. \ref{def:d-toeplitz-graph})   &  $\subset$ & $d$-level Toeplitz matrix \\
		\rotatebox{270}{$\subset$}	&           &       \rotatebox{270}{$\subset$}    \\
	adjacency matrix of	$d$-level diamond Toeplitz graph (Def. \ref{def:d-level-diamond})   &  $\subset$ & $d$-level block Toeplitz matrix \\
};
\end{tikzpicture}
\end{figure}
 
\subsection{Toeplitz graphs and $d$-level Toeplitz graphs}
We first focus on a particular type of graphs, namely Toeplitz graphs.
These are graphs whose adjacency matrices are Toeplitz matrices.
\begin{defn}[Toeplitz graph]\label{def:toeplitz-graph}
Let $n,m,t_1,\ldots,t_m$ be positive integers such that $0< t_1<t_2<\ldots<t_m< n$, and fix $m$ nonzero real numbers $w_{t_1},\ldots,w_{t_m}$. A \textnormal{Toeplitz graph}, denoted by $T_n\langle (t_1,w_{t_1}),\ldots,(t_m,w_{t_m})\rangle$,
is an undirected graph defined by a node set $V_n=\{v_1,\ldots,v_n\}$ and a weight function $w$ such that
$$
w(v_i,v_j) = \begin{cases}
w_{t_k} & \mbox{if } |i-j|=t_k,\\
0 & \mbox{otherwise}.
\end{cases}
$$
In the case of simple graphs, i.e., $w_{t_k}=1$ for every $k$, then we will indicate the Toeplitz graph just as $T_n\langle t_1,\ldots,t_m\rangle$. The number of edges in a Toeplitz graph is equal to $\sum_{k=1}^m(n-t_k)$. By construction, the adjacency matrix $W_n=(w_{i,j})_{i,j=1}^n$ of a  Toeplitz graph has a symmetric Toeplitz structure, i.e., $w_{i,j}=w_{|i-j|}$.
\end{defn}

If we assume that $m$, $t_1,\ldots,t_m$, are fixed (independent of
$n$) and we let the size $n$ grow, then  the sequence of adjacency
matrices $W_n$ can be related to a  unique real integrable function $f$
(the {\itshape symbol}) defined on $[-\pi,\pi]$ and expanded
periodically on $\mathbb{R}$. In this case, according to
(\ref{fourier}), the entries $w_{ij}=\hat{f}_{i-j}$ of the matrix $W_n$ are
defined via the Fourier coefficients of $f$, where the $k$-th
Fourier coefficient of $f$ is given by 
\[
\hat{f}_k=\frac1{2\pi}\int_{-\pi}^\pi
f(\theta)e^{-\i k \theta}d\theta, \qquad  k\in\mathbb Z.
\]

We know that the
Fourier coefficients $\hat{f}_k$ are all in $\{0,w_{t_1},\ldots,w_{t_m}\}$ and that the matrix is
symmetric. Note that obviously any such graph is uniquely defined
by the first row of its adjacency matrix. On the other hand, we know that $w_{j-1}=w_{1,j}=w(v_1,v_j)$ for $j=1,\ldots,n$, namely, $w_{j-1}\neq 0$ iff $j-1 \in \{t_1,\ldots,t_m\}$. From this conditions we
can infer that the symbol has a special polynomial structure and in fact it is equal to
\begin{equation}\label{symbol-toeplitz-graph}
f(\theta)=\sum_{j=1-n}^{n-1} w_{|j|} e^{\i j \theta}=\sum_{k=1}^{m} 2w_{t_k}\cos(t_k  \theta).
\end{equation}

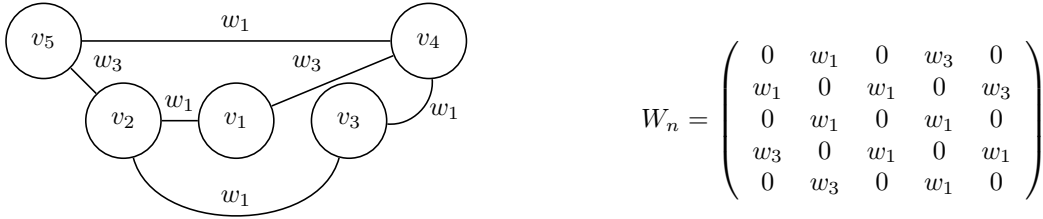
\begin{figure}[H]
	\begin{minipage}{.5\textwidth}
		\centering
		\begin {tikzpicture}[auto ,node distance =1.5 cm,on grid ,
		semithick ,
		state/.style ={ circle ,top color =white , bottom color =white,
			draw,black , text=black , minimum width =1 cm}]
		\node[state] (C){$v_1$};
		\node[state] (D) [left=of C] {$v_2$};
		\node[state] (E) [right=of C] {$v_3$};
		\node[state] (G) [above right=of E] {$v_4$};
		\node[state] (H) [above left=of D] {$v_5$};
		\path (H) edge [black] node[] {$w_3$} (D);
		\path (D) edge [black] node[] {$w_1$} (C);
		\path (C) edge [black] node[] {$w_3$} (G);
		\path (E) edge[black, bend right=50] node[right] {$w_1$} (G);
		\path (G) edge[black] node [above] {$w_1$} (H);
		\path (D) edge[black, bend right=75] node [] {$w_1$} (E);
	\end{tikzpicture}
\end{minipage}
\begin{minipage}{.5\textwidth}
	\begin{equation*}
W_n=\left(
\begin{array}{ccccc}
0      &   w_1  & 0& w_3 &0\\
w_1    &    0   & w_1 & 0& w_3\\
0 & w_1 & 0 & w_1& 0\\
w_3& 0 &   w_1  & 0&w_1\\
0 & w_3 & 0 &w_1 &0\\
\end{array}
\right)
	\end{equation*}
	\end{minipage}
\caption{Example of a $1$-level Toeplitz graph $T_5\langle (1,w_1), (3,w_3)\rangle$. On the left there is a visual representation of the graph while on the right it is explicated the associated adjacency matrix $W_n$ which presents the typical Toeplitz structure. In particular, $W_n$ has symbol function $f(\theta)= 2w_1\cos(\theta) + 2w_3\cos(3\theta)$.}\label{fig:example_toeplitz_graph}
\end{figure}

In such a way, according to (\ref{symbol-toeplitz-graph}), our adjacency matrix $W_n$ is the matrix $T_n(f)$ (real and symmetric) having the following  structure
$$W_n=T_n(f)=\left(
\begin{array}{cccc}
0      &   w_1  & \cdots & w_{n-1}\\
w_1    &    0   & \ddots & \vdots \\
\vdots & \ddots & \ddots & w_1 \\
w_{n-1} & \cdots & w_1     & 0\\
\end{array}
\right)   \qquad w_j = \begin{cases}
w_{t_k} & \mbox{if } j=t_k,\\
0 &\mbox{otherwise},
\end{cases}
$$
and, as expected, the symbol $f$ is real-valued and such that $f(\theta)=f(|\theta|)$ for every $\theta$. See Figure \ref{fig:example_toeplitz_graph} for an example.

Along the same lines, we can define $d$-level Toeplitz graphs as a generalization of the Toeplitz graphs, but beforehand we need to define the set of directions associated to a $d$-index.

\begin{defn}
Given a $d$-index $\bftt_k=((t_k)_1,\ldots,(t_k)_d)$ such that $\bfzero \trianglelefteq \bftt_k$ and $\bftt\neq \bfzero$, then define
$$
I_k := \left\{  \bfii \in \mathbb{Z}^d\, | \,\bfii = \left(\pm (t_k)_1, \ldots, \pm (t_k)_d \right)  \right\},
$$
$$
[\bftt_k] := \quot{I_k}{\sim}, \qquad \mbox{where } \bfii \sim \bfjj \mbox{ iff } \bfii = \pm \bfjj
$$
Trivially, it holds that 
$$
|I_k|= 2^{\sum_{r=1}^{d}\mathds{1}_{(0,\infty)}(|(t_k)_r|)}, \qquad \mbox{where } \mathds{1}_{(0,\infty)}(x)=\begin{cases}
1 &\mbox{if } x \in  (0,\infty)\\
0 & \mbox{otherwise},
\end{cases}
$$
and $\left|[\bftt_k]\right|= |I_k|/2$. We call $[\bftt_k]$ the set of \textnormal{directions} associated to $\bftt_k$. For $\alpha=1,\ldots,\left|[\bftt_k]\right|$, the elements $[\bftt_k]_{\alpha}\in[\bftt_k]$ are called \textnormal{directions} and clearly $\left|[\bftt_k]_\alpha\right|=2$. We will indicate with $[\bftt_k]_\alpha^+$ the element in $[\bftt_k]_\alpha$ such that has positive the first nonzero component and with $[\bftt_k]_\alpha^-$ the other one. 
\end{defn}

\begin{defn}[$d$-level Toeplitz graphs]\label{def:d-toeplitz-graph}
Let $\bfnn,\bftt_1,\ldots,\bftt_m$ be $d$-indices such that $\bf{0}<\bfnn$, let 
$$
0\vartriangleleft\bftt_1\triangleleft \bftt_2 \triangleleft\ldots \triangleleft \bftt_m \triangleleft \bfnn-\bfuno,
$$
and fix $m$ nonzero real vectors $\bfww_1,\ldots, \bfww_m$, such that $\bfww_k \in \R^{s_k}$ with $c_k=\left|[\bftt_k]\right|$ for every $k=1,\ldots,m$, where $[\bftt_k]=\{[\bftt_k]_1,\ldots,[\bftt_k]_{c_k}\}$ is the set of directions associated to $\bftt_k$. We then indicate the components of the vectors $\bfww_k$ by the following index notation,
$$
\bfww_k = \left( w_{[\bftt_k]_1^+}, w_{[\bftt_k]_2^+}, \ldots,  w_{[\bftt_k]_{c_k}^+} \right).
$$
A $d$-\textnormal{level
Toeplitz graph}, denoted by 
$$
T_\bfnn\langle
\{[\bftt_1],\bfww_1\},\ldots,\{[\bftt_m],\bfww_m\}\rangle,
$$ 
is an undirected graph defined by a node set
$V_\bfnn=\left\{v_{\bfkk}\,|\, \bf{1}\trianglelefteq \bfkk\trianglelefteq\bfnn \right\}$ and a weight function $\omega$ such that 
\begin{equation}\label{eq:weight_dToeplitz}
w(v_\bfii,v_\bfjj) = \begin{cases}
w_{[\bftt_k]_\alpha^+} & \mbox{if } |\bfii-\bfjj|=\bftt_k \mbox{ and } (\bfii-\bfjj) \in [\bftt_k]_\alpha=\left\{ [\bftt_k]_\alpha^+, [\bftt_k]_\alpha^-\right\}, \mbox{ for some }\alpha=1,\ldots, c_k,\\
0 & \mbox{otherwise}.
\end{cases}
\end{equation}
If there exist $m$ nonzero real numbers such that $\bfww_k = w_k \bfuno$ for every $k=1,\ldots,m$, then the above relation translates into
	$$
	w(v_\bfii,v_\bfjj) = \begin{cases}
	w_k & \mbox{if } |\bfii-\bfjj|=\bftt_k, \\
	0 & \mbox{otherwise},
	\end{cases}
	$$
and we will indicate the $d$-level Toeplitz graph as $T_\bfnn\langle \{\bftt_1,w_1\},\ldots,\{\bftt_m,w_m\}\rangle$. In the case of simple graph, i.e., $\bfww_k=\bfuno$ for every $k$, then we will indicate the $d$-level Toeplitz graph just as $T_n\langle \bftt_1,\ldots, \bftt_m\rangle$. The number of nodes in a $d$-level Toeplitz graph is equal to $D(\bfnn)$ with
$D(\bfnn)=\prod_{r=1}^d n_r$, while
the number of edges is equal to $\sum_{r=1}^m D(\bfnn-\bftt_r)$.
\end{defn}

\begin{lem}
	A Toeplitz graph is a $1$-level Toeplitz graph as in Definition \ref{def:d-toeplitz-graph}.
\end{lem}
{\bf Proof}
We simply note that, for $d=1$, $\bfnn, \bftt_1,\ldots,\bftt_m$ and the associated $\bfww_1,\ldots, \bfww_m$ are scalars, so that the resulting graph has $n$ points and weight function given by
$$
w(v_i,v_j) = \begin{cases}
w_{t_k} & \mbox{if } |i-j|=t_k,\\
0 & \mbox{otherwise}.
\end{cases}
$$
as in definition \ref{def:toeplitz-graph}, completing the proof.
\hfill \ \, $\bullet$ \ \\

If we assume that $m$, $\{[\bftt_1],\bfww_1\},\ldots,\{[\bftt_m],\bfww_m\}$, are fixed
(independent of $\bfnn$) and we let the sizes $n_j$ grow,
$j=1,\ldots,d$, then  the sequence of adjacency matrices  can be
related to a  unique real integrable function $f: [-\pi,\pi]^d\to \R$ (the
{\itshape symbol}) and expanded
periodically on $\mathbb{R}^d$. In this case, the entries $w_{\bfii,
\bfjj}=\hat{f}_{\bfii-\bfjj}$ of the adjacency matrix are defined via the Fourier
coefficients of $f$, where the $\bfkk$-th Fourier coefficient of
$f$ is defined according to the equations in (\ref{fourier}). Following the same considerations which led to Equation \eqref{symbol-toeplitz-graph}, we can summarize everything said till now in the following proposition.
\begin{prop}\label{prop:d-level-distribution}
	Fix a $d$-level Toeplitz graph $T_\bfnn\langle
	\{[\bftt_1],\bfww_1\},\ldots,\{[\bftt_m],\bfww_m\}\rangle$, and assume that $m$,\\ $\{[\bftt_1],\bfww_1\},\ldots,\{[\bftt_m],\bfww_m\}$ are fixed and
	independent of $\bfnn$. Then the adjacency matrix $W_\bfnn$ of the graph is a symmetric matrix with a $d$-level Toeplitz structure (see Section \ref{sec:toeplitz_matrices}),
	\begin{equation}
	W_\bfnn = \left(w_{\bfii-\bfjj}\right)_{\bfii,\bfjj=\bfuno}^{\bfnn},\qquad \mbox{where }  w_{\bfii-\bfjj}=\begin{cases}
	w_{[\bftt_k]_\alpha} & \mbox{if } |\bfii-\bfjj|=\bftt_k \mbox{ and } (\bfii-\bfjj) \in [\bftt_k]_\alpha, \mbox{ for some }\alpha=1,\ldots, c_k,\\
	0 & \mbox{otherwise}.
	\end{cases}
	\end{equation}
	In particular $W_\bfnn=T_\bfnn(f)$ with symbol function $f: [-\pi,\pi]^d\to \R$ given by
	\begin{equation}\label{eq:symbol_d_toeplitz}
	f(\bftheta)=\sum_{k=1}^{m}\sum_{\alpha=1}^{c_k} 2w_{[\bftt_k]_\alpha}\cos([\bftt_k]_\alpha \cdot \bftheta), \qquad \mbox{with } c_k=\left|[\bftt_k] \right| \mbox{ and } \bftheta=(\theta_1,\ldots,\theta_d),
	\end{equation} 
	that is,
	\begin{equation*}
	\left\{ W_\bfnn \right\}_{\bfnn} \sim_\lambda \bfff.
	\end{equation*}
	
\end{prop}
{\bf Proof}
The fact that
\begin{equation*}
W_\bfnn = \left(w_{\bfii-\bfjj}\right)_{\bfii,\bfjj=\bfuno}^{\bfnn}
\end{equation*}
is clear by definition \ref{def:d-toeplitz-graph}, while, by direct computation of the Fourier coefficients of $f$, we see that $\hat f_{\bfii-\bfjj} = w_{\bfii-\bfjj}$, so that $W_\bfnn=T_\bfnn(f)$.
\hfill \ \, $\bullet$ \ \\

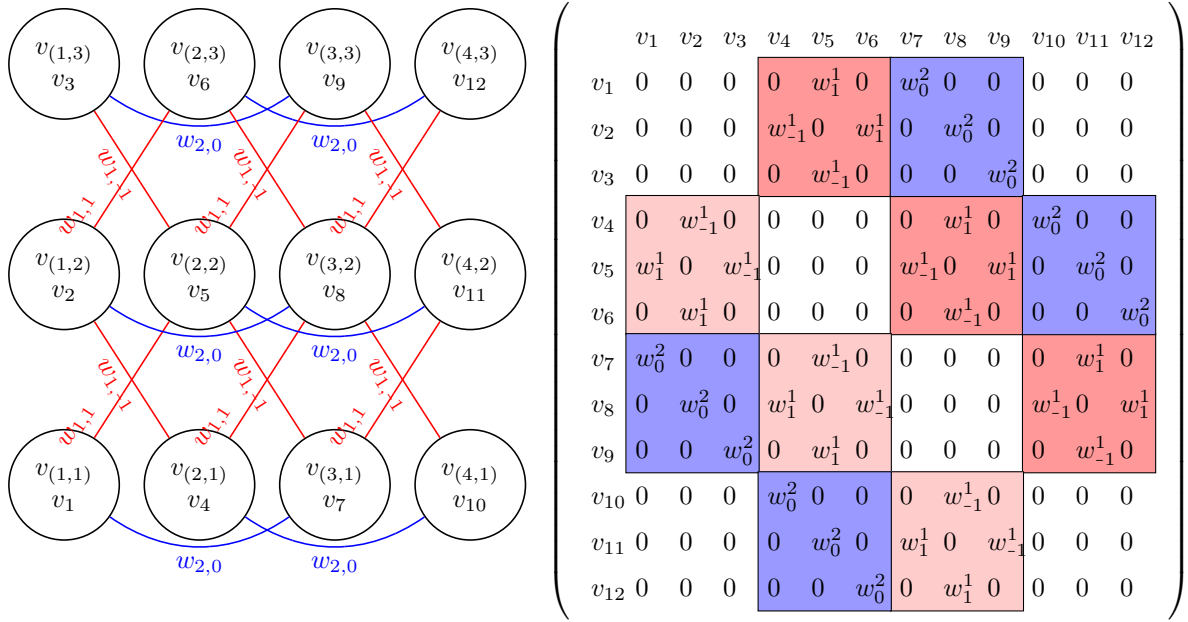
\begin{figure}[H]
	\begin{minipage}{.45\textwidth}
		\centering
		\begin {tikzpicture}[auto ,node distance =2.8 cm and 1.8 cm,on grid ,
		semithick ,
		state/.style ={ circle ,top color =white , bottom color =white,
			draw,black , text=black , minimum width =1 cm}]
		\node[state] (A) {\makecell[t]{$v_{(1,1)}$\\$v_1$}};
		\node[state] (B) [right=of A] {\makecell[t]{$v_{(2,1)}$\\$v_4$}};
		\node[state] (C) [right=of B] {\makecell[t]{$v_{(3,1)}$\\$v_7$}};
		\node[state] (D) [right=of C] {\makecell[t]{$v_{(4,1)}$\\$v_{10}$}};
		\node[state] (E) [above=of A] {\makecell[t]{$v_{(1,2)}$\\$v_2$}};
		\node[state] (F) [above=of B] {\makecell[t]{$v_{(2,2)}$\\$v_5$}};
		\node[state] (G) [above=of C] {\makecell[t]{$v_{(3,2)}$\\$v_8$}};
		\node[state] (H) [above=of D] {\makecell[t]{$v_{(4,2)}$\\$v_{11}$}};
		\node[state] (I) [above=of E] {\makecell[t]{$v_{(1,3)}$\\$v_3$}};
		\node[state] (L) [above=of F] {\makecell[t]{$v_{(2,3)}$\\$v_{6}$}};
		\node[state] (M) [above=of G] {\makecell[t]{$v_{(3,3)}$\\$v_{9}$}};
		\node[state] (N) [above=of H] {\makecell[t]{$v_{(4,3)}$\\$v_{12}$}};
		\path (B) edge [black,bicolor={red}{red}] node[sloped, near start] {$w_{1,1}$} (G);
		\path (C) edge [black,bicolor={red}{red}] node[sloped, near start] {$w_{1,\shortminus 1}$} (F);
		\path (A) edge [black,bicolor={red}{red}] node[sloped, near start] {$w_{1,1}$} (F);
		\path (B) edge [black,bicolor={red}{red}] node[sloped, near start] {$w_{1,\shortminus 1}$} (E);
		\path (C) edge [black,bicolor={red}{red}] node[sloped, near start] {$w_{1,1}$} (H);
		\path (D) edge [black,bicolor={red}{red}]node[sloped, near start] {$w_{1,\shortminus 1}$} (G);
		\path (E) edge [black, bend right=35,bicolor={blue}{blue}] node[sloped,below, midway] {$w_{2,0}$} (G);
		\path (F) edge [black, bend right=35,bicolor={blue}{blue}] node[sloped,below, midway] {$w_{2,0}$} (H);
		\path (A) edge [black, bend right=35,bicolor={blue}{blue}] node[sloped,below, midway] {$w_{2,0}$} (C);
		\path (B) edge [black, bend right=35,bicolor={blue}{blue}] node[sloped,below, midway] {$w_{2,0}$} (D);
		\path (I) edge [black, bend right=35,bicolor={blue}{blue}] node[sloped,below, midway] {$w_{2,0}$} (M);
		\path (L) edge [black, bend right=35,bicolor={blue}{blue}] node[sloped,below, midway] {$w_{2,0}$} (N);
		\path (F) edge [black,bicolor={red}{red}] node[sloped, near start] {$w_{1,\shortminus 1}$} (I);
		\path (E) edge [black,bicolor={red}{red}] node[sloped, near start] {$w_{1,1}$} (L);
		\path (G) edge [black,bicolor={red}{red}] node[sloped, near start] {$w_{1,\shortminus 1}$} (L);
		\path (F) edge [black,bicolor={red}{red}] node[sloped, near start] {$w_{1,1}$} (M);
		\path (H) edge [black,bicolor={red}{red}] node[sloped, near start] {$w_{1,\shortminus 1}$} (M);
		\path (G) edge [black,bicolor={red}{red}] node[sloped, near start] {$w_{1,1}$} (N);
	\end{tikzpicture}
\end{minipage}
\begin{minipage}{.5\textwidth}
	\begin{equation*}
	\begin{tikzpicture}[baseline=-\the\dimexpr\fontdimen22\textfont2\relax]
	\matrix (m)[matrix of math nodes,left delimiter=(,right delimiter=), every node/.style={anchor=base,text depth=.5ex,text height=2ex,text width=1em}]
	{ &$v_1$&$v_2$&$v_3$&$v_4$&$v_5$&$v_6$&$v_7$&$v_8$ & $v_{9}$&$v_{10}$&$v_{11}$&$v_{12}$\\
		$v_1$& 0& 0 & 0 & 0 & $w^{1}_{1}$ &0  &$w^2_0$ &0 & 0 &0 &0 &0\\
		$v_2$&0 & 0 & 0 & $w^1_{\shortminus 1}$ & 0 & $w^1_1$ & 0&$w^2_0$& 0 &0 &0 &0\\
		$v_3$&0 & 0 & 0 & 0 & $w^1_{\shortminus 1}$ & 0 &0 &0& $w^2_0$ &0 &0 &0\\
		$v_4$&0 & $w^1_{\shortminus 1}$ & 0 & 0 & 0 & 0 &0 &$w^1_1$& 0 &$w^2_0$ &0 &0\\
		$v_5$&$w^1_1$ & 0 & $w^1_{\shortminus 1}$ &0  & 0 & 0 &$w^1_{\shortminus 1}$ &0& $w^1_1$ & 0&$w^2_0$ &0\\
		$v_6$&0 & $w^1_1$ & 0 &0 &0& 0 &0 &$w^1_{\shortminus 1}$ &0 &0 &0 &$w^2_0$\\
		$v_7$&$w^2_0$ & 0 & 0 &0 &$w^1_{\shortminus 1}$& 0 &0 &0& 0 &0 &$w^1_1$ &0\\
		$v_8$&0 & $w^2_0$ & 0 &$w^1_1$ &0& $w^1_{\shortminus 1}$ &0 &0& 0 &$w^1_{\shortminus 1}$ &0 &$w^1_1$\\
		$v_9$&0 & 0 & $w^2_0$ &0 &$w^1_1$& 0 &0 &0& 0 &0 &$w^1_{\shortminus 1}$ &0\\
		$v_{10}$&0 & 0 & 0 &$w^2_0$ &0& 0 &0 &$w^1_{\shortminus 1}$& 0 &0 &0 &0\\
		$v_{11}$&0 & 0 & 0 &0 &$w^2_0$&0 &$w^1_1$ &0& $w^1_{\shortminus 1}$ &0 &0 &0\\
		$v_{12}$&0 & 0 & 0 &0 &0& $w^2_0$ &0 &$w^1_1$& 0 &0 &0 &0\\
	};
\begin{pgfonlayer}{myback}
\fhighlight[red!40]{m-2-5}{m-4-7}
\fhighlight[red!40]{m-5-8}{m-7-10}
\fhighlight[red!40]{m-8-11}{m-10-13}
\fhighlight[red!20]{m-5-2}{m-7-4}
\fhighlight[red!20]{m-8-5}{m-10-7}
\fhighlight[red!20]{m-11-8}{m-13-10}
\fhighlight[blue!40]{m-2-8}{m-4-10}
\fhighlight[blue!40]{m-5-11}{m-7-13}
\fhighlight[blue!40]{m-8-2}{m-10-4}
\fhighlight[blue!40]{m-11-5}{m-13-7}
\end{pgfonlayer}
	\end{tikzpicture}
	\end{equation*}
\end{minipage}
\caption{Example of a $2$-level Toeplitz graph $T_{\bfnn}\langle \{[\bftt_1],\bfww_1\}, \{[\bftt_2],\bfww_2\}\rangle$, where $\bfnn=(4,3)$, $[\bftt_1]=[(1,1)]$, $\bfww_1=(w_{1,\shortminus 1}, w_{1,1})$, $[\bftt_2]=[(2,0)]$ and $\bfww_2=w_{2,0}$. In particular, $[\bftt_1]_1=[(1,1)]_1=\{\pm(1,-1)\}$, $[\bftt_1]_2=[(1,1)]_2=\{\pm(1,1)\}$ and $[\bftt_2]_1=[(2,0)] =\{\pm(2,0)\}$. On the left there is a visual representation of the graph while on the right there is the associated adjacency matrix $W_\bfnn$, where we used the standard lexicographic ordering to sort the nodes $\{v_\bfkk\,|\, (1,1)\trianglelefteq (k_1,k_2) \trianglelefteq (4,3)\}$. Comparing the adjacency matrix of this example with $W_\bfnn$ of \eqref{eq:example_2d}, the red blocks correspond to $\bfww_1$ and $\bfww_1^T$, and the blue blocks correspond to $\bfww_2=\bfww_2^T$. All the blocks are Toeplitz matrices, and indeed $W_\bfnn$ is a matrix which possesses a block Toeplitz with Toeplitz blocks (BTTB) structure. In particular, $W_\bfnn$ has symbol function $f(\theta_1,\theta_2)= 2w_{2,0}\cos(2\theta_1) + 2w_{1,1}\cos(\theta_1+\theta_2) + 2w_{1,\shortminus 1}\cos(\theta_1-\theta_2)$: notice that the coefficient of the variable $\theta_1$ refers to the diagonals block while the coefficient of $\theta_2$ refers to the diagonals inside the block. Finally, observe that $W_\bfnn$ is not connected, since the graph can be decomposed into two disjoint subgraphs $G_1$ and $G_2$ having $\{v_1,v_3,v_5,v_{7},v_{9},v_{11}\}$ and $\{v_{2},v_{4},v_{6},v_{8},v_{10},v_{12}\}$ as vertex sets, respectively.}\label{fig:2l-toeplitz}
\end{figure}

As an example, the adjacency matrix of a $2$-level Toeplitz graph has the form
\begin{equation}\label{eq:example_2d}
W_{\bfnn}=T_{\bfnn}(f)=\left(
\begin{array}{cccc}
\bfww_0      &   \bfww_1^T  & \cdots & \bfww_{n_1-1}^T\\
\bfww_1    &    \bfww_0   & \ddots & \vdots \\
\vdots & \ddots & \ddots & \bfww_1^T \\
\bfww_{n_1-1} & \cdots & \bfww_1     & \bfww_0\\
\end{array}
\right), \;
\bfww_{j_1}=\left(
\begin{array}{ccccc}
w_{j_1,0}      & \cdots    &  w_{j_1,\shortminus j_2} & \cdots &w_{j_1,1\shortminus n_2}\\
\vdots   &    \ddots   & & \ddots & \\
w_{j_1,j_2} &  & w_{j_1,0} &   &w_{j_1,\shortminus j_2}\\
\vdots & \ddots &       & \ddots & \\
w_{j_1,n_2\shortminus 1} & \cdots & w_{j_1,j_2}     & \cdots&w_{j_1,0}\\
\end{array}
\right),
\end{equation}
with $j_1=0,\ldots,n_1-1$, $j_2=0,\ldots,n_2-1$ and $\bfww_0=\bfww_0^T$. See Figure \ref{fig:2l-toeplitz} for an explicit graphic example.

\subsection{Graphs with uniform local structure: introducing the  ``diamond''}\label{ssec:local graphs-diamond}

The idea here is that each node in Definition \ref{def:d-toeplitz-graph} is replaced by a subgraph of fixed
dimension $\nu$. For instance,
fix a reference simple graph 
	$$
	G=\left([\nu],E\right)
	$$
	with adjacency matrix $W$ and where $[\nu]$ is the standard set of cardinality $\nu\in \N$. Consider $0<n\in \N$ copies of such a graph, i.e., $G(k)=(V(k),E(k))$ such that $G(k)\simeq G$ for every $k=1,\ldots,n$. Indicating the distinct elements of each $V(k)$, $k=1,\ldots,n$, with the notation $v_{(k,r)}$, for $r=1,\ldots,\nu$, we can define a new node set $V_{n}$ as the disjoint union of the sets $V(k)$, i.e.,
$$
V_{n}:=\bigsqcup_{k=1}^n V(k)=\left\{ v_{(k,r)}\, : \, (1,1)\trianglelefteq (k,r) \trianglelefteq (n,\nu) \right\}.
$$

Fix now $m$ integers $0<t_1<\ldots<t_m$ with $1\leq m\leq n-1$, and moreover fix $L_{t_1},\ldots,L_{t_m}$ simple linking-graph operators for the reference node set $[\nu]$, as in Definition \ref{def:linking-graph}, along with their uniquely determined edge sets $E_{t_1}, \ldots, E_{t_m}\subseteq [\nu]\times [\nu]$. Let us define the edge set $E_{n} \subseteq V_{n}\times V_{n}$,
	$$
	(v_{(i,r)}, v_{(j,s)}) \in E_\bfnn \qquad iff \qquad \begin{cases}
	i=j & \mbox{and } (r,s) \in E, \mbox{ or}\\
	i-j= t_k \mbox{ for some } k=1,\ldots,m & \mbox{and } (r,s) \in E_{t_k},\\
	i-j= -t_k \mbox{ for some } k=1,\ldots,m & \mbox{and } (s,r) \in E_{t_k}.
	\end{cases}
	$$
	Namely, $E_{n}$ is the disjoint union of all the edge sets $E(k)$ plus all the edges which possibly connect nodes in a graph $G(i)$ with nodes in a graph $G(j)$: two graphs $G(i),G(j)$ are connected iff $|i-j| \in \{t_1,\ldots,t_m\}$ and in that case the connection between the nodes of the two graphs is determined by the linking-graph operator $L_{t_k}$ (and by its transpose $L_{t_k}^T$). We can define then a kind of symmetric \textquoteleft weight-graph function \textquoteright
	$$
	\bfww : \{V(k)\, | \, k=1,\ldots,n\}\times \{V(k)\, | \, k=1,\ldots,n\} \to \R^{\nu\times\nu}
	$$
	such that
	$$
	\bfww \left[V(i),V(j)\right]:=\begin{cases}
	W & \mbox{if } i=j,\\
	L_{t_K} & \mbox{if } i-j \in \{t_1,\ldots,t_m\},\\
	L_{t_k}^T & \mbox{if } i-j \in \{-t_1,\ldots,-t_m\},\\
	0 & \mbox{otherwise.}
	\end{cases}
	$$
	It is not difficult then to prove that the adjacency matrix $W_{n,\nu}^G$ of the graph $(V_{n},E_{n})$ is of the form
	$$
	W_{n,\nu}^G= \begin{pmatrix}
	\bfww _0 &\bfww_{1}^T & \cdots & \bfww_{n-1}^T\\
	\bfww _{1} &  \bfww _{0}& \ddots & \vdots\\
	\vdots & \ddots & \ddots & \bfww_{1}^T\\
	\bfww _{n-1} & \cdots & \bfww_{1}& \bfww _{0}
	\end{pmatrix}, \quad \mbox{where } \bfww _{j} = \begin{cases}
	W \in \R^{\nu\times \nu}& \mbox{if } j=0,\\
	L_{t_k} \in \R^{\nu\times \nu} & \mbox{if } j=t_k,\\
	\bfzero & \mbox{otherwise.}
	\end{cases}
	$$
	Trivially, $W_{n,\nu}^G$ is a symmetric matrix with a block-Toeplitz structure and symbol function $\bfff$ given by
	$$
	\bfff(\theta)= W + \sum_{k=1}^m\left(L_{t_k}+L_{t_k}^T\right)\cos(t_k\theta) + \sum_{k=1}^m\left(L_{t_k}-L_{t_k}^T\right)\i \sin(t_k\theta). 
	$$
	Let us observe that $\bfff(\theta)$ is an Hermitian matrix in $\mathbb{C}^{\nu\times\nu}$ for every $\theta \in [0,\pi]$, and therefore $\lambda_j\left(\bfff(\theta)\right)$ are real for every $j=1,\ldots,\nu$, as we requested at the end of Subsection \ref{ssec:weyl_distribution}. We will call $T_{n}^G\left\langle \left(t_1,L_{t_1}\right),\ldots,\left(t_m,L_{t_m}\right) \right\rangle :=(V_{n},E_{n})$ a (simple) {\itshape diamond Toeplitz graph} associated to the graph $G$. A copy $G(k)$ of the graph $G$ will be called $k$-th diamond.

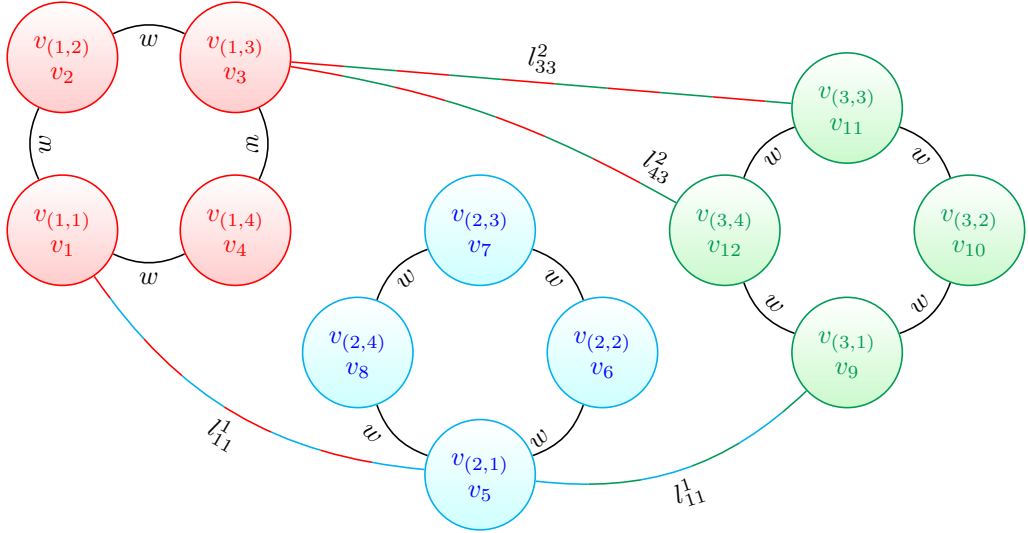
\begin{figure}[H]
\begin{minipage}[tc]{0.4\textwidth}
$W=$\begin{tikzpicture}[baseline=-\the\dimexpr\fontdimen22\textfont2\relax]
\matrix (m)[matrix of math nodes,left delimiter=(,right delimiter=), every node/.style={anchor=base,text depth=.5ex,text height=2ex,text width=1em}]
{ 0 & $w$ & 0 & $w$\\
	$w$ & 0 & $w$ & 0\\
	0& $w$ & 0 & $w$\\
	$w$ & 0 & $w$ & 0\\
};
\begin{pgfonlayer}{myback}
\fhighlight[Brown!20]{m-1-1}{m-4-4}
\end{pgfonlayer}
\end{tikzpicture}\\
$L_1=$\begin{tikzpicture}[baseline=-\the\dimexpr\fontdimen22\textfont2\relax]
\matrix (m)[matrix of math nodes,left delimiter=(,right delimiter=), every node/.style={anchor=base,text depth=.5ex,text height=2ex,text width=1em}]
{ $l_{11}^1$ & 0 & 0 & 0\\
	0 & 0 & 0 & 0\\
	0& 0 & 0 & 0\\
	0 & 0 & 0 & 0\\
};
\begin{pgfonlayer}{myback}
\fhighlight[violet!20]{m-1-1}{m-4-4}
\end{pgfonlayer}
\end{tikzpicture}\\
$L_2=$\begin{tikzpicture}[baseline=-\the\dimexpr\fontdimen22\textfont2\relax]
\matrix (m)[matrix of math nodes,left delimiter=(,right delimiter=), every node/.style={anchor=base,text depth=.5ex,text height=2ex,text width=1em}]
{ 0 & 0 & 0 & 0\\
0 & 0 & 0 & 0\\
	0& 0 & $l_{33}^2$ &0\\
	0 & 0 & $l_{43}^2$& 0\\
};
\begin{pgfonlayer}{myback}
\fhighlight[Goldenrod!20]{m-1-1}{m-4-4}
\end{pgfonlayer}
\end{tikzpicture}
\end{minipage}
\begin{minipage}[t]{0.4\textwidth}
	\centering
\begin{tikzpicture}[baseline=-\the\dimexpr\fontdimen22\textfont2\relax]
\matrix (m)[matrix of math nodes,left delimiter=(,right delimiter=), every node/.style={anchor=base,text depth=.5ex,text height=2ex,text width=1em}]
{        &\textcolor{red}{$v_1$}&\textcolor{red}{$v_2$}&\textcolor{red}{$v_3$}&\textcolor{red}{$v_4$}&\textcolor{blue}{$v_5$}&\textcolor{blue}{$v_6$}&\textcolor{blue}{$v_7$}&\textcolor{blue}{$v_8$} & \textcolor{ForestGreen}{$v_{9}$}&\textcolor{ForestGreen}{$v_{10}$}&\textcolor{ForestGreen}{$v_{11}$}&\textcolor{ForestGreen}{$v_{12}$}\\
	\textcolor{red}{$v_1$}& 0& $w$ & 0 & $w$ & $l^1_{11}$ &0  &0 &0 & 0 &0 &0 &0\\
	\textcolor{red}{$v_2$}&$w$ & 0 & $w$ & 0 & 0 & 0 & 0&0& 0 &0 &0 &0\\
	\textcolor{red}{$v_3$}&0 & $w$ & 0 & $w$ & 0 & 0 &0 &0& 0 &0 &$l^2_{33}$ &$l^2_{43}$\\
	\textcolor{red}{$v_4$}&$w$ & 0 & $w$ & 0 & 0 & 0 &0 &0& 0 &0 &0 &0\\
	\textcolor{blue}{$v_5$}&$l^1_{11}$ & 0 & 0 &0  & 0 & $w$ &0 &$w$ & $l^1_{11}$ & 0&0 &0\\
	\textcolor{blue}{$v_6$}&0 & 0 & 0 &0 &$w$& 0 &$w$ &0 &0 &0 &0 &0\\
	\textcolor{blue}{$v_7$}&0 & 0 & 0 &0 &0& $w$ &0 &$w$& 0 &0 &0 &0\\
	\textcolor{blue}{$v_8$}&0 & 0 & 0 &0 &$w$& 0 &$w$ &0& 0 &0 &0 &0\\
	\textcolor{ForestGreen}{$v_9$}&0 & 0 & 0 &0 &$l^1_{11}$& 0 &0 &0& 0 &$w$ &0 &$w$\\
	\textcolor{ForestGreen}{$v_{10}$}&0 & 0 & 0 &0 &0& 0 &0 &0& $w$ &0 &$w$ &0\\
	\textcolor{ForestGreen}{$v_{11}$}&0 & 0 & $l^2_{33}$ &0 &0&0 &0 &0& 0 &$w$ &0 &$w$\\
	\textcolor{ForestGreen}{$v_{12}$}&0 & 0 & $l^2_{43}$ &0 &0& 0 &0 &0& $w$ &0 &$w$ &0\\
};
\begin{pgfonlayer}{myback}
\fhighlight[Brown!20]{m-2-2}{m-5-5}
\fhighlight[Brown!20]{m-6-6}{m-9-9}
\fhighlight[Brown!20]{m-10-10}{m-13-13}
\fhighlight[violet!20]{m-2-6}{m-5-9}
\fhighlight[violet!20]{m-6-10}{m-9-13}
\fhighlight[violet!20]{m-6-2}{m-9-5}
\fhighlight[violet!20]{m-10-6}{m-13-9}
\fhighlight[Goldenrod!60]{m-2-10}{m-5-13}
\fhighlight[Goldenrod!20]{m-10-2}{m-13-5}
\end{pgfonlayer}
\end{tikzpicture}
\end{minipage}\par\medskip
\centering
\begin {tikzpicture}[auto ,node distance =2.3 cm,on grid ,
semithick ,
state/.style ={ circle ,top color =white , bottom color = processblue!20 ,
	draw,processblue , text=blue , minimum width =1 cm},
state2/.style ={ circle ,top color =white , bottom color = ForestGreen!20 ,
	draw,ForestGreen , text=ForestGreen , minimum width =1 cm},
state3/.style ={ circle ,top color =white , bottom color = red!20 ,
	draw,red , text=red , minimum width =1 cm},
stateW/.style ={ circle ,top color =white , bottom color = white ,
	draw,white , text=white , minimum width =1 cm}]
\node[state] (A)  {\makecell[t]{$v_{(2,1)}$\\$v_5$}};
\node[state] (C) [above left=of A] {\makecell[t]{$v_{(2,4)}$\\$v_8$}};
\node[state] (B) [above right=of A] {\makecell[t]{$v_{(2,2)}$\\$v_6$}};
\node[state] (D) [above right=of C] {\makecell[t]{$v_{(2,3)}$\\$v_{7}$}};

\node[state2] (E) [above right=of B] {\makecell[t]{$v_{(3,4)}$\\$v_{12}$}};
\node[state2] (F) [above right=of E] {\makecell[t]{$v_{(3,3)}$\\$v_{11}$}};
\node[state2] (G) [below right=of F] {\makecell[t]{$v_{(3,2)}$\\$v_{10}$}};
\node[state2] (H) [below left=of G] {\makecell[t]{$v_{(3,1)}$\\$v_{9}$}};

\node[state3] (I) [above left=of C] {\makecell[t]{$v_{(1,4)}$\\$v_{4}$}};
\node[state3] (L) [left=of I] {\makecell[t]{$v_{(1,1)}$\\$v_{1}$}};
\node[state3] (M) [above=of L] {\makecell[t]{$v_{(1,2)}$\\$v_{2}$}};
\node[state3] (N) [right=of M] {\makecell[t]{$v_{(1,3)}$\\$v_{3}$}};

\path (A) edge [black, bend left=25] node[sloped,midway] {$w$} (C);
\path (A) edge [black, bend right=25] node[sloped,midway] {$w$} (B);
\path (B) edge [black, bend right=25] node[sloped,midway,below] {$w$} (D);
\path (C) edge [black, bend left=25] node[sloped,midway,below] {$w$} (D);

\path (E) edge [black, bend left=25] node[sloped,midway,below] {$w$} (F);
\path (F) edge [black, bend left=25] node[sloped,midway,below] {$w$} (G);
\path (G) edge [black, bend left=25] node[sloped,midway,above] {$w$} (H);
\path (H) edge [black, bend left=25] node[sloped,midway,above] {$w$} (E);

\path (I) edge [black, bend left=25] node[sloped,midway] {$w$} (L);
\path (L) edge [black, bend left=25] node[sloped,midway,below] {$w$} (M);
\path (M) edge [black, bend left=25] node[sloped,midway,below] {$w$} (N);
\path (N) edge [black, bend left=25] node[sloped,midway,below] {$w$} (I);

\draw[-,bicolor={processblue}{ForestGreen},bend right=25]
(A) to node[sloped,midway, below] {$\textcolor{black}{l^1_{11}}$} (H);
\draw[-,bicolor={processblue}{red},bend left=25]
(A) to node[sloped,midway, below] {$\textcolor{black}{l^1_{11}}$} (L);

\draw[-,bicolor={red}{ForestGreen}]
(N) to node[sloped, above] {$\textcolor{black}{l^2_{33}}$}  (F);
\draw[-,bicolor={red}{ForestGreen},bend left=10]
(N) to node[sloped,very  near end] {$\textcolor{black}{l^2_{43}}$} (E);
\end{tikzpicture}
\caption{Example of a $1$-level diamond Toeplitz graph $T^G_{3}\left\langle \left(1,L_1\right),\left(2,L_2\right) \right\rangle$, with mold graph $G=T_4\langle (1,w),(3,w) \rangle$. The adjacency matrix of $G$ is $W$. The node sets of the diamond graphs $G(1), G(2), G(3)$ are $V(1)=\{v_{(1,1)},v_{(1,2)}, v_{(1,3)},v_{(1,4)} \}, V(2)=\{v_{(2,1)},v_{(2,2)},v_{(2,3)},v_{(2,4)}\}$ and $V(3)=\{v_{(3,1)},v_{(3,2)},v_{(3,3)},v_{(3,4)}\}$, respectively. Clearly, all the diamond graphs are characterized by the same adjacency matrix $W$. The adjacency matrix $W_n$ of the whole graph is a $1$-level block Toeplitz and has symbol function $\bfff(\theta)= W + 2L_1\cos(\theta) + \left(L_2+L_2^T\right)\cos(2\theta) + \left(L_2-L_2^T\right)\i \sin(2\theta)$.}\label{fig:diamond_graph}
\end{figure}

See Figure \ref{fig:diamond_graph} for an example. We can now generalize everything said till now.

\begin{defn}[$d$-level diamond Toeplitz graph]\label{def:d-level-diamond}
	Let $d,m,\nu$ be fixed integers and let $G\simeq \left([\nu], E,w\right)$ be a fixed undirected graph which we call \textnormal{mold graph}. 
	
	Let $\bfnn,\bftt_1,\ldots,\bftt_m$ be $d$-indices such that $\bf{0}<\bfnn$, and $0\vartriangleleft\bftt_1\triangleleft \bftt_2 \triangleleft\ldots \triangleleft \bftt_m \triangleleft \bfnn-\bfuno.$ For $k=1,\ldots,m$, let $\boldsymbol{L}_k$ be a collection of linking-graph operators of the standard set $[\nu]$ and such that $\left|\boldsymbol{L}_k \right|=c_k$, with $c_k=\left|[\bftt_k]\right|$ for every $k=1,\ldots,m$, where $[\bftt_k]=\{[\bftt_k]_1,\ldots,[\bftt_k]_{c_k}\}$ is the set of directions associated to $\bftt_k$. We then indicate the elements of the sets $\boldsymbol{L}_k$ by the following index notation,
	$$
	\boldsymbol{L}_k = \left\{ L_{[\bftt_k]_1^+}, L_{[\bftt_k]_2^+}, \ldots,  L_{[\bftt_k]_{c_k}^+} \right\}, \qquad \R^{\nu\times\nu}\ni L_{[\bftt_k]_\alpha^+}=\left(l_{[\bftt_k]_\alpha^+}(r,s)\right)_{r,s=1}^\nu \quad \mbox{for } \alpha=1,\ldots, c_k.
	$$
Finally, consider $\bfnn$ copies $G(\bfkk)\simeq G$ of the mold graph, which we will call \textnormal{diamonds}.
	
	A $d$-\textnormal{level diamond
		Toeplitz graph}, denoted by $T_{\bfnn,\nu}^G\left\langle \left\{\bftt_1,\boldsymbol{L}_1\right\},\ldots,\left\{\bftt_m,\boldsymbol{L}_m\right\} \right\rangle$, is an undirected graph with
	$$
	V_\bfnn=\left\{v_{(\bfkk,r)}\,|\, (\bf{1},1)\trianglelefteq ({\bfkk}, r)\trianglelefteq (\bfnn,\nu) \right\}
	$$
	and characterized by the weight function $w_{\bfnn}: V_\bfnn \times V_\bfnn\to \R$ such that
	$$
	w_{\bfnn}\left(v_{(\bfii,r)},v_{(\bfjj,s)}\right) := \begin{cases}
	w(r,s) & \mbox{if }\bfii=\bfjj,\\
	l_{[\bftt_k]_\alpha^+}(r,s) & \mbox{if } |\bfii-\bfjj|=\bftt_k \mbox{ and } (\bfii-\bfjj) = [\bftt_k]_\alpha^+, \mbox{ for some }\alpha=1,\ldots, c_k,\\
	l_{[\bftt_k]_\alpha^+}(s,r) & \mbox{if } |\bfii-\bfjj|=\bftt_k \mbox{ and } (\bfii-\bfjj) = [\bftt_k]_\alpha^-,\mbox{ for some }\alpha=1,\ldots, c_k,\\
	0 & \mbox{otherwise}.
	\end{cases}
	$$
The number of nodes in a $d$-level diamond Toeplitz graph is equal to $\nu D(\bfnn)$ with
	$D((\bfnn)=\prod_{r=1}^d n_r$, while
	the number of edges is equal to $\nu\sum_{r=1}^m D(\bfnn-\bftt_r)$. 
\end{defn}

\begin{cor}
A $d$-level Toeplitz graph is a special case of a $d$-level diamond Toeplitz graph.
\end{cor}
{\bf Proof}
We simply need to notice that, for $\nu=1$, i.e. in the case of a diamond with only one element, the two definitions \ref{def:d-toeplitz-graph} and \ref{def:d-level-diamond} coincide with $\boldsymbol{L}_k=w_k$.
\hfill \ \, $\bullet$ \ \\

\begin{prop}\label{prop:diamond-distribution}
Fix a $d$-level diamond
	Toeplitz graph $T_{\bfnn,\nu}^G\left\langle \left\{\bftt_1,\boldsymbol{L}_1\right\},\ldots,\left\{\bftt_m,\boldsymbol{L}_m\right\} \right\rangle$ with $G\simeq \left([\nu], E,w\right)$ and $W$ the adjacency matrix of $G$. Let $d,m,\nu,\{\bftt_k,\boldsymbol{L}_k\},G$ be fixed and independent of $\bfnn$. Then the adjacency matrix $W_{\bfnn,\nu}^G$ of  $T_{\bfnn,\nu}^G\left\langle \left\{\bftt_1,\boldsymbol{L}_1\right\},\ldots,\left\{\bftt_m,\boldsymbol{L}_m\right\} \right\rangle$ is a symmetric matrix with a $d$-level block Toeplitz structure (see Section \ref{sec:toeplitz_matrices} and Equation \eqref{toeplitz-symbol}),
	\begin{equation}
	W^G_{\bfnn,\nu} = \left[\bfww_{\bfii-\bfjj}\right]_{\bfii,\bfjj=\bfuno}^{\bfnn},\qquad \mbox{where } \R^{\nu\times\nu}\ni \bfww_{\bfii-\bfjj}=\begin{cases}
	W & \mbox{if } \bfii = \bfjj,\\
	L_{[\bftt_k]_\alpha^+} & \mbox{if } |\bfii-\bfjj| = \bftt_k \mbox{ and } (\bfii-\bfjj)= [\bftt_k]_\alpha^+,\\
	L_{[\bftt_k]_\alpha^+}^T & \mbox{if } |\bfii-\bfjj| = \bftt_k \mbox{ and } (\bfii-\bfjj)= [\bftt_k]_\alpha^-,\\
	\bfzero & \mbox{otherwise}.
	\end{cases}
	\end{equation}
In particular $W^G_{\bfnn,\nu}=T_{\bfnn,\nu}(\bfff)$ with symbol function $\bfff : [-\pi,\pi]^d \to \mathbb{C}^{\nu\times \nu}$ given by
\begin{equation}\label{eq:d-diamond-symbol}
\bfff(\bftheta)=W + \sum_{k=1}^{m}\sum_{\alpha=1}^{c_k}\left[ \left(L_{[\bftt_k]_\alpha^+} + L_{[\bftt_k]_\alpha^+}^T\right)\cos(\bftt_k \cdot \bftheta) +  \left(L_{[\bftt_k]_\alpha^+} - L_{[\bftt_k]_\alpha^+}^T\right)\i\sin(\bftt_k \cdot \bftheta)\right], \qquad \boldsymbol{\theta}=(\theta_1,\ldots,\theta_d),
\end{equation} 
that is,
\begin{equation*}
\left\{ W^G_{\bfnn,\nu} \right\}_{\bfnn} \sim_\lambda \bfff.
\end{equation*}
The symbol function $\bfff$ is Hermitian-matrix valued for every $\bftheta\in [0,\pi]^d$.
\end{prop}
{\bf Proof}
 We note that $W^G_\bfnn = \left[\bfww_{\bfii-\bfjj}\right]_{\bfii,\bfjj=\bfuno}^{\bfnn}$ is immediate by definition \ref{def:d-level-diamond} and that the symbol $\bfff$ is a Hermitian matrix for every $\bftheta$, so that it has real eigenvalues. Moreover we see that, as in Proposition \ref{prop:d-level-distribution}, $\hat f_{\bfii-\bfjj} = w_{\bfii-\bfjj}$. Now Theorem \ref{thm:symbol_d-block-toeplitz} concludes the proof.
\hfill \ \, $\bullet$ \ \\

\section{Grid graphs with uniform local structure and main spectral results}\label{sec:main}

The section is divided into two parts. In the first we give the definition of grid graphs with uniform
local structure. In the second part we show the links of the above notions with Toeplitz and GLT sequences and we use the latter for proving the main spectral results.

\subsection{Sequence of grid graphs with uniform local structure}\label{ssec:local graphs}
The main idea in this section is to immerse the graphs presented in Section \ref{sec:toeplitz_matrices} inside a bounded regular domain $\Omega\subset \R^d$. We start with a series of definitions in order to give a
mathematical rigor to our derivations.

\begin{defn}[$d$-level Toeplitz grid graphs in the cube]\label{local gr cube} Given a continuous almost everywhere (a.e.) function $p: [0,1]^d \to \R$, choose a $d$-level Toeplitz graph 
	$$
	T_{\bfnn}\langle \{[\bftt_1],\bfww_1\},\ldots,\{[\bftt_m],\bfww_m\}\rangle,
$$
 and consider the $d$-dimensional vector
	
		$$
		\bfhh:=(h_1,\ldots,h_d)=\left(\frac{1}{n_1+1},\ldots, \frac{1}{n_d+1}\right).
		$$
	We introduce a bijective correspondence between the nodes $v_\bfjj$ of $T_{\bfnn}\langle \{[\bftt_1],\bfww_1\},\ldots,\{[\bftt_m],\bfww_m\}\rangle$ and the interior points $\bfxx$ of the cube $[0,1]^d$ by the immersion map $\iota : V_\bfnn \to (0,1)^d$ such that 
	$$
	\iota(v_\bfjj):= \bfjj \circ \bfhh = \left(j_1h_1,\ldots, j_dh_d \right)
	$$
	with $\circ$ being the Hadamard (component-wise) product. The $d$-level Toeplitz graph induces then a grid graph in $[0,1]^d$, $G=\left(\tilde{V}_\bfnn,\tilde{E}_\bfnn, w^p\right)$ with 
		$$
		\tilde{V}_\bfnn:=\left\{ \bfxx_\bfkk = \iota(v_\bfkk)\,|\, \bf{1}\trianglelefteq \bfkk \trianglelefteq\bfnn  \right\}, \qquad \tilde{E}_\bfnn :=\left\{ (\bfxx_\bfii,\bfxx_\bfjj) \, | \, w^p(\bfxx_\bfii,\bfxx_\bfjj)\neq 0  \right\},
		$$
	where 
$$
w^p (\bfxx_\bfii,\bfxx_\bfjj) := 
p\left(\frac{\bfxx_\bfii + \bfxx_\bfjj}{2}\right)w(v_\bfii,v_\bfjj),
$$
and $w$ is the weight function defined in \eqref{eq:weight_dToeplitz}. 
	With abuse of notation we will identify $V_\bfnn = \tilde{V}_\bfnn$ and we will write 
		$$
		T_{\bfnn}\langle \{[\bftt_1],\bfww^p_1\},\ldots,\{[\bftt_m],\bfww^p_m\}\rangle, 
		$$ 
		for a $d$-level grid graph in $[0,1]^d$.
\end{defn}
Observe that now $\bfww^p_k$, for $k=1,\ldots,m$, are not anymore constant vectors as $\bfww_k$ but vector-valued functions $\bfww^p_k: [0,1]^d\times [0,1]^d \to \R^{c_k}$, with $c_k=|[\bftt_k]|$, such that
	$$
(\bfww^p_k)_\alpha(\bfxx_\bfii,\bfxx_\bfjj) = \begin{cases}
p\left(\frac{\bfxx_\bfii + \bfxx_\bfjj}{2}\right)w_{[\bftt_k]_\alpha^+} & \mbox{if } |\bfii-\bfjj|=\bftt_k \mbox{ and } (\bfii-\bfjj) \in [\bftt_k]_\alpha=\left\{ [\bftt_k]_\alpha^+, [\bftt_k]_\alpha^-\right\},\\
0 & \mbox{otherwise},
\end{cases}  
$$ 
for $\alpha=1,\ldots,c_k$. It is then not difficult to see  that we can express the weight function $w^p$ as 
$$
w^p(\bfxx_\bfii,\bfxx_\bfjj) = \sum_{k=1}^m\sum_{\alpha=1}^{c_k} (\bfww^p_k)_\alpha(\bfxx_\bfii,\bfxx_\bfjj). 
$$
In other words, taking in mind the role of the reference domain
$[0,1]^d$, $\bfxx_\bfii$ can be connected to $\bfxx_\bfjj$ only if
$|(\bfxx_\bfjj)_r-(\bfxx_\bfii)_r|=O(h_r)$, for all $r=1,\ldots,d$.
From this property we derive the name of 'grid graphs with local
structure'. Naturally, the above notion can be generalized to any
domain $\Omega\subset [0,1]^d$: as we will see in the next
subsection, the only restriction in order to have meaningful
spectral properties of the related sequences, is that $\Omega$ is
regular.

\begin{defn}[$d$-level Toeplitz grid graphs in $\Omega$]\label{def:local gr Omega}
Given a regular domain $\Omega \subseteq [0,1]^d$ and a continuous a.e. function $p: [0,1]^d \to \R$, choose a $d$-level Toeplitz graph 
$$
T_{\bfnn}\langle \{[\bftt_1],\bfww_1\},\ldots,\{[\bftt_m],\bfww_m\}\rangle,
$$
and consider its associated $d$-level Toeplitz grid graph $	T_{\bfnn}\langle \{[\bftt_1],\bfww^p_1\},\ldots,\{[\bftt_m],\bfww^p_m\}\rangle$. We define then the $d$-level Toeplitz grid graph immersed in $\Omega$ as the graph $G=(V_{n'}^\Omega, E_{n'}^\Omega, w^{\Omega,p})$  such that
$$
V_{n'}^\Omega:= V_\bfnn \cap \Omega, \qquad w^{\Omega,p}:= w^p_{|V_\bfnn^\Omega\times V_\bfnn^\Omega}.  
$$
Clearly, $\left|V_{n'}^\Omega\right|=n'\leq \prod_{r=1}^d n_r= \left| V_\bfnn\right|$. Nevertheless, $n'=n'(\bfnn) \to \infty$ as $\bfnn \to \infty$. Therefore, with abuse of notation we will keep writing $\bfnn$ instead of $n'$. We will indicate such a graph with the notation
$$
T_{\bfnn}^{\Omega}\langle \{[\bftt_1],\bfww^p_1\},\ldots,\{[\bftt_m],\bfww^p_m\}\rangle.
$$
\end{defn}
\begin{rem}
In the application, as we will see in Section \ref{sec:appl}, once it is chosen the domain $\Omega$ and the kind of discretization technique to solve a differential equation, then the weight function $w$ is fixed accordingly, and consequently the coefficients $\bfww_1,\ldots,\bfww_m$. In particular, it is important to remark that the weight function of $T_{\bfnn}\langle \{[\bftt_1],\bfww_1\},\ldots,\{[\bftt_m],\bfww_m\}\rangle,$ will depend on the differential equation and on the discretization technique.
\end{rem}

Finally, we immerse the diamond graphs in the cube $[0,1]^d$ (and then in a generic regular domain $\Omega \subset [0.1]^d$).
   
\begin{defn}[$d$-level diamond Toeplitz grid graphs in the cube]\label{def:diamind-Toeplitz-grid} 
	The same definition as in Definition \ref{local gr cube} where the $d$-level Toeplitz graph is replaced by a $d$-level diamond Toeplitz graph. The only difference now is that  	
	$$
	\bfhh:=(h_1,\ldots,h_d)=\left(\frac{1}{\nu n_1+1},\frac{1}{ n_2+1},\ldots, \frac{1}{ n_d+1}\right),
	$$
	and
\begin{equation*}
	\iota(v_{(\bfjj,r)}):= (\bfjj,r) \circ \bfhh = \left((j_1+r-1)h_1,j_2h_2,\ldots,j_dh_d \right), \qquad r=1,\ldots,\nu.
\end{equation*}
With abuse of notation we will write 
$$
T_{\bfnn,\nu}^G\langle \{[\bftt_1],\boldsymbol{L}^p_1\},\ldots,\{[\bftt_m],\boldsymbol{L}^p_m\}\rangle, 
$$ 
for a $d$-level diamond Toeplitz grid graph in $[0,1]^d$.
\end{defn}
\begin{rem}
While in the case of a $d$-level Toeplitz graph the immersion map $\iota$ was introduced naturally as the Hadamard product between the indices of the graph nodes and the natural Cartesian representation of points in $\R^d$, the diamond Toeplitz graphs grant another degree of freedom for the immersion map. In Definition \ref{def:diamind-Toeplitz-grid} we decided for the simplest choice, namely lining-up all the nodes of the diamonds along the first axis. Clearly, other choices of the immersion map $\iota$ would be able to describe more complex grid geometries.  
\end{rem}

\begin{defn}[$d$-level  diamond Toeplitz grid graphs in $\Omega$]\label{def:local_diamond_gr_Omega}
The same definition as in Definition \ref{def:local gr Omega} where the $d$-level Toeplitz grid graph is replaced by a $d$-level diamond Toeplitz grid graph. We will indicate such a graph with the notation
$$
T_{\bfnn,\nu}^{G,\Omega}\langle \{[\bftt_1],\boldsymbol{L}^p_1\},\ldots,\{[\bftt_m],\boldsymbol{L}^p_m\}\rangle.
$$
\end{defn}

\subsection{Asymptotic spectral results}\label{ssec:spectral results}

We start this section containing the spectral results, by giving
the distribution theorem in the Weyl sense in its maximal
generality, i.e. for a sequence of weighted (diamond) local grid
graphs in $\Omega$, according to the case depicted in Definition
\ref{def:local_diamond_gr_Omega}.

\begin{thm}\label{thm:d-grid-graphs-symbol}
	Given a regular domain $\Omega\subseteq [0,1]^d$ and a continuous a.e. function $p:\Omega \to \R$, fix a $d$-level Toeplitz grid graph $T^\Omega_{\bfnn}\langle
	\{[\bftt_1],\bfww_1^p\},\ldots,\{[\bftt_m],\bfww_m^p\}\rangle$ as in Definition \ref{def:local gr Omega}, and assume that $m$, $\{[\bftt_1],\bfww_1\},\ldots,\{[\bftt_m],\bfww_m\}$ are fixed and
	independent of $\bfnn$. Then, indicating with $\{W_\bfnn^{\Omega,p}\}_\bfnn$ the sequence of adjacency matrix of the $d$-level Toeplitz grid graph as $\bfnn\to \infty$, it holds that
	\begin{equation}\label{eq:d-grid-graph-symbol}
	\left\{ W_\bfnn^{\Omega,p} \right\}_{\bfnn} \sim_\lambda g,
	\end{equation}
	with $g : \Omega\times [-\pi,\pi]^d\subset \R^{2d} \to \R$ and
	$$
	g(\bfxx,\bftheta)= p(\bfxx)f(\bftheta),
	$$
	where $f(\bftheta)$ is the symbol function defined in \eqref{eq:symbol_d_toeplitz}. 
\end{thm}
{\bf Proof}
We note that, in the case where $\nu=1$, a $d$-level diamond Toeplitz grid graph reduces to a $d$-level Toeplitz grid graph according to definition \ref{def:local gr Omega}. The conclusion of the theorem is then obvious once we prove our next result, Theorem \ref{main-theorem}.
\hfill \ \, $\bullet$ \ \\

\begin{thm}\label{main-theorem}
		Given a regular domain $\Omega\subseteq [0,1]^d$ and a continuous a.e. function $p:\Omega \to \R$, fix a $d$-level diamond Toeplitz grid graph $T^{G,\Omega}_{\bfnn,\nu}\langle
	\{[\bftt_1],\boldsymbol{L}_1^p\},\ldots,\{[\bftt_m],\boldsymbol{L}_m^p\}\rangle$ as in Definition \ref{def:local gr Omega}, and assume that $m$, $\{[\bftt_1],\boldsymbol{L}_1\},\ldots,\{[\bftt_m],\boldsymbol{L}_m\}$ are fixed and
	independent of $\bfnn$. Then, indicating with $\{W_{\bfnn,\nu}^{G,\Omega,p}\}_\bfnn$ the sequence of adjacency matrix of the $d$-level diamond Toeplitz grid graph as $\bfnn\to \infty$, it holds that
	\begin{equation}\label{eq:diamond-grid-graph-symbol}
	\left\{ W_{\bfnn,\nu}^{G,\Omega,p} \right\}_{\bfnn} \sim_\lambda \boldsymbol{g},
	\end{equation}
	with $\boldsymbol{g} : \Omega\times [0,\pi]^d \subset \R^{2d}\to \mathbb{C}^{\nu\times\nu}$ a matrix-valued function  and
	\begin{equation}\label{eq:diamond-grid-graph-symbolB}
	\boldsymbol{g}(\bfxx,\bftheta)= p(\bfxx)\bfff(\bftheta),
	\end{equation}
	where $\bfff(\bftheta)$ is the symbol function defined in \eqref{eq:d-diamond-symbol}. 
\end{thm}
{\bf Proof}
First of all we observe that our assumption of $\Omega$ regular is equivalent to require $\Omega$ to be measurable according to Peano-Jordan, which is the fundamental assumption to apply the GLT theory (see \cite{glt-book-1}).\\
Assume $\Omega=[0,1]^d$ and $p(x)\equiv 1$ over $\Omega$. Then our sequence of graphs reduces to a sequence of $d$-level diamond Toeplitz graphs and the proof is over using Proposition \ref{prop:diamond-distribution}.\\
Assuming now that $\Omega=[0,1]^d$ and $p$ is just a Riemann-integrable function over $\Omega$, we decompose the adjacency matrix $ W_{\bfnn,\nu}^{G,p}$ as $ W_{\bfnn,\nu}^{G,p}= \diag_\bfnn(p)T_\bfnn(\bfff)+E_\bfnn$. The only observation needed here is that $\left\{\diag_\bfnn(p)\right\}$ is a multilevel block GLT with symbol function $p$, while $\left\{T_\bfnn(\bfff)\right\}$ is a multilevel block GLT with symbol function $\bfff$ (see item $\bf GLT 2$ in \cite{GMS18}). Moreover, by direct calculation, we see that $E_\bfnn$, for $\bfnn$ large, can be written as a term of small spectral norm, plus a term of relatively small rank. Therefore, $E_\bfnn$ is a zero-distributed sequence of matrices and hence a multilevel block GLT with symbol function $0$. Summing up we have
$$
\left\{ T_\bfnn(\bfff) \right\}_{\bfnn} \sim \bfff \text{ over } [0,\pi]^d, \left\{\rm Diag_\bfnn(p) \right\}_{\bfnn} \sim p \text{ over } \Omega, \left\{ E_\bfnn \right\}_{\bfnn} \sim 0.
$$
Now, by the structure of algebra of multilevel block GLT sequences and using the symmetry of the sequence (see item $\bf GLT 3$ in \cite{GMS18}), we conclude that $\left\{  W_{\bfnn,\nu}^{G,p} \right\}_{\bfnn} \sim p(x)\bfff(\bftheta)$ over $[0,1]^d \times [0,\pi]^d$.\\
For the general case where $\Omega$ is a generic regular subset of $[0,1]^d$, we simply notice that, using Definition \ref{def:local_diamond_gr_Omega}, we can see $ W_{\bfnn,\nu}^{G,p}$ as a principal sub-matrix of $W_{\bfnn,\nu}^{G,\Omega,p}$, where $\left\{  W_{\bfnn,\nu}^{G,\Omega,p} \right\}$ is constructed according to Definition \ref{def:diamind-Toeplitz-grid} and the function $p$ is substituted by $p_{|\Omega}:=p(i(x))$, where $i : \Omega \hookrightarrow [0,1]^d$ is the inclusion map. Since $\Omega$ is regular we conclude that $p_{|\Omega}$ is Riemann-integrable, $p_{|\Omega}\equiv p$ over $\Omega$, and that 
$$
\left\{  W_{\bfnn,\nu}^{G,\Omega,p} \right\}_{\bfnn} \sim p(x) \bfff(\bftheta),
$$
with $p(x)\bfff(\bftheta)$  defined over $\Omega \times [0,\pi]^d$, as required to complete the proof.
\hfill \ \, $\bullet$ \ \\

\begin{cor}\label{cor:graph-Laplacian-distribution}
Let $G_{\bfnn,\nu}= \left(T_{\bfnn,\nu}^{G,\Omega}\left\langle\left([\bftt_1],\boldsymbol{L}_{1}^p\right),\ldots,\left([\bftt_m],\boldsymbol{L}_{m}^p\right)  \right\rangle, \kappa_n  \right)$ be a $d$-level diamond Toeplitz grid graph as in Definition \ref{def:local_diamond_gr_Omega}, plus a potential term $\kappa_n : V_\bfnn \to [0,\infty)$.  Let $D,K$ be as in Definition \ref{def:graph-Laplacian}. If 
\begin{equation*}
D+K= cI\cdot \diag_{(\bfuno,1)\trianglelefteq (\bfii,r) \trianglelefteq (\bfnn,\nu)}\{p(\bfxx_{(\bfii,r)})\} + o(\bfuno),
\end{equation*}
where $c \in \mathbb{R}$ is a fixed constant, $I$ is the identity matrix and $p : \Omega \subset [0,1]^d \to \R$ is a continuous a.e. function as in Definition \ref{def:local_diamond_gr_Omega}. Then it holds that
\begin{equation*}
\left\{ \Delta_{G_{\bfnn,\nu}} \right\} \sim_\lambda p(\bfxx) \left( cI - \bfff(\bftheta) \right), \qquad (\bfxx,\bftheta) \in \Omega \times [-\pi,\pi]^d,
\end{equation*} 
where $\Delta_{G_{\bfnn,\nu}}$ is the graph-Laplacian as in Definition \ref{def:graph-Laplacian} and $\bfff(\bftheta)$ is defined in \eqref{eq:d-diamond-symbol}.
\end{cor}
{\bf Proof}\, \  From Definition  \ref{def:graph-Laplacian}, $\Delta_{G_{\bfnn,\nu}}= D+K-W_{\bfnn,\nu}^{G,\Omega,p}$. By assumption, note that $ D+K $ is a GLT sequence with symbol function $ cp $, so that the conclusion follows once again by the algebra structure of the GLT sequences.
\hfill \ \, $\bullet$ \ \\ 

\section{Spectral gaps}\label{sec:appl-2}
In this section we report general results on the gaps between the extreme eigenvalues of adjacency and graph-Laplacian matrices of the type considered so far.

The spectral properties of a Toeplitz matrix $T_n(f)$ are well understood by considering $f$; in fact, it is known (see e.g. \cite{GS,Kac1953}) that the spectrum
of $T_n(f)$ is contained in $(m_f,M_f)$, where $m_f =\min f$ and $M_f=\max f$, and moreover
$$\lim_{n\rightarrow\infty}\lambda_1^{(n)}=m_f$$
and
$$\lim_{n\rightarrow\infty}\lambda_n^{(n)}=M_f.$$

Recall the following result due to Kac, Murdoch, Szeg\"o, Parter, Widom and later third author, B\"ottcher, Grudsky, Garoni (see also \cite{BS} and references therein for more details and for the history of such results).
\begin{thm}
	Let $f\in C_{2\pi}[-\pi,\pi]$ be a continuous function on $[-\pi,\pi]$ extended periodically on $\mathbb{R}$, let $m_f=\min f$ and $M_f=\max f$.
	Let $T_n(f)$ the associated Toeplitz matrix and let
	$$
	\lambda_1^{(n)}\leq\lambda_2^{(n)}\leq\ldots\leq\lambda_{n-1}^{(n)}\leq \lambda_n^{(n)}
	$$
	be the eigenvalues of $T_n(f)$ ordered in non-decreasing order.
	
	Then, for all fixed $j$ we have that
	$$\lim_{n\rightarrow\infty}\lambda_j^{(n)}=m_f \qquad \lim_{n\rightarrow\infty}\lambda_{n+1-j}^{(n)}=M_f$$
	and moreover
	$$\lim_{n\rightarrow\infty}\frac{M_f-\lambda_{n+1-j}^{(n)}}{c(j^{-1}n)^\alpha}=1,\qquad c\in\mathbb{R}$$
	where $\alpha\in\mathbb{R}$ is such that, if $f(x_0)=M_f$ then $|f(x_0)-f(x)|\sim c|x-x_0|^\alpha$ as $x\to x_0$.
\end{thm}

At the light of some of the new results presented in \cite{Bianchi} and reported in the Appendix \ref{appendix}, we can make the following statements within the scope of this paper.

\begin{thm}\label{thm:diamond-spectral-gap}
	 Let $\{W_{\bfnn,\nu}^{G,\Omega,p}\}_\bfnn$ be a sequence of adjacency matrices of $d$-level diamond Toeplitz grid graphs as in Theorem \ref{main-theorem}, with $\left\{\lambda_k^{(\bfnn)}\right\}_{k=1}^{d_\bfnn}$ the collection of their eigenvalues sorted in non-decreasing order and where $d_\bfnn = \nu\prod_{j=1}^d n_j$ is the dimension of the matrices. Let $p(\bfxx)$ be piecewise continuous and let $\tilde{g}:[0,1]\to [\min R_\bfgg, \max R_\bfgg]$ be the monotone rearrangement of the symbol function $\bfgg$ of $\{W_{\bfnn,\nu}^{G,\Omega,p}\}$, as in Definition \ref{def:rearrangment}. Finally, let $\tau:[\min R_\bfgg, \max R_\bfgg] \to \R$ be a function that is differentiable in $\max R_\bfgg$ and define $x_\bfnn:= 1-\frac{1}{d_\bfnn}$. If
	 \begin{enumerate}[(i)]
	 	\item $\lambda_{d_\bfnn\shortminus 1}^{(\bfnn)}<\lambda_{d_\bfnn}^{(\bfnn)} \leq \max R_\bfgg$ definitely for $\bfnn \to \infty$;\label{hp1}
	 	\item $\tilde{g}$ is piecewise Lipschitz and differentiable at $x=1$;\label{hp2}
	 \end{enumerate}
 then it holds that  
	\begin{equation*}
	\lim_{\bfnn \to \infty} \frac{d_\bfnn\left[\tau\left(\lambda_{d_\bfnn}^{(\bfnn)}\right) - \tau\left(\lambda_{d_\bfnn-1}^{(\bfnn)} \right)\right]}{\tau'(\tilde{g}(x_\bfnn))\tilde{g}'(x_\bfnn)}=1.
	\end{equation*}
\end{thm}
{\bf Proof}
It is immediate from Corollary \ref{cor:differentiability}.
\hfill \ \, $\bullet$ \ \\

\begin{cor}\label{cor:diamond-spectral-gap}
In the case $p(\bfxx)\equiv c$ and $d=1$ then $\tilde{g}$ is differentiable at $x=1$, and 
\begin{equation}\label{eq:convergence_zero1}
\tilde{g}'(x_n)\to 0,
\end{equation}
for any choice of $\{(t_1,L_{t_1}),\ldots, (t_m,L_{t_m}) \}$ (or $\{(t_1,w_{t_1}),\ldots, (t_m,w_{t_m}) \}$ in the case $\nu=1$). In particular, if $\tau\equiv id$ then
\begin{equation}\label{eq:convergence_zero2}
\nu n\left[\lambda_{\nu n}^{(n)} - \lambda_{\nu n-1}^{(n)}\right]\to 0 \qquad \mbox{as }n\to \infty.
\end{equation}
\end{cor}
{\bf Proof}
Let us fix $c=1$. Since all the components $f_{i,j}(\theta)$ of $\bfff(\theta)$ are analytic with $\bfff(\theta)$ Hermitian for every $\theta \in [-\pi,\pi]$, then there exists an ordering (not necessarily the usual ordering by increasing magnitude) such that  $\hat{\lambda}_k(\bfff(\theta))$ are analytic real-valued functions for every $k=1,\ldots,\nu$ (see for example \cite[Chapter 2]{Kato}). Clearly, this kind of ordering does not affect the validity of Proposition \ref{prop:d-block-symbol_rearrangment}. In particular,  the function $g(\theta) = \sum_{k=1}^\nu \hat{\lambda}_k\left(\bfff_k(\theta)\right)$, as in Definition \ref{def:rearrangment} is piecewise analytical: it can have a finite number of jumps or generic points of non-analyticity at $\theta_k= \frac{(2k-\nu)\pi}{\nu}$ for $k=1,\ldots,\nu-1$, and so its image set $R_\bfgg$ can be a finite union of disjoint closed intervals. Then, by an appropriate modification of \cite[Lemma 2.3]{Kawohl}, the monotone rearrangement $\tilde{g}$ is piecewise Lipschitz continuous, and hypothesis of Theorem \ref{thm:discrete_Weyl_law} holds. Clearly, $\tilde{g}$ is differentiable almost everywhere on $[0,1]$ and it is differentiable in $x=1$ iff $\lim_{x\to 1}\tilde{g}'(x)<\infty$. Let us prove then that $\lim_{x\to 1}\tilde{g}'(x)=0$. Suppose for the moment that there exists only one interior point $\theta_M\in (-\pi,\pi), \theta_M\neq \theta_k$, where $g$ achieves its maximum $M$, i.e., such that $g(\theta_M)=M=\max_{\theta\in[-\pi,\pi]}g(\theta)=\tilde{g}(1)$. By regularity of $g$, there exists then $\delta>0$ such that $g$ is invertible on $I_l=(\theta_M-\delta,\theta_M]$ and $I_r=[\theta_M,\theta_M+\delta)$.  By equations \eqref{eq:rearrangment}-\eqref{eq:rearrangment2}, there exists $\epsilon=\epsilon(\delta)$ such that 
\begin{align*}
\phi(t) &= \phi(M-\epsilon) + g^{-1}_{|I_l}(t) + g^{-1}_{|I_r}(t),\\
&= \phi(M-\epsilon) + (\hat{\lambda}_k\circ \bfff_k)^{-1}_{|I_l}(t)+ (\hat{\lambda}_k\circ \bfff_k)^{-1}_{|I_r}(t) \qquad \mbox{for every } t \in [M-\epsilon,M], \mbox{  for some } k=1,\ldots,\nu,
\end{align*}
and then,
$$
\tilde{g}'(x) = \frac{1}{[(g^{-1})'_{|I_l}(t) + (g^{-1})'_{|I_r}(t)]_{|t=\tilde{g}^{-1}(x)}}
$$
for every $x \in (1-\eta,1]$ with $\eta=\eta(\epsilon)$ small enough. Since $t \to M$ as $x\to 1$ and  $g'_{|I_l}(\theta_M)=g'_{|I_r}(\theta_M)=g'(\theta_M)=0$, it follows that $\tilde{g}'(1)=0$. In the case that $\theta_M=\theta_k$ or $\theta_M=\pm \pi$, then
$$
g^{-1}_{|I_l}(t)= (\hat{\lambda}_k\circ \bfff_k)^{-1}_{|I_l}(t), \qquad g^{-1}_{|I_r}(t)= (\hat{\lambda}_{k+1}\circ \bfff_{k+1})^{-1}_{|I_r}(t),
$$ 
and we can conclude again that $\tilde{g}'(1)=0$ if we show that 
$$
[\hat{\lambda}_{k}(\bfff_k(\theta_k))]'=[\hat{\lambda}_{k}(\bfff(\pi))]'= 0 = [\hat{\lambda}_{k+1}(\bfff(-\pi))]'=[\hat{\lambda}_{k+1}(\bfff_{k+1}(\theta_k))]'.
$$
Due to the peculiar structure of the matrix function $\bfff$ from \eqref{eq:d-diamond-symbol}, for $d=1$, by direct computation it can be shown that the coefficients of its characteristic polynomial $\chi(\lambda)$ are functions of $(\cos(t_1\theta),\ldots, \cos(t_m\theta))$. Therefore, $\hat{\lambda}_k(\bfff(\theta))=h_k\left((\cos(t_1\theta),\ldots, \cos(t_m\theta))\right)$ with $h_k$ an analytic function, and then it follows easily that $[\hat{\lambda}_k(\bfff(\pm \pi))]'=0$. The generalization to the case of existence of countable many points $\theta_{M,j}$ such that $g(\theta_{M,j})=M$ is straightforward.

Finally, to prove \eqref{eq:convergence_zero2} we suppose $\lambda_{\nu n}^{(n)} > \lambda_{\nu n-1}^{(n)}$, otherwise the thesis would be trivial. By \cite[Theorem 2.2]{Serra1998B} we have that $\lambda_{\nu n}^{(n)}\leq \max_{\theta\in[-\pi,\pi]}g(\theta)$. Therefore all the hypotheses of Corollary \ref{cor:differentiability} are satisfied and the thesis follows at once.   
\hfill \ \, $\bullet$ \ \\

\begin{rem}
In general, $\tilde{g}$ is absolute continuous which means that it is a.e. differentiable. Therefore, it is differentiable in $x=1$ iff $\sup_{x\in[0,1]}\tilde{g}'(x)<\infty$. If it happens that $\lim_{x\to 1}\tilde{g}'(x)=\infty$, then $d_\bfnn(\lambda_{d_\bfnn}^{(\bfnn)}-\lambda_{d_\bfnn-1}^{(\bfnn)})\to \infty$ but it could diverge at a different rate with respect to $\tilde{g}'(x_\bfnn)$.
\end{rem}

\subsection{Examples}

\noindent\textit{Example 1.}\\

\begin{figure}[th]
	\centering
	\subfloat[Comparison between the distributions of $f$ and $\lambda_k^{(n)}$]{
		\includegraphics[height=32mm]{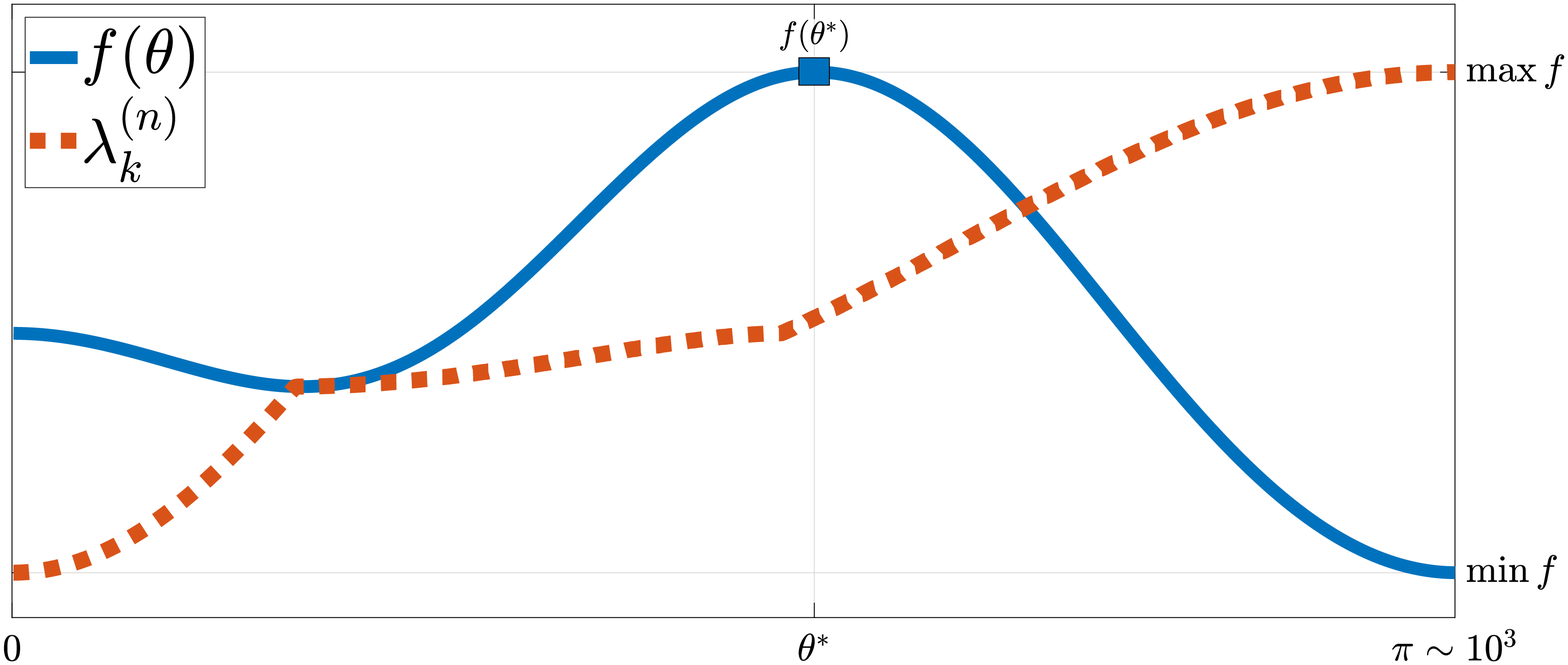}  
		\label{subfig:example_1A}
	}
	\centering
	\subfloat[Comparison between the distributions of $\tilde{f}$ and $\lambda_k^{(n)}$]{
		\includegraphics[height=32mm]{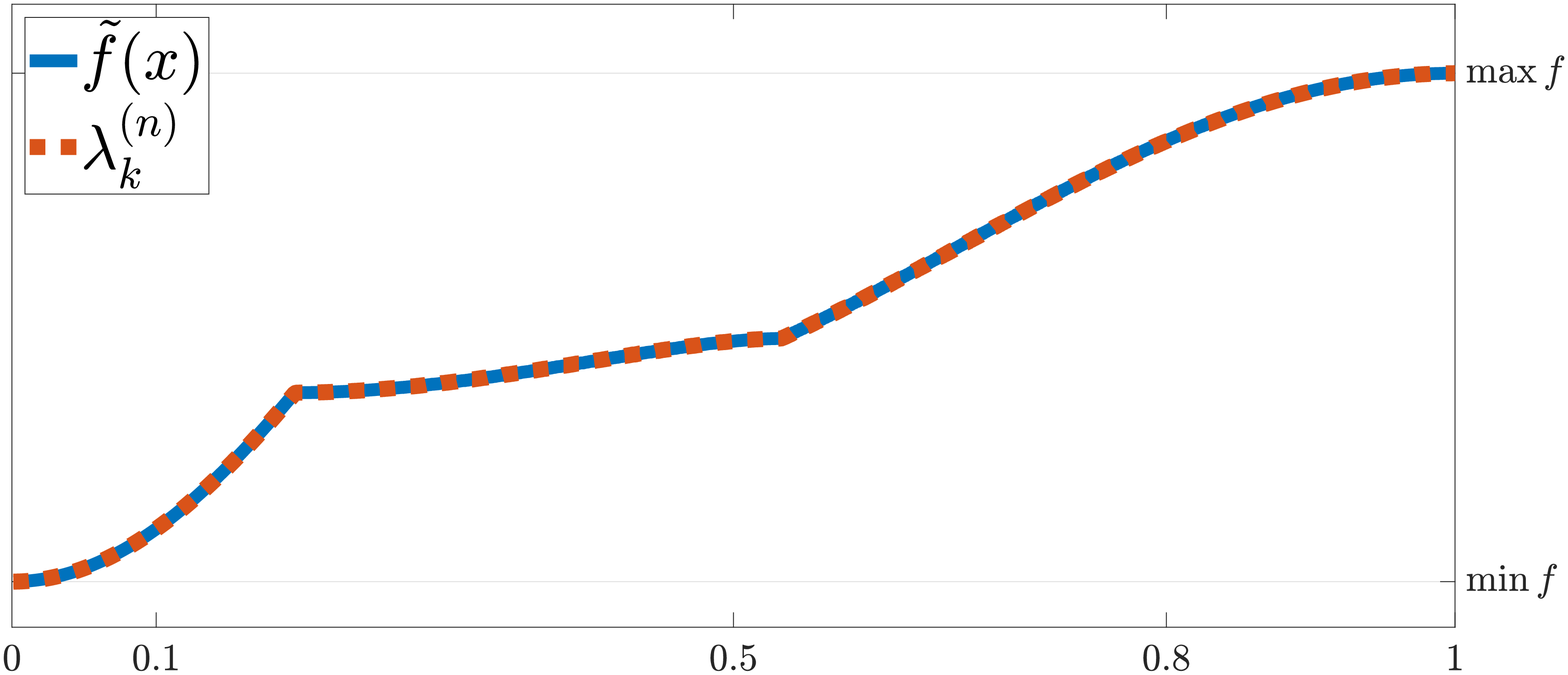}  
		\label{subfig:example_1B}
	}
	\caption{In Figure \ref{subfig:example_1A} we compare the distributions of $f(\theta)$, for $\theta \in [0,\pi]$ uniformly sampled over $10^3$ sampling points, and of the eigenvalues $\lambda_k^{(n)}$ of the adjacency matrix of $T_n\left\langle (1,1), (2,-6), (3,1), (4,1) \right\rangle$,  for $k=1,\ldots,10^3$. With the notations of Definition \ref{def:rearrangment}, notice that $\lambda_k^{(n)} \in [\min R_f, \max R_f]=[\min_{[0,\pi]}f(\theta), \max_{[0,\pi]}f(\theta)]$. In Figure \ref{subfig:example_1B}, instead, it is possible to observe the validity of Theorem \ref{thm:discrete_Weyl_law}, comparing the distribution of an approximation of the monotone rearrangement $\tilde{f}$, for $x\in[0,1]$, and the distribution of $\lambda_k^{(n)}$, for $k(n)/n \in [0,1]$ and $n=10^3$.}\label{fig:example_1}
\end{figure}

As a first example we consider the Toeplitz graph $T_n\left\langle (1,1), (2,-6), (3,1), (4,1) \right\rangle$ with corresponding symbol function $f(\theta)= 2\cos(\theta) - 12 \cos(2\theta) + 2\cos(3\theta) + 2\cos(4\theta)$, according to Proposition \ref{prop:d-level-distribution}. In this case, by symmetry of the symbol $f$ over $[-\pi,\pi]$ we can restrict it to $\theta\in [0,\pi]$ without affecting the validity of the identity \eqref{distribution:sv-eig}. It is easy to verify that the graph is connected, and in Figure \ref{fig:example_1} and Table \ref{table:ex_1} it is possible to check the numerical validity of Theorem \ref{thm:diamond-spectral-gap}, Corollary \ref{cor:diamond-spectral-gap} and Theorem \ref{thm:discrete_Weyl_law}.

\begin{table}[H]
	\centering
	\refstepcounter{table}
	\label{table:ex_1}
	\begin{tabular}{|c|c|cccc|} 
		\cline{3-6}
		\multicolumn{1}{c}{}               & \multicolumn{1}{l|}{}          & \multicolumn{4}{c|}{relative errors}                                                  \\ 
		\cline{3-6}
		\multicolumn{1}{c}{}               &                                & $n=10^2$             & $n=5\cdot 10^2$      & $n=10^3$             & $n=2\cdot 10^3$  \\ 
		\hline
		\multirow{4}{*}{$\frac{k(n)}{n}$ } & 0.1                            &  0.0039                    &             8.1567e-04         &              4.0990e-04        &  2.0547e-04                \\ 
		\cline{2-6}
		& 0.5        & 0.0013 & 6.7743e-04 & 0.0025 &     1.0263e-04             \\ 
		\cline{2-6}
		& 0.8                            &   0.0502                   &      0.0097                &               0.0035       &        0.0019          \\ 
		\cline{2-6}
		& 1                              &   0.0028                   &   1.1539e-04                   &           3.0532e-05           &     6.1804e-06             \\ 
		\hhline{|==|====|}
		\multicolumn{2}{|c|}{$n[\lambda_{n}^{(n)}-\lambda_{n-1}^{(n)}]$ }   & 0.1662 & 0.0045 & 0.0019 &        2.0245e-04          \\
		\hline
	\end{tabular}
\caption{The first four rows present the relative errors between the eigenvalue $\lambda_{k(n)}^{(n)}$ and the sampling of the monotone rearrangement $\tilde{f}\left(\frac{k(n)}{n}\right)$, i.e.: $\left|\frac{\lambda_{k(n)}^{(n)}}{\tilde{f}(k(n)/n)}-1\right|$. The index $k(n)$ of the eigenvalue $\lambda_{k(n)}^{(n)}$ is chosen such that $k(n)/n$ is constant for every fixed $n=10^2, 5\cdot10^3, 10^3, 2\cdot10^3$. As it can be seen, the relative errors decrease as $n$ increases, in accordance with Theorem \ref{thm:discrete_Weyl_law}. The convergence to zero is slow and not uniform: one of the reasons is that we are using a linear approximation of $\tilde{g}$ instead of $\tilde{g}$ itself. In the last row we show the computation of the gap between the biggest and the second-biggest eigenvalues, confirming the prediction of Corollary \ref{cor:diamond-spectral-gap}.}
\end{table}

\noindent\textit{Example 2.}\\

For this example we consider a $2$-level Toeplitz graph on $[0,1]^2$ given by 
$$
T_\bfnn \langle \{[\bftt_1], \bfww_1 \}, \{[\bftt_2], \bfww_2^p \}, \{[\bftt_3], \bfww_3 \},\{[\bftt_4], \bfww_4 \} \rangle,
$$
 like in Definition \ref{local gr cube}. We set $\bfnn=(n,n)$ and
\begin{align*}
 [\bftt_1]= [(1,0)]=\{[(1,0)]_1\}, \quad [\bftt_2]= [(0,1)]=\{[(0,1)]_1\}, \quad &[\bftt_3]= [(1,1)]=\{[(1,1)]_1, [(1,1)]_2\}, \\
 &[\bftt_4]= [(2,2)]=\{[(2,2)]_1, [(2,2)]_2\},
\end{align*}
where
\begin{align*}
[(1,0)]_1=\left\{ \pm (1,0) \right\}, \quad [(0,1)]_1=\left\{ \pm (0,1) \right\}, \quad &[(1,1)]_1=\{\pm(1,-1)\}, \quad [(1,1)]_2=\{\pm(1,1)\},\\
&[(2,2)]_1=\{\pm(2,-2)\}, \quad [(2,2)]_2=\{\pm(2,2)\}.
\end{align*}
Finally,
\begin{equation*}
  \bfww_1=w_{1,0}=1, \quad \bfww_2=w_{0,1}=2, \quad \bfww_3=(w_{1,\shortminus 1}, w_{1,1})=(-3,-3), \quad \bfww_4=(w_{2,\shortminus 2}, w_{2,2})=(1,1)
\end{equation*}
By Theorem \ref{thm:d-grid-graphs-symbol}, the sequence of adjacency matrices $\{W_\bfnn\}_\bfnn$ has symbol $f:  [-\pi,\pi]^2 \to \R$,
\begin{equation*}
f(\theta_1,\theta_2) = 2\cos(\theta_1) + 4\cos(\theta_2) -6\cos(\theta_1-\theta_2) -6\cos(\theta_1+\theta_2) +2\cos(2\theta_1-2\theta_2)+2\cos(2\theta_1+2\theta_2).
\end{equation*}
Due to the symmetry of the symbol $f$ over $[-\pi,\pi]^2$ we can restrict it to $(\theta_1,\theta_2)\in [0,\pi]^2$ without affecting the validity of the identity \eqref{distribution:sv-eig}. In the following Table \ref{table:ex_2} and Figure \ref{fig:ex_2} it is possible to check numerically the validity of Theorem \ref{thm:discrete_Weyl_law} and Theorem \ref{thm:diamond-spectral-gap}.

\begin{table}[H]
	\begin{minipage}{0.52\textwidth}
		\centering
	\begin{tabular}[m]{|c|c|ccc|} 
		\cline{3-5}
		\multicolumn{1}{l}{}               & \multicolumn{1}{l|}{}          & \multicolumn{3}{c|}{relative errors}                  \\ 
		\cline{3-5}
		\multicolumn{1}{l}{}               &                                & $n=10$              & $n=50$  & $n=100$   \\ 
		\hline
		\multirow{4}{*}{$\frac{k(\bfnn)}{d_\bfnn}$ } & 0.2                            &    0.0422                 &         0.0053         &    0.0029        \\ 
		\cline{2-5}
		& 0.5                            &  0.2553                  &    0.1211          &  0.0705       \\ 
		\cline{2-5}
		& 0.7                           &     0.0396                 &  0.0096             &         0.0089   \\ 
		\cline{2-5}
		& 1                              &   0.0071                  &      7.8515e-05            &     1.0184e-05      \\ 
		\hhline{=====}
		\multicolumn{2}{|c|}{$|\frac{\lambda_{d_\bfnn}^{(\bfnn)}-\lambda_{d_\bfnn\shortminus 1}^{(\bfnn)}}{\tilde{f}^{(\bfnn)}_{d_\bfnn}-\tilde{f}^{(\bfnn)}_{d_\bfnn\shortminus 1}}-1|$ }   & 0.0487 &         0.0145         &     0.0075       \\
		\hline
	\end{tabular}\caption{The first four rows present the relative errors between the eigenvalue $\lambda_{k(\bfnn)}^{(\bfnn)}$ and the sampling of the monotone rearrangement $\tilde{f}\left(\frac{k(\bfnn)}{d_\bfnn}\right)$. The index $k(\bfnn)$ is chosen such that $k(\bfnn)/d_\bfnn$ is constant for every fixed $\bfnn=(n,n)$ with $n=10, 50, 100$. The relative errors decrease as $n$ increases, in accordance with Theorem \ref{thm:discrete_Weyl_law}. The convergence is slow and not uniform: one of the reasons is that we are using a simple linear approximation of $\tilde{f}$. In the last row we show relative errors between the gap of the biggest and the second-biggest eigenvalues and the gap between the maximum and second-maximum values of a uniform sampling of $\tilde{f}$, confirming the prediction of Theorem \ref{thm:diamond-spectral-gap}. In particular, $d_\bfnn\left(\lambda_{d_\bfnn}^{(\bfnn)}-\lambda_{d_\bfnn\shortminus 1}^{(\bfnn)}\right)\to \tilde{f}'(1)\in (0,\infty)$.}\label{table:ex_2}
\end{minipage}\hfill
\begin{minipage}[m]{0.45\textwidth}
	\centering
	\includegraphics[height=32mm]{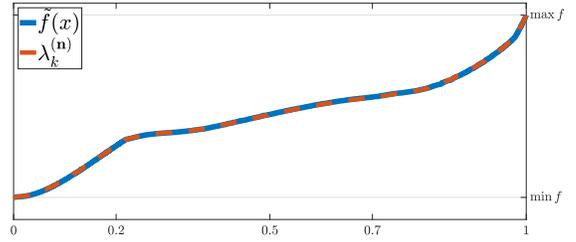}  
	\captionof{figure}{It is possible to observe the validity of Theorem \ref{thm:discrete_Weyl_law}, comparing the distribution of an approximation of the monotone rearrangement $\tilde{f}$, for $x\in[0,1]$, and the distribution of $\lambda_k^{(\bfnn)}$, for $k(\bfnn)/d_\bfnn \in [0,1]$ and $\bfnn=(60,60)$, i.e., $d_\bfnn=36\cdot10^2$.}
		\label{fig:ex_2}
	\end{minipage}
\end{table}

\noindent\textit{Example 3.}\\

\begin{figure}
	\centering
	\subfloat[Comparison between the distributions of $\bfff$ and $\lambda_k^{(n)}$]{
		\includegraphics[height=32mm]{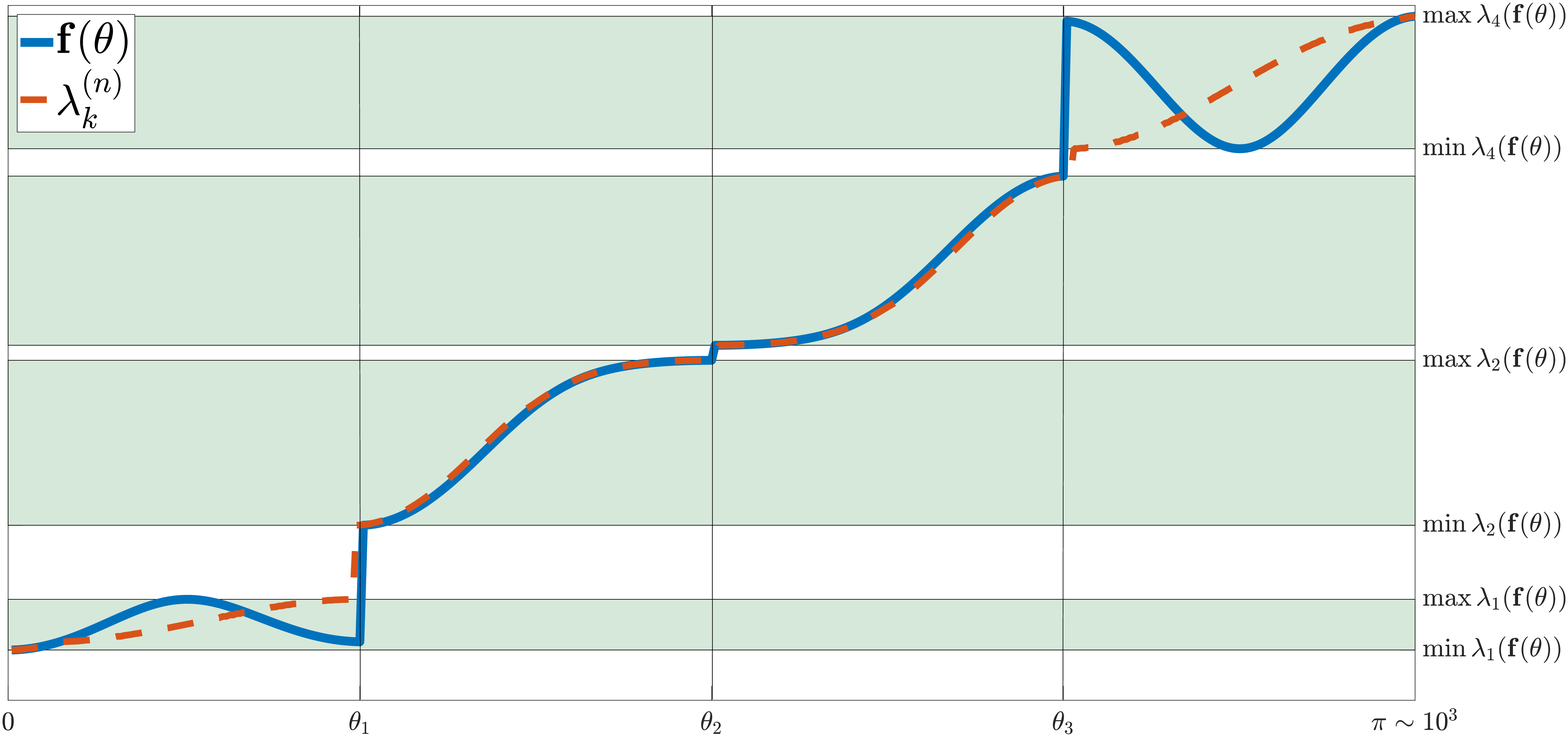}  
		\label{subfig:example_3A}
	}
	\centering
	\subfloat[Comparison between the distributions of $\tilde{f}$ and $\lambda_k^{(n)}$]{
		\includegraphics[height=32mm]{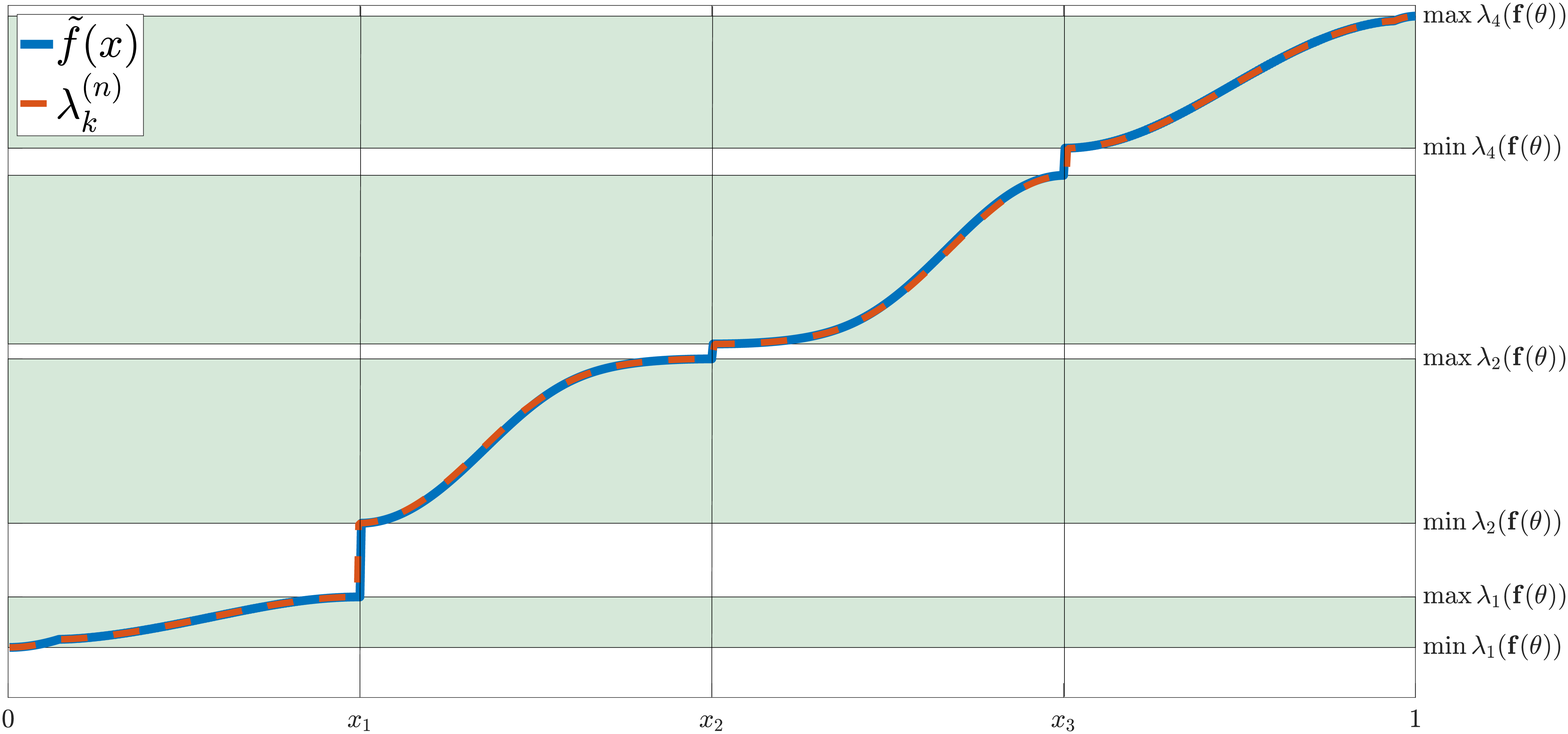}  
		\label{subfig:example_3B}
	}

\subfloat[Maginification of distribution graphs of $\tilde{f}$ and $\lambda_k^{(n)}$ in a neighborood of $x_2$.]{
	\includegraphics[height=32mm]{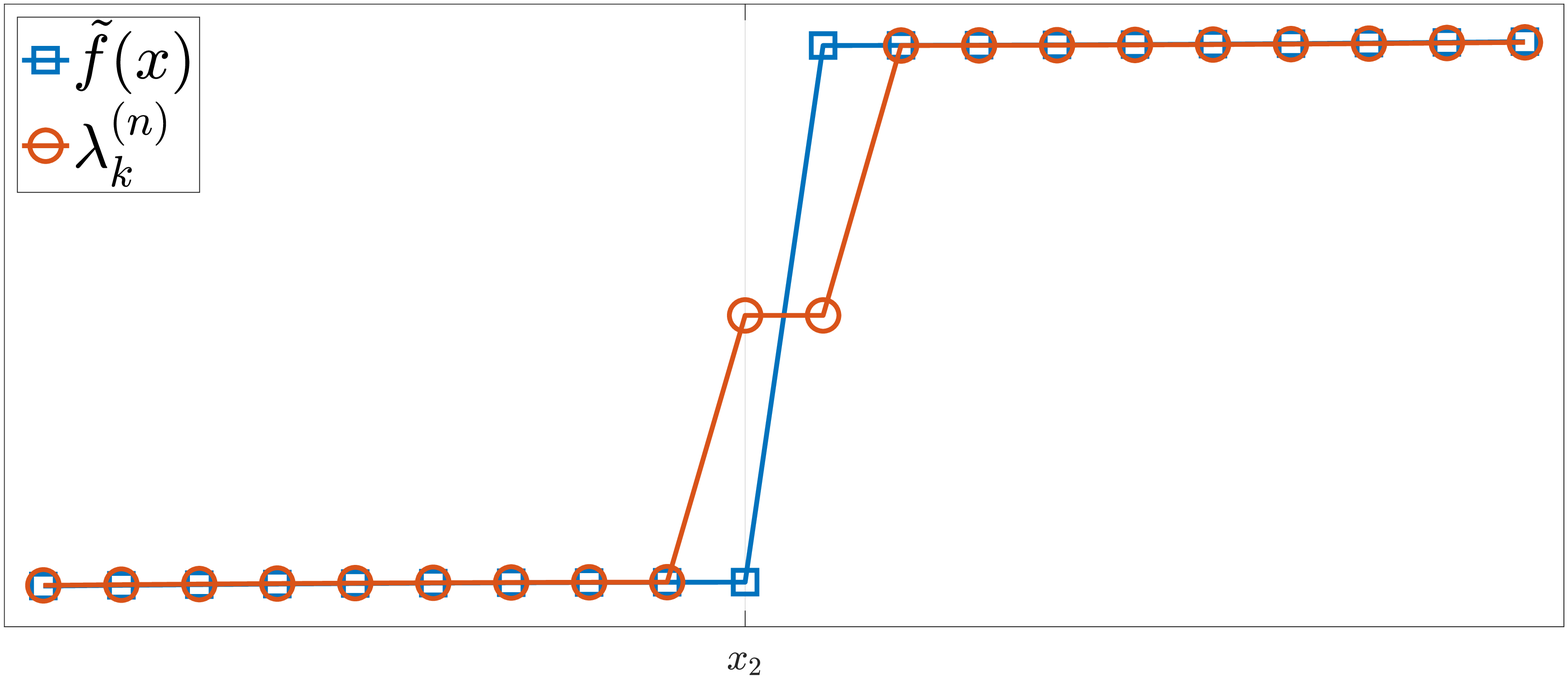}  
	\label{subfig:example_3C}
}
	\caption{In Figure \ref{subfig:example_3A} we compare the distributions of $\bfff(\theta)$ from \eqref{eq:real_symbol_diamond}, for $\theta \in [0,\pi]$, and of the eigenvalues $\lambda_k^{(n)}$ of the adjacency matrix $W_n^G$,  for $k=1,\ldots,10^3$. We have that $\theta_1=\pi/4$, $\theta_2=\pi/2$, $\theta_3=3\pi/4$ are points of discontinuity since the inequalities in \eqref{eq:discontinuities} are strict.  With the notations of Definition \ref{def:rearrangment}, notice that in this example $R_\bfff$ is the union of four disjoint intervals and that $\lambda_k^{(n)} \in [\min R_\bfff, \max R_\bfff]=[\min_{[0,\pi]}\lambda_1(\bfff(\theta)), \max_{[0,\pi]}\lambda_4(\bfff(\theta))]$. In Figure \ref{subfig:example_1B}, instead, it is possible to observe the validity of Theorem \ref{thm:discrete_Weyl_law}, comparing the distribution of an approximation of the monotone rearrangement $\tilde{f}$, for $x\in[0,1]$, and the distribution of $\lambda_k^{(n)}$, for $k(n)/(4n) \in [0,1]$ and $4n=10^3$. The points of discontinuity are now located at $x_1=1/4, x_2=1/2, x_3=3/4$. Finally, in figure \ref{subfig:example_3C} we zoom in a neighborhood of the discontinuity point $x_2$: it is possible to observe the presence of two outliers which do not belong to $R_\bfff$. This is not in contradiction to Theorem \ref{thm:discrete_Weyl_law} since by Corollary \ref{cor:weak_clustering} it admits at most a number of $o(4n)$ of outliers.}\label{fig:example_3}
\end{figure}

We consider a sequence of adjacency matrices $\{W_n^G\}_n$ from the $1$-level diamond Toeplitz graph given in Figure \ref{fig:diamond_graph}, namely $T_n^G\langle (1,L_1), (2,L_2) \rangle$ with mold graph $G=T_4\langle (1,1),(3,1) \rangle$ and 
\begin{equation*}
W=\begin{tikzpicture}[baseline=-\the\dimexpr\fontdimen22\textfont2\relax]
\matrix (m)[matrix of math nodes,left delimiter=(,right delimiter=), every node/.style={anchor=base,text depth=.5ex,text height=2ex,text width=1em}]
{ 0 & $1$ & 0 & $1$\\
	$1$ & 0 & $1$ & 0\\
	0& $1$ & 0 & $1$\\
	$1$ & 0 & $1$ & 0\\
};
\end{tikzpicture}
\qquad 
L_1=\begin{tikzpicture}[baseline=-\the\dimexpr\fontdimen22\textfont2\relax]
\matrix (m)[matrix of math nodes,left delimiter=(,right delimiter=), every node/.style={anchor=base,text depth=.5ex,text height=2ex,text width=1em}]
{ $\shortminus 2$ & 0 & 0 & 0\\
	0 & 0 & 0 & 0\\
	0& 0 & 0 & 0\\
	0 & 0 & 0 & 0\\
};
\end{tikzpicture}\qquad
L_2=\begin{tikzpicture}[baseline=-\the\dimexpr\fontdimen22\textfont2\relax]
\matrix (m)[matrix of math nodes,left delimiter=(,right delimiter=), every node/.style={anchor=base,text depth=.5ex,text height=2ex,text width=1em}]
{ 0 & 0 & 0 & 0\\
	0 & 0 & 0 & 0\\
	0& 0 & $\frac{1}{2}$ &0\\
	0 & 0 & $6$& 0\\
};
\end{tikzpicture},
\end{equation*}
where $W$ is the adjacency matrix of the mold graph $G$. From Proposition \ref{prop:diamond-distribution}, $\{W_n^G\}_n \sim_\lambda \bfff(\theta)$, where
\begin{equation}\label{eq:original_symbol_diamond}
\bfff(\theta)=W + 2L_1\cos(\theta) + \left(L_2+L_2^T\right)\cos(2\theta) + \left(L_2-L_2^T\right)\i\sin(2\theta).
\end{equation}
In this case, by symmetry of the symbol $\bfff$ over $[-\pi,\pi]$ we can restrict it to $\theta\in [0,\pi]$ without affecting the validity of the identity \eqref{distribution:sv-eig}. Moreover, due to Proposition \ref{prop:d-block-symbol_rearrangment} and taking into account that we restricted $\theta$ to $[0,\pi]$, with abuse of notation we write
\begin{equation}\label{eq:real_symbol_diamond}
\bfff(\theta) = \sum_{k=1}^4 \mathds{1}_{I_k}(\theta) \lambda_k\left(\bfff_k(\theta)\right),
\end{equation}
where 
\begin{equation*}
\bfff_k(\theta)=
\bfff(4\theta - (k-1)\pi) \qquad 
I_k=\left[\frac{(k-1)\pi}{4},\frac{k\pi}{4}\right].
 \end{equation*}
The map $\theta \mapsto 4\bftheta - (k-1)\pi$ is a diffeomorphism between $I_k$ and $[0,\pi]$, and the eigenvalues functions are ordered by magnitude, namely, $\lambda_k\left(\bfff_k(\theta)\right)$ is the $k$-th eigenvalue function of the matrix \eqref{eq:original_symbol_diamond} over $\theta \in [0,\pi]$. In particular, in this case it holds that 
\begin{equation}\label{eq:discontinuities}
\max_{\theta\in I_k}\lambda_k(\bfff_k(\theta)) = \max_{\theta\in[0,\pi]}\lambda_k(\bfff(\theta))< \min_{\theta\in[0,\pi]}\lambda_{k+1}(\bfff(\theta))= \min_{\theta\in I_k} \lambda_{k+1}(\bfff_{k+1}(\theta)), \quad \mbox{for }k=1,2,3.
\end{equation}
 In the following Table \ref{table:ex_3} and Figure \ref{fig:example_3} it is possible to check numerically the validity of Theorem \ref{thm:discrete_Weyl_law} and Corollary \ref{cor:diamond-spectral-gap}.

\begin{table}[H]
	\centering
	\begin{tabular}{|c|c|ccc|} 
		\cline{3-5}
		\multicolumn{1}{l}{}               & \multicolumn{1}{l|}{}          & \multicolumn{3}{c|}{relative errors}                  \\ 
		\cline{3-5}
		\multicolumn{1}{l}{}               &                                & $4n=10^2$              & $4n=10^3$  & $4n=5\cdot 10^3$   \\ 
		\hline
		\multirow{4}{*}{$\frac{k(n)}{4n}$ } & 0.1                            &    5.2984e-04                 &        3.8322e-05         &    1.0184e-05        \\ 
		\cline{2-5}
		& 0.4                            &     0.0285                  &       0.0029          &     5.7782e-04      \\ 
		\cline{2-5}
		& 0.8                            &        0.0074             &     6.2919e-04            &          1.2811e-04  \\ 
		\cline{2-5}
		& 1                              &            4.5189e-06          &        4.6087e-09          &       3.6928e-11      \\ 
		\hhline{=====}
		\multicolumn{2}{|c|}{$4n[\lambda_{4n}^{(n)}-\lambda_{4n-1}^{(n)}]$ }   & 3.2287 &   0.3253               &     0.0651       \\
		\hline
	\end{tabular}
\caption{The first four rows present the relative errors between the eigenvalue $\lambda_{k(n)}^{(n)}$ and the sampling of the monotone rearrangement $\tilde{f}\left(\frac{k(n)}{4n}\right)$, i.e.: $\left|\frac{\lambda_{k(n)}^{(n)}}{\tilde{f}(k(n)/4n)}-1\right|$. The index $k(n)$ of the eigenvalue $\lambda_{k(n)}^{(n)}$ is chosen such that $k(n)/(4n)$ is constant for every fixed $4n=10^2,10^3, 5\cdot10^3$. As it can be seen, the relative errors decrease as $n$ increases, in accordance with Theorem \ref{thm:discrete_Weyl_law}. The convergence is slow and not uniform:  one of the reasons is that we are using a linear approximation of $\tilde{g}$ instead of $\tilde{g}$ itself. In the last row we show the computation of the gap between the biggest and the second-biggest eigenvalues, confirming the prediction of Corollary \ref{cor:diamond-spectral-gap}.}\label{table:ex_3}
\end{table}

\section{Applications to PDEs approximation}\label{sec:appl}
The section is divided into three parts where we show that the approximation of a model differential problem by three celebrated approximation techniques leads to sequences of matrices that fall in the theory developed in the previous sections. 

\subsection{Approximations of PDEs vs sequences of weighted d-level grid graphs: FD}\label{ssec:appl-0}

As a first example we consider the discretization of a self-adjoint operator $\mathcal{L}$ with (homogeneous) Dirichlet boundary conditions (BCs) on the disk $B_{1/2}\subset[0,1]^2$ by an equispaced Finite Difference (FD) approximation with $(4m+1)$-points. That is, our model operator with Dirichlet BCs is given by
\begin{align}
&\mathcal{L}:W^{1,2}_0\left(B_{\frac{1}{2}}\right) \to \textnormal{L}^2\left(B_{\frac{1}{2}}\right), \qquad B_{\frac{1}{2}}=\left\{(x,y)\in \R^2\, : \, 4\left(x-\frac{1}{2}\right)^2+4\left(y-\frac{1}{2}\right)^2<1\right\}, \label{eq:SLO1}\\
&\mathcal{L}[u](x,y):= -\textnormal{div}\left[p(x,y)\nabla u(x,y)\right] + q(x,y)u(x,y) \qquad (x,y)\in B_{\frac{1}{2}}.\label{eq:SLO2}
\end{align}

Fixing the diffusion term $p(x,y)=1+(x-1/2)^2+(y-1/2)^2$ and the potential term $q(x,y)=\textrm{e}^{xy}$, then the operator $\mathcal{L}$ is self-adjoint and has purely discrete spectrum, see \cite{Davies}. 

Now, if we fix $m=1$, $n\in \N$ and $i,j \in \mathbb{Z}$, then the uniform second-order $5$-point FD approximation of the (negative) Laplacian operator (i.e., $\Delta[\cdot]=-(\partial^2_{x^2}+\partial^2_{y^2})[\cdot]$) is given by
\begin{equation*}
\Delta[u](x_i,y_j) \approx h^{-2}\left( -u(x_{i},y_{j-1}) - u(x_{i-1},y_j) + 4u(x_i,y_j) - u(x_{i+1},y_j)  -u(x_{i},y_{j+1})  \right)
\end{equation*}
for every $u\in C^\infty\left(\R^2\right)$, where $h={(n+1)^{-1}},x_i={\frac{i}{n+1}}, y_j={\frac{j}{n+1}}$. The same approximation applies for every $u\in C^\infty\left(B_{1/2}\right)$ and $i,j,n$ such that $\{(x_i,y_j),(x_{i\pm1},y_{j\pm1})\}\subset B_{1/2}$. Notice that the weight of the central point $u(x_i,y_j)$ is the sum of all the other weights, changed of sign. 

\begin{figure}[H]
\begin{center}
\begin{tikzpicture}[whitestyle/.style={circle,draw,fill=white!40,minimum size=4},
graystyle/.style ={ circle ,top color =white , bottom color = gray ,
	draw,black, minimum size =4},	state2/.style ={ circle ,top color =white , bottom color = gray ,
	draw,black , text=black , minimum width =1 cm},
state3/.style ={ circle ,top color =white , bottom color = white ,
draw,white , text=white , minimum width =1 cm},
state/.style={circle ,top color =white , bottom color = white,
draw,black , text=black , minimum width =1 cm}]
\draw[step=0.5] (-3,-3) grid (3,3);
\draw[step=0.5, red,very thick] (-2.5,-2.5) grid (2.5,2.5);
\draw[step=0.5, ForestGreen,very thick] (-2,-2) grid (2,2);
\draw[step=0.5, ForestGreen, very thick] (-1.5,-2.5) grid (1.5,2.5);
\draw[step=0.5, ForestGreen, very thick] (-2.5,-1.5) grid (2.5,1.5);
\draw[very thick] (0,0) circle (3); 
\draw[very thick, blue] (2.25,-2) ellipse (0.4cm and 0.8cm); 
\draw[blue, very thick, ->] (2.65,-2) to (4,-0.5);
\draw[very thick, blue] (-2.5,0) ellipse (0.4cm and 0.8cm);

\foreach \y in {-1.5,...,1.5}
\draw[red,very thick] (2.5,\y) to (3,\y);
\foreach \y in {-1,...,1}
\draw[red,very thick] (2.5,\y) to (3,\y);

\foreach \y in {-1.5,...,1.5}
\draw[red,very thick] (-2.5,\y) to (-3,\y);
\foreach \y in {-1,...,1}
\draw[red,very thick] (-2.5,\y) to (-3,\y);

\foreach \x in {-1.5,...,1.5}
\draw[red,very thick] (\x,-2.5) to (\x,-3);
\foreach \x in {-1,...,1}
\draw[red,very thick] (\x,-2.5) to (\x,-3);

\foreach \x in {-1.5,...,1.5}
\draw[red,very thick] (\x,2.5) to (\x,3);
\foreach \x in {-1,...,1}
\draw[red,very thick] (\x,2.5) to (\x,3);

\foreach \x in {-1,...,1}
\node [whitestyle] at (\x,2.5) {};
\foreach \x in {-0.5,...,0.5}
\node [whitestyle] at (\x,2.5) {};

\foreach \x in {-2,...,2}
\node [whitestyle] at (\x,2) {};
\foreach \x in {-1.5,...,1.5}
\node [whitestyle] at (\x,2) {};

\foreach \x in {-2,...,2}
\node [whitestyle] at (\x,1.5) {};
\foreach \x in {-1.5,...,1.5}
\node [whitestyle] at (\x,1.5) {};

\foreach \x in {-2.5,...,2.5}
\node [whitestyle] at (\x,1) {};
\foreach \x in {-2,...,2}
\node [whitestyle] at (\x,1) {};

\foreach \x in {-2.5,...,2.5}
\node [whitestyle] at (\x,0.5) {};
\foreach \x in {-2,...,2}
\node [whitestyle] at (\x,0.5) {};

\foreach \x in {-2.5,...,2.5}
\node [whitestyle] at (\x,0) {};
\foreach \x in {-2,...,2}
\node [whitestyle] at (\x,0) {};

\foreach \x in {-2.5,...,2.5}
\node [whitestyle] at (\x,-0.5) {};
\foreach \x in {-2,...,2}
\node [whitestyle] at (\x,-0.5) {};

\foreach \x in {-2.5,...,2.5}
\node [whitestyle] at (\x,-1) {};
\foreach \x in {-2,...,2}
\node [whitestyle] at (\x,-1) {};

\foreach \x in {-2,...,2}
\node [whitestyle] at (\x,-1.5) {};
\foreach \x in {-1.5,...,1.5}
\node [whitestyle] at (\x,-1.5) {};

\foreach \x in {-2,...,2}
\node [whitestyle] at (\x,-2) {};
\foreach \x in {-1.5,...,1.5}
\node [whitestyle] at (\x,-2) {};

\foreach \x in {-1,...,1}
\node [whitestyle] at (\x,-2.5) {};
\foreach \x in {-0.5,...,0.5}
\node [whitestyle] at (\x,-2.5) {};

\node [whitestyle] at (1.5,2.5) {};
\node [whitestyle] at (1.5,-2.5) {};
\node [whitestyle] at (-1.5,2.5) {};
\node [whitestyle] at (-1.5,-2.5) {};

\node [whitestyle] at (-2.5,-1.5) {};
\node [whitestyle] at (2.5,1.5) {};
\node [whitestyle] at (-2.5,1.5) {};
\node [whitestyle] at (2.5,-1.5) {};
\node [graystyle] at (2.5,2.5) {};
\node [graystyle] at (2,2.5) {};
\node [graystyle] at (2.5,2) {};

\node [graystyle] at (-2.5,-2.5) {};
\node [graystyle] at (-2,-2.5) {};
\node [graystyle] at (-2.5,-2) {};

\node [graystyle] at (-2.5,2.5) {};
\node [graystyle] at (-2,2.5) {};
\node [graystyle] at (-2.5,2) {};

\node [graystyle] at (2.5,-2.5) {};
\node [graystyle] at (2,-2.5) {};
\node [graystyle] at (2.5,-2) {};
\node[state] at (5,2) (A) {$\kappa_0$};
\node[state]  (B) [right=of A] {$\kappa_2$};
\node[state]  (D) [below=of A] {$\kappa_2$};
\node[state2]  (E) [right=of D] {};
\node[state2]  (G) [below=of D] {};
\node[state2]  (H) [right=of G] {};
\path[-] (A) edge [ForestGreen] node[above] {$\textcolor{black}{w^p}$} (B);
\path[-] (G) edge [red] node[above] {$\textcolor{black}{w^{\bar{p}}}$} (H);
\path[-] (D) edge [red] node[above] {$\textcolor{black}{w^{\bar{p}}}$} (E);
\path[-] (H) edge [red] node[right] {$\textcolor{black}{w^{\bar{p}}}$} (E);

\path[-] (A) edge [ForestGreen] node[right] {$\textcolor{black}{w^p}$} (D);
\path[-] (B) edge [red] node[right] {$\textcolor{black}{w^{\bar{p}}}$} (E);
\path[-] (D) edge [red] node[right] {$\textcolor{black}{w^{\bar{p}}}$} (G);

\draw[ForestGreen,dashed] (B) to (7,3);
\draw[red,dashed] (B) to (8,2);
\draw[ForestGreen,dashed] (A) to (5,3);
\draw[ForestGreen,dashed] (A) to (4,2);
\draw[ForestGreen,dashed] (D) to (4,0);
\node[state] at (-5,0) (A1) {$\kappa_1$};
\node[state]  (B1) [above=of A1] {$\kappa_1$};
\node[state]  (C1) [below=of A1] {$\kappa_1$};
\node[state2]  (D1) [left=of A1] {};
\node[state2]  (E1) [left=of B1] {};
\node[state2]  (F1) [left=of C1] {};
\draw[blue, very thick, ->] (-2.9,0) to (-3.8,0);
\path[-] (A1) edge [ForestGreen] node[right] {$\textcolor{black}{w^p}$} (B1);
\path[-] (A1) edge [ForestGreen] node[right] {$\textcolor{black}{w^p}$} (C1);
\path[-] (A1) edge [red] node[below] {$\textcolor{black}{w^{\bar{p}}}$} (D1);
\path[-] (B1) edge [red] node[below] {$\textcolor{black}{w^{\bar{p}}}$} (E1);
\path[-] (C1) edge [red] node[below] {$\textcolor{black}{w^{\bar{p}}}$} (F1);
\path[-] (D1) edge [red] node[right] {$\textcolor{black}{w^{\bar{p}}}$} (F1);
\path[-] (D1) edge [red] node[right] {$\textcolor{black}{w^{\bar{p}}}$} (E1);

\draw[ForestGreen,dashed] (B1) to (-5,3);
\draw[ForestGreen,dashed] (C1) to (-5,-3);
\draw[ForestGreen,dashed] (C1) to (-4,-2);
\draw[ForestGreen,dashed] (A1) to (-4,0);
\draw[ForestGreen,dashed] (B1) to (-4,2);
\end{tikzpicture}
\end{center}\caption{Immersion of a $2$-level grid graph inside the disk $B_{1/2}\subset [0,1]^2$. The white nodes are the nodes of $V_\bfnn$ while the gray nodes belong to $\bar{V}_\bfnn\setminus V_\bfnn$.  The potential term $\bar{\kappa}$ of the mother graph $\bar{G}_n$ is determined only by the potential term $q$ from \eqref{eq:SLO1} while the potential term $\kappa$ of the sub-graph $G_\bfnn$ is influenced by the nodes in $\bar{V}_\bfnn\setminus V_\bfnn$ on the boundary set $\partial V_\bfnn$. This influence is due to the presence of Dirichlet BCs in \eqref{eq:SLO2}. The green connections represent the weighted edges whose end-nodes are both interior nodes of $V_\bfnn$, while the red connections represent the weighted edges which have at least one end-node that belongs to $\bar{V}_\bfnn\setminus V_\bfnn$. The potential term $\kappa$ sums the weight of a red edge to every of its end-nodes which belong to $V_\bfnn$.}\label{fig:immersion_in_the_disk}
\end{figure}
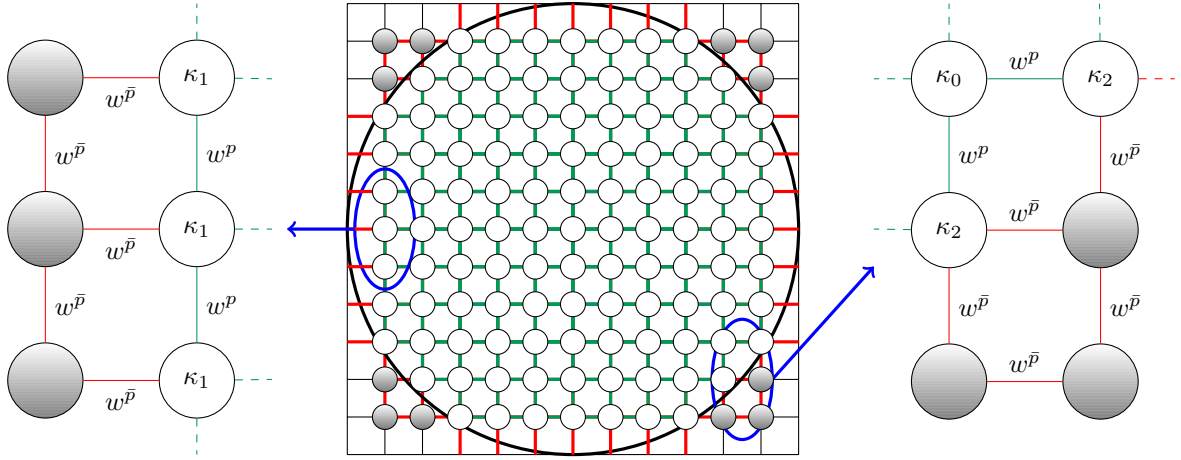

Let us consider the $2$-level Toeplitz graph $T_{\bfnn}\left\langle \left\{[1,0],1  \right\},\left\{[0,1],1  \right\}\right\rangle$ with $\bfnn=(n,n)$ as in Definition \ref{def:d-toeplitz-graph} and immerse it in $B_{1/2}$ as in Definition \ref{def:local gr Omega}, i.e., $T_{\bfnn}^{B_{1/2}}\left\langle \left\{[1,0],w^p  \right\},\left\{[0,1],w^p  \right\}\right\rangle$ such that
\begin{align*}
&V_{\bfnn}= \left\{(x_i,y_j) \in [0,1]^2 \, : \, (x_i,y_j) \in B_{\frac{1}{2}} \right\},\\
&w^p((x_i,y_j),(x_r,y_s))= 1+\left(\frac{x_i+x_r}{2}-\frac{1}{2}\right)^2 +\left(\frac{y_j+y_s}{2}-\frac{1}{2}\right)^2 \quad \mbox{if }(|i-r|,|j-s|)\in\{(1,0),(0,1)\}.
\end{align*}
Extend now continuously the diffusion term $p(x,y)$ outside $B_{1/2}$, that is,
\begin{equation}\label{eq:FD_extended_p}
\overline{p}(x,y) = \begin{cases}
\frac{5}{4} & \mbox{if } (x,y)\in \R^2\setminus B_{\frac{1}{2}},\\
p(x,y) &\mbox{if } (x,y)\in B_{\frac{1}{2}},
\end{cases}
\end{equation}
and define the graph
\begin{equation}\label{eq:5-FD-graph}
G_{\bfnn}=\left(T_{\bfnn}^{B_{1/2}}\left\langle \left\{[1,0],w^p  \right\},\left\{[0,1],w^p  \right\}\right\rangle, \kappa\right),
\end{equation}
as a sub-graph of 
\begin{align}\label{eq:5-FD-graph_extended}
\bar{G}_{\bfnn}=\left(T_{\bfnn}\left\langle \left\{[1,0],w^{\bar{p}}  \right\},\left\{[0,1],w^{\bar{p}}  \right\}\right\rangle, \bar{\kappa}\right) \qquad \mbox{where} \quad &\bar{V}_{\bfnn}= \left\{(x_i,y_j) \in [0,1]^2 \right\}\\
&\bar{p} \mbox{ as in }\eqref{eq:FD_extended_p},\nonumber\\
&\bar{\kappa} (x_i,y_j) =h^2q(x_i,y_j), \nonumber
\end{align}
like in Definition \ref{def:sub-graph}. Namely, the mother graph $\bar{G}_\bfnn$ is the $2$-level grid graph on $[0,1]^2$ obtained by extending continuously the diffusivity function $p$ to $[0,1]^2$ and adding a nontrivial potential term $\bar{\kappa}$ which naturally depends on $q$. On the other hand, the potential term $k$ of the  sub-graph $G_{\bfnn}$ describes the edge deficiency of nodes in $G_{\bfnn}$ compared to the same nodes in $\bar{G}_{\bfnn}$,
\begin{align*}
&\kappa(x_i,y_j)=h^2q(x_i,y_j) + \sum_{\substack{(x_r,y_s)\sim (x_i,y_j)\\(x_r,y_s)\in \bar{V}_{\bfnn}\setminus V_{\bfnn}}} w^{\bar{p}}((x_r,y_s),(x_i,x_j)),\\
&|\kappa(x_i,y_j)- \bar{\kappa}(x_i,y_j)|= \sum_{\substack{(x_r,y_s)\sim (x_i,y_j)\\(x_r,y_s)\in \bar{V}_\bfnn\setminus V_\bfnn}} w^{\bar{p}}((x_r,y_s),(x_i,x_j)),
\end{align*}
see \cite[pg.197]{KL}. It is easy to check that the boundary points $(x_i,x_j)\in \partial V_\bfnn$ are connected at most with two points of $\bar{V}_\bfnn\setminus V_\bfnn$, therefore the potential term $\kappa$ can be split into three terms
\begin{equation*}
\kappa(x_i,y_j)=\begin{cases}
\kappa_0(x_i,y_j)= h^2q(x_i,y_j) &\mbox{if } (x_i,y_j) \in \mathring{V}_\bfnn,\\
\kappa_1(x_i,y_j)= h^2q(x_i,y_j) + \deg((x_i,y_j)), &\mbox{if } \left|\left\{(x_r,y_s) \in \bar{V}_\bfnn\setminus V_\bfnn \, : \, (x_r,y_s)\sim  (x_i,y_j)\right\}\right|=1,\\
\kappa_2(x_i,y_j)= h^2q(x_i,y_j) + \deg((x_i,y_j)),&\mbox{if } \left|\left\{(x_r,y_s) \in \bar{V}_\bfnn\setminus V_\bfnn \, : \, (x_r,y_s)\sim  (x_i,y_j)\right\}\right|=2.
\end{cases}
\end{equation*}
See Figure \ref{fig:immersion_in_the_disk}.

The given graph-Laplacian $\Delta_{G_\bfnn}$ approximates the weighted operator $h^2\mathcal{L}$. Moreover, by Corollary \ref{cor:graph-Laplacian-distribution} it holds that 
\begin{equation*}
\{ \Delta_{G_\bfnn} \} \sim_\lambda g(x,y,\theta_1,\theta_2), \qquad (x,y,\theta_1,\theta_2) \in B_{\frac{1}{2}}\times[-\pi,\pi]^2
\end{equation*}
where
\begin{equation}\label{eq:distribution_FD_example}
g(x,y,\theta_1,\theta_2)=\left[1+(x-1/2)^2+(y-1/2)^2\right]\left(4 - 2\cos(\theta_1)  -2\cos(\theta_2)\right).
\end{equation}
By the symmetry of $g$ over $[-\pi, \pi]^2$ we can restrict it to $B_{1/2}\times [0,\pi]^2$ without affecting the validity of the identity \eqref{distribution:sv-eig}. If we consider now the monotone rearrangement $\tilde{g} : [0,1] \to [0, 10]$ of the symbol $g$ as in Definition \ref{def:rearrangment}, then we can see from  Table \ref{table:FD_table2} and Figure \ref{fig:FD_figure} that 
$$
\lim_{n\to\infty} \lambda_{k(\bfnn)}^{(\bfnn)} \to \tilde{g}\left(\frac{k(\bfnn)}{d_\bfnn}\right)
$$
for any index sequence $\{k(\bfnn)\}$, $1\leq k(\bfnn)\leq d_\bfnn$, such that $\frac{k(\bfnn)}{d_\bfnn}\to x \in [0,1]$ as $\bfnn\to \infty$, where $d_\bfnn$ is the dimension of the graph-Laplacian $\Delta_{G_\bfnn}$. We want to stress out that $d_\bfnn< n^2$ since $V_{\bfnn} \subsetneq \bar{V}_{\bfnn}$, but clearly it holds that $d_\bfnn \to \infty$ as $n\to \infty$, and the Hausdorff distance between the node set $V_\bfnn$ and the disk $B_{1/2}$ is going to zero.

Since $\tilde{g}$ does not posses an easy analytical expression to calculate, it has been approximated by an equispaced sampling of $g$ over $B_{1/2}\times [0,\pi]^2$ by $d_\bfnn$-points and then rearranging it in non-decreasing order. The approximation converges to $\tilde{g}$ as $n\to \infty$, see for example \cite{T86}. Finally, see Remark \ref{rem:rearrangment}.

For other applications of the GLT techniques for approximation of partial differential operators see \cite{BS18}

\begin{table}[H]
	\begin{minipage}[m]{0.54\textwidth}
\begin{tabular}{|l|c|ccc|} 
	\cline{3-5}
	\multicolumn{1}{l}{}               & \multicolumn{1}{l|}{} & \multicolumn{3}{c|}{\textbf{relative errors} }                       \\ 
	\cline{3-5}
	\multicolumn{1}{l}{}               &                       & $n=10$             & $n=50$             & $n=80$        \\ 
	\hline
	\multirow{4}{*}{$\frac{k(\bfnn)}{d_\bfnn}$ } & 0.1                   & 0.0788               &    0.0094            &     0.0020            \\ 
	\cline{2-5}
	& 0.5                   & 0.0055           & 3.3995e-04          & 1.2565e-04           \\ 
	\cline{2-5}
	& 0.8                   & 0.0100               & 7.2173e-04               & 3.9353e-05             \\ 
	\cline{2-5}
	& 1                     & 0.0443               & 0.0440              & 0.0360                 \\ 
	\hline
\end{tabular}\captionof{table}{Relative errors between the eigenvalue $\lambda_{k(\bfnn)}^{(\bfnn)}$ and the sampling of the monotone rearrangement $\tilde{g}\left(\frac{k(\bfnn)}{d_\bfnn}\right)$, i.e.: $\left|\frac{\lambda_{k(\bfnn)}^{(\bfnn)}}{\tilde{g}(k(\bfnn)/d_\bfnn)}-1\right|$. The index $k(\bfnn)$ of the eigenvalue $\lambda_{k(\bfnn)}^{(\bfnn)}$ is chosen such that $k(\bfnn)/d_\bfnn$ is constant for every fixed $n=10, 50, 80$. As it can be seen, the relative errors decrease as $n$ increases, in accordance with Theorem \ref{thm:discrete_Weyl_law}. The convergence to zero is not uniform and slower in some subintervals.}
\label{table:FD_table2}
\end{minipage}\hfill
\begin{minipage}[m]{0.42\textwidth}
\includegraphics[width=\textwidth]{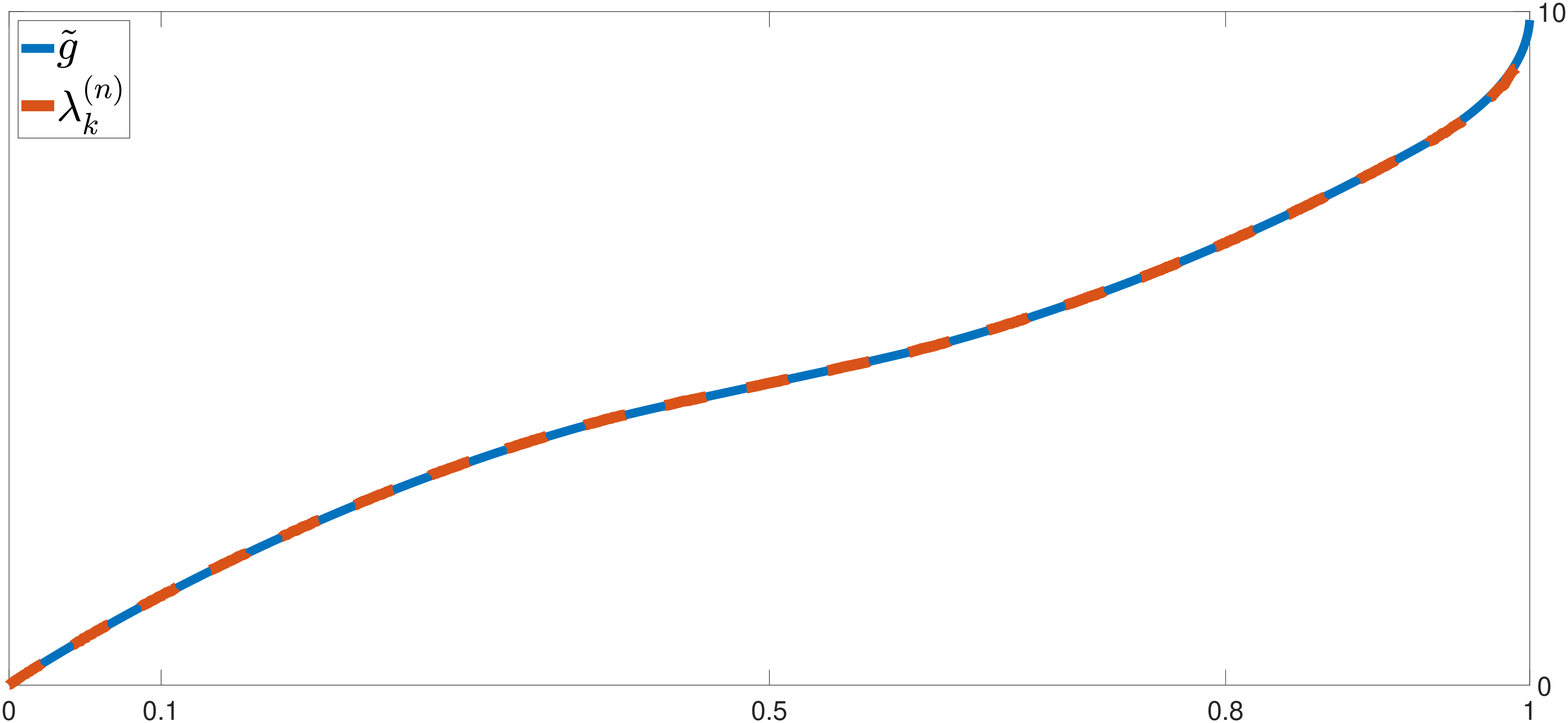}
\captionof{figure}{Plots of the monotone rearrangement $\tilde{g}$ (blue-continuous line) and the eigenvalues $\{\lambda_k^{(\bfnn)}\}$ (orange-dotted line), for $n=80$, of the graph-Laplacian $\Delta_{G_\bfnn}$ associated to the graph $G_\bfnn$ defined in \eqref{eq:5-FD-graph}. In this case we have that $d_\bfnn = 5140 < 80^2$.}
\label{fig:FD_figure}
\end{minipage}
\end{table}

 \subsection{Approximations of PDEs vs sequences of weighted diamond graphs: FEM}\label{ssec:appl-1}

Consider the model boundary value problem
\begin{equation}\label{problem}
\left\{\begin{array}{ll}
-\Delta u = f, & \mbox{in}\ \Omega,\\
u=0, & \rm{on\ }\partial\Omega,
\end{array}\right.
\end{equation}
where $\Omega=(0,1)$ and $f\in L^2(\Omega)$. We approximate \eqref{problem} by using the quadratic $C^0$ B-spline discretization on the uniform mesh with stepsize $ \frac{1}{n} $, where the basis functions are chosen as the $ C^0 $ B-spline of degree $ 2 $ defined on the knot sequence $ \{ \frac{1}{n+1}, \frac{1}{n+1}, \frac{2}{n+1}, \frac{2}{n+1},..., \frac{n}{n+1}, \frac{n}{n+1} \} $. 

Proceeding as in \cite{FEM-paper}, we trace the problem back to solving a linear system whose stiffness matrix reads as follows:
$$
A_{n}=\frac{n}{3}\left(\begin{array}{cccccccccc}
                                \cline{1-4} \multicolumn{1}{|c}{4} & -2 & \multicolumn{1}{|c}{} &  & \multicolumn{1}{|}{}      &        &    \\ 
                                   \multicolumn{1}{|c}{-2} &  8 & \multicolumn{1}{|c}{-2} & -2 & \multicolumn{1}{|}{0}   \\
                                   
                           \cline{1-6}        \multicolumn{1}{|c}{} & -2 & \multicolumn{1}{|c}{4} & -2 & \multicolumn{1}{|c}{}      &        &  \multicolumn{1}{|}{} \\ 
                                   \multicolumn{1}{|c}{} &  -2 & \multicolumn{1}{|c}{-2} & 8 & \multicolumn{1}{|c}{-2} & -2 &  \multicolumn{1}{|c}{}   &   &    \\  \cline{1-6}                                               
                                      &  & \ddots & & \ddots & & \ddots    \\ 
                                      &  &  & \ddots & & \ddots & & \ddots  \\ 
                                \cline{5-10}      &  &  &  &  \multicolumn{1}{|c}{} & -2 & \multicolumn{1}{|c}{4} & -2 & \multicolumn{1}{|c}{} & \multicolumn{1}{c|}{}    \\ 
                                      &  &  &  &  \multicolumn{1}{|c}{} &  -2 & \multicolumn{1}{|c}{-2} & 8 & \multicolumn{1}{|c}{-2} &  \multicolumn{1}{c|}{-2}     \\  \cline{5-10}                        
                                      &  &  &  &  &  & \multicolumn{1}{|c}{} & \multicolumn{1}{c|}{-2} & 4 & \multicolumn{1}{c|}{-2} \\
                                      &  &  &  &  &  & \multicolumn{1}{|c}{} & \multicolumn{1}{c|}{-2} & -2 & \multicolumn{1}{c|}{8} \\ \cline{7-10}
                                                        
                                  \end{array}\right)_{2n\times 2n}.
$$

We note that $\frac{1}{n} A_{n}$ can be seen as the graph-Laplacian of a $1$-level diamond Toeplitz graph with nonzero killing term. Namely, according to Definition \ref{def:graph-Laplacian}, we have that
$$
\frac{1}{n} A_n = K_n-W_{n,2}^G
$$
where $ K_n $ is the diagonal matrix given by
$$ \left(K_n \right)_{ii}= \begin{cases}
\frac{4}{3} & \text{ if } i \text{ is even} \\
\frac{8}{3} & \text{ otherwise}
\end{cases}
$$
and $ W_{n,2}^G $ is the adjacency matrix of the $1$-level diamond Toeplitz graph $ T_3^{G} \langle \left(1,L_1 \right) \rangle $, with
$$
G= \left( \begin{array}{cc}
0 & 2 \\
2 & 0
\end{array} \right)
$$ 
and
$$
L_1=\left( \begin{array}{cc}
0 & 0 \\
2 & 2
\end{array} \right).
$$
Combining now Proposition \ref{prop:diamond-distribution} and Corollary \ref{cor:graph-Laplacian-distribution} we get that $ \frac{1}{n} A_n $ has asymptotic spectral distribution with symbol function $ \bfff : [-\pi,\pi] \to \mathbb{C}^{2\times 2} $ given by
$$
\bfff(\theta) = \frac{1}{3} \left\{ \left( \begin{array}{cc}
4 & -2 \\
-2 & 8
\end{array} \right) + \left( \begin{array}{cc}
0 & -2 \\
-2 & -4
\end{array} \right) \cos (\theta) + \left( \begin{array}{cc}
0 & -2 \\
-2 & -4
\end{array} \right) i \sin (\theta)  \right\}.
$$

\begin{rem}
It is possible to study the multi-dimensional case of the problem in the example above using the fact that for every $ \bfmm, \bfss \in \mathbb{N}^d $ there exists a permutation matrix $ \Gamma_{\bfmm, \bfss} $ of dimension $ \prod_{j=1}^{d} m_j s_j $ such that
$$
T_{m_1}(\bfpp_1) \otimes T_{m_2}(\bfpp_2) \cdots \otimes T_{m_d}(\bfpp_d) = \Gamma_{\bfmm, \bfss} \left[ T_{\bfmm}(\bfpp_1(\theta_1) \otimes \cdots \otimes \bfpp(\theta_d)) \right] \Gamma_{\bfmm, \bfss}^T
$$
for any choice of trigonometric polynomials $ \bfpp : [-\pi,\pi] \to \mathbb{C}^{s_j \times s_j}, j=1,...,d$, as stated in Lemma 4 of \cite{FEM-paper}, where $ \otimes $ denotes the usual tensor product.
In the $d$-dimensional case (see, again, \cite{FEM-paper} and references therein), the discretizing matrix is given by
\begin{equation*}
A_{\bfnn} = \sum_{k=1}^{d} \left( \bigotimes_{r=1}^{k-1} \frac{1}{n_r} M_{n_r} \right) \otimes A_{n_k} \otimes \left( \bigotimes_{r=k+1}^{d} \frac{1}{n_r} M_{n_r} \right),
\end{equation*}
where $A_n$ is defined as in the example above and $M_n$ is a $2n \times 2n $ matrix with the same block Toeplitz structure as $A_n$ and, hence, an analogous symbol function, which we denote by $ \bfhh $. Assuming now that the multi-index $ \bfnn= \bfnnu n = \left(\nu_1 n,\nu_2 n,...,\nu_d n \right)$ for a fixed $ \bfnnu \in \mathbb{Q}^{d}_{>0} := \left \{ \left( \nu_1,...,\nu_d \right) \in \mathbb{Q} \, : \, \nu_1,...,\nu_d >0 \right \}$, it is immediate to see by the considerations above that $n^{d-2} A_\bfnn$ is the linear combination of graph-Laplacians of $d$-level diamond Toeplitz graphs with spectral distribution given by the following symbol function $ \bfgg:\left[ -\pi, \pi \right]^d \to \mathbb{C}^{2^d \times 2^d} $,
$$
\bfgg(\bftheta) = \sum_{k=1}^{d} c_k(\bfnnu) \left( \bigotimes_{r=1}^{k-1} \bfhh (\theta_{r}) \right) \otimes \bfff (\theta_{k}) \otimes \left( \bigotimes_{r=k+1}^{d} \bfhh (\theta_{r}) \right),
$$
with
$$
c_k(\bfnnu)= \frac{\nu_k}{\nu_1 \cdots \nu_{k-1} \nu_{k+1} \cdots \nu_d}, \, k=1,...,d.
$$
\end{rem}

 \subsection{Approximations of PDEs vs sequences of weighted d-level graphs: IgA approach}\label{ssec:appl-1-bis}

In this section we consider the approximation of the same differential equation \eqref{problem} by using the Galerkin B-splines approximation of degree $\nu$ in every direction: while the approximation of standard $\mathbb Q_\nu$ Lagrangian FEM  considered
in Section \ref{ssec:appl-1} leads to a symbol which is a linear trigonometric polynomial in every variable, but $\mathbb{C}^{\nu^d\times \nu^d}$ Hermitian matrix valued, here the symbol is scalar valued, but the degree of the trigonometric polynomial is much higher. For example, by means of cubic $C^2$ B-spline discretization the (normalized) stiffness matrix is given by

\begin{equation*}
\R^{n\times n} \ni A_n = \frac{1}{240}\begin{tikzpicture}[baseline=-\the\dimexpr\fontdimen22\textfont2\relax]
\matrix (m)[matrix of math nodes,left delimiter=(,right delimiter=), every node/.style={anchor=base,text depth=.5ex,text height=2ex,text width=2em}]
{	$360$ & $9$  & $\shortminus 60$ & $\shortminus 3$ & $0$  &$\cdots$ & &  & & & $0$ \\
	$9$ & $162$ & $\shortminus 8$ & $\shortminus 47$ & $\shortminus 2$ & $0$ & $\cdots$& &  & &0 \\
 $\shortminus 60$ & $\shortminus 8$ & $160$ &$\shortminus 30$  & $\shortminus 48$ & $\shortminus 2$ &0 &$\cdots$ &  & &0 \\
 $\shortminus 3$ & $\shortminus 47$ & $\shortminus 30$ & $160$ & $\shortminus 30$ & $\shortminus 48$ &$\shortminus 2$ &$0$& $\cdots$ & &0 \\
 $0$ & $\shortminus 2$ & $\shortminus 48$ &$\shortminus 30$  & $160$ & $\shortminus 30$ &$\shortminus 48$ &$\shortminus 2$ & $0$  & $\cdots$&0 \\
	 &  & $\ddots$ &$\ddots$ &$\ddots$& $\ddots$ &$\ddots$ &$\ddots$ &$\ddots$ & & \\
	  & &  & $\shortminus 2$ & $\shortminus 48$ &$\shortminus 30$  & $160$ & $\shortminus 30$ &$\shortminus 48$ &$\shortminus 2$ & $0$ \\
	  $0$ &  & $\cdots$ & $0$& $\shortminus 2$ & $\shortminus 48$ &$\shortminus 30$  & $160$ & $\shortminus 30$ & $\shortminus 47$ &$\shortminus 3$ \\
	$0$ &  &  &$\cdots$ &$0$ & $\shortminus 2$ & $\shortminus 48$ &$\shortminus 30$  & $160$ &$\shortminus 8$ &$\shortminus 60$ \\
	$0$ &  &  & &$\cdots$& $0$ & $\shortminus 2$ & $\shortminus 47$ &$\shortminus 8$  & $162$ &$9$  \\
	0 &  & & & &$\cdots$&0 & $\shortminus 3$ &$\shortminus 60$  & $9$ &$360$ \\
};
\begin{pgfonlayer}{myback}
\fhighlight[yellow!40]{m-1-1}{m-4-7}
\fhighlight[blue!40]{m-5-2}{m-7-10}
\fhighlight[yellow!40]{m-8-5}{m-11-11}
\end{pgfonlayer}
\end{tikzpicture}.
\end{equation*}
Observe that the principal $\R^{(n-4)\times (n-4)}$-submatrix (highlighted in blue) is an exact symmetric Toeplitz matrix while globally $A_n$ is not Toeplitz due to the presence of perturbations near the boundary points (highlighted in yellow). This behavior is influenced by the presence of BCs and the specific choice for the test-functions (B-spline with $C^3$ local regularity and $C^2$ global regularity). For those reasons, $A_n$ can not be representative of the graph-Laplacian of a graph in the form $G_n=\left(T_n\langle(t_1,w_1),\ldots,(t_m,w_m)\rangle, \kappa\right)$, even if it is clearly the graph-Laplacian for another kind of graph $G_n$ which does not own globally the symmetries of the graphs studied in sections \ref{sec:Diamond_Toeplitz} and \ref{sec:main}. Nevertheless, since the perturbations are local, it happens that $\{A_n\}_n \sim_\lambda f(\theta)$ where $f$ is the same symbol function of the graph-Laplacian of 
$$
G_n=  \left(T_n\left\langle \left(1,\frac{30}{240}\right),\left(1,\frac{48}{240}\right),\left(1,\frac{2}{240}\right)\right\rangle, \kappa\right), \quad \kappa\equiv 0,
$$
namely, $f(\theta)= \frac{160}{240} - \frac{60}{240}\cos(\theta) - \frac{96}{240}\cos(2\theta) - \frac{4}{240}\cos(3\theta)$. For a complete treatment of the study of the eigenvalue distribution for Galerkin B-splines approximations we refer to \cite{GSERSH18}, where many examples are provided along the exposure.

\begin{rem}
There is another quite important difference between this case and the case of Section \ref{ssec:appl-0}. In the FD case, the nodes of the graph were representative of the physical domain while in this case, even if the node set can be immersed in $[0,1]$, the nodes represent the base functions of the test-functions set. That said, given $A_n$, it will be of interest to calculate the corresponding weight function $w$ on the node set of the physical domain: this approach could lead some insight about the problem of the presence of a fixed number of outliers in the spectrum of $A_n$, see \cite[Chapter 5.1.2 p. 153]{IgA-book}.  
\end{rem}

\section{Conclusions, open problems, and future work}\label{sec:final}

We have defined general classes of graph sequences  having a grid
geometry with a uniform local structure  in a domain
$\Omega\subset [0,1]^d$, $d\ge 1$. With the only weak requirement
that $\Omega$ is Lebesgue measurable with boundary of zero Lebesgue measure, we have shown that  the
underlying sequences of adjacency matrices have a canonical
eigenvalue distribution, in the Weyl sense, with a symbol $f$
being a trigonometric polynomial in the $d$ Fourier variables: as
specific cases, we mention standard Toeplitz graphs, when
$\Omega=[0,1]$, and  $d$-level Toeplitz graphs when
$\Omega=[0,1]^d$, but also matrices coming from the approximation
of differential operators by local techniques, including Finite
Differences, Finite Elements, Isogeometric Analysis etc. In such a
case we considered block structures and weighted graphs, from the
perspective of GLT sequences, where the tools taken from the
latter field have resulted crucial for deducing all the asymptotic
spectral results. In particular, the knowledge of the symbol and
of basic analytical features have been employed for deducing a lot
of information on the eigenvalue structure, including precise
asymptotics on the gaps between the largest eigenvalues.

Many open problems remain, ranging from a deeper analysis of the matrix sequences arising from different families of Finite Element approximations of multidimensional differential problems to the study of the convergence features of the ordered asymptotic spectra to the rearrangement of the corresponding symbol (see also the study and discussions in \cite{GSERSH18,Bianchi}). 

\appendix
\section{Appendix}\label{appendix}
Fix a square matrix sequence $\{X_{\bfnn,\nu}\}_\bfnn$ of dimension $d_\bfnn$, with symbol function $\bfgg : D\subset \R^m \to \mathbb{C}^{\nu\times \nu}$ as Definition \ref{def:eig-distribution}. Observe that $\bfgg$ is not unique and in general not univariate. To avoid this, we will soon introduce in Definition \ref{def:rearrangment} the notion of monotone rearrangement of the symbol. In order to simplify the notations and since all the cases we investigate in this paper can be lead back to this situation, we make the following assumptions:

\vspace{0.4cm}
\noindent\textbf{Assumptions}
\begin{enumerate}[label={\upshape(\bfseries AS\arabic*)},wide = 0pt, leftmargin = 3em]
	\item $\overline{D}$ is compact and of the form $\overline{\Omega} \times [-\pi,\pi]^d$ with $\overline{\Omega}\subseteq [0,1]^{d}$, and therefore $m=2d$;\label{AS1}
	\item $\bfgg(\boldsymbol{y})= \bfgg(\bfxx,\bftheta)= p(\bfxx)\bfff(\bftheta),$ with $(\bfxx,\bftheta)\in \Omega\times (-\pi,\pi)^d$ and $p: \Omega \to \R$, $\bfff : [-\pi,\pi]^d \to \mathbb{C}^{\nu\times\nu}$;\label{AS2}
	\item $p: \Omega \to \R$ piecewise continuous;\label{AS3}
	\item every component $f_{i,j}: [-\pi,\pi]^d \to \mathbb{C}$ of $\bfff$ is continuous;\label{AS4}
	\item $\bfff(\bftheta)$ has real eigenvalues for every $\bftheta \in [-\pi,\pi]^d$.\label{AS5}
\end{enumerate}

Since we are assuming the eigenvalues to be real, then for convenience notation we will order the eigenvalue functions $\lambda_k\left(p(\bfxx)\bfff(\bftheta)\right)$ by magnitude, namely $\lambda_1\left(p(\bfxx)\bfff(\bftheta)\right)\leq \ldots\leq \lambda_\nu\left(p(\bfxx)\bfff(\bftheta)\right)$. This kind of ordering could affect the global regularity of the eigenvalue functions, but it does not affect the global regularity of the monotone rearrangement of the symbol. Nevertheless, by well-known results (see \cite{Kato}), items \ref{AS3} and \ref{AS4} imply that $\lambda_k\left(p(\bfxx)\bfff(\bftheta)\right)$ is at least piecewise continuous for every $k=1,\ldots,\nu$. We have the following result.

\begin{prop}\label{prop:d-block-symbol_rearrangment}
Suppose that $\{X_{\bfnn,\nu}\}_\bfnn \sim_\lambda \bfgg(\bfxx,\bftheta)= p(\bfxx)\bfff(\bftheta)$ as in Definition \ref{distribution:sv-eig}, that is
	\begin{equation*}
\lim_{\bfnn\rightarrow \infty}\Sigma_{\lambda}(F,X_{\bfnn,\nu})=\frac1{\mu_m(D)}\iint_D \sum_{k=1}^\nu  F\left(\lambda_k(p(\bfxx)\bfff(\boldsymbol{\theta}))\right)d\mu_m(\bfxx,\boldsymbol{\theta}),\qquad\forall F\in C_c(\mathbb R),
\end{equation*}
where $\lambda_k(p(\bfxx)\bfff(\boldsymbol{\theta})),\ k=1,\ldots,\nu,$ are the eigenvalues of $p(\bfxx)\bfff(\boldsymbol{\theta})$. Then
\begin{equation}\label{distribution:sv-eig2}
\lim_{\bfnn\rightarrow \infty}\Sigma_{\lambda}(F,X_{\bfnn,\nu})=\frac{1}{\mu_{m}(\hat{D})}\iint_{\hat{D}}   F\left(p(\bfxx)\sum_{k=1}^\nu \lambda_k\left(\bfff_k(\bftheta)\right)\right)d\mu_{m}(\bfxx,\bftheta),\qquad\forall F\in C_c(\mathbb R),
\end{equation} 
where
\begin{subequations}
\begin{equation}\label{subeq:prop_d-block-symbol_rearrangmentA}
 \hat{D}= \Omega \times \left( \bigcup_{k=1}^\nu I_k\right), I_k=\underbrace{\left[\frac{(2(k-1)-\nu)\pi}{\nu},\frac{(2k-\nu)\pi}{\nu}\right]\times\cdots \times \left[\frac{(2(k-1)-\nu)\pi}{\nu},\frac{(2k-\nu)\pi}{\nu}\right]}_{d-\mbox{times}};
 \end{equation}
 \begin{equation}\label{subeq:prop_d-block-symbol_rearrangmentB}
	 \bfff_k(\bftheta)=\begin{cases}
	\bfff(\nu\bftheta - (2k-1-\nu)\boldsymbol{\pi}) & \mbox{if } \bftheta\in I_k,\\
	\bfzero & \mbox{otherwise}.
	\end{cases}
	\end{equation}
\end{subequations}
\end{prop}
{\bf Proof}
By the monotone convergence theorem, since $F$ is limit of a monotone sequence of step functions, it is sufficient to prove the statement for $F=\mathds{1}_{E}$, with $E$ a measurable subset of $\textnormal{supp}(F)$. Then it holds trivially that
\begin{equation*}
\frac1{\mu_m(D)}\iint_D \sum_{k=1}^\nu  \mathds{1}_{E}\left(\lambda_k(p(\bfxx)\bfff(\boldsymbol{\theta}))\right)d\mu_m(\bfxx,\boldsymbol{\theta}) = \sum_{k=1}^\nu\frac1{\mu_m(D)}\iint_{D}   \mathds{1}_{E_k}(\bfxx,\bftheta)d\mu_m(\bfxx,\boldsymbol{\theta}),
\end{equation*}
with $E_k=\left\{(\bfxx,\bftheta) \in D \, : \, \lambda_k(p(\bfxx)\bfff(\boldsymbol{\theta})) \in E \right\}$. The rest follows from the changes of variables 
$$
(\bfxx,\bftheta) \mapsto (\bfxx, \nu^{-1}\left[\bftheta + (2k-1-\nu)\boldsymbol{\pi}\right]), \qquad \hat{E}_k=\left\{(\bfxx,\bftheta) \in \Omega \times I_k \, : \, \lambda_k(p(\bfxx)\bfff(\boldsymbol{\theta})) \in E \right\}
$$
and from the easy fact that $$
\sum_{k=1}^\nu  \mu_m\left( \hat{E}_k\right)=  \mu_m\left( (\bfxx,\bftheta) \in \hat{D} \, : \, p(\bfxx)\sum_{k=1}^\nu \lambda_k(\bfff_k(\boldsymbol{\theta})) \in E\right).
$$
\hfill \ \, $\bullet$ \ \\

\begin{rem}\label{rem:appendix_1}
Trivially, the maps $I_k \ni \bftheta \mapsto \nu\bftheta - (2k-1-\nu)\boldsymbol{\pi}$ are diffeomorphism between $I_k$ and $[-\pi,\pi]^d$. Therefore, the image set of $\sum_{k=1}^\nu \lambda_k\left(\bfff_k(\bftheta)\right)$ over $I_k$ is exactly the image set of $\lambda_k\left(\bfff(\bftheta)\right)$ over $[-\pi,\pi]^d$.
\end{rem}

The next definition of {\itshape{monotone rearrangement}} \`{a} la Hardy-Littlewood is crucial for the understanding of the asymptotic distribution of the eigenvalues of $\{X_{\bfnn,\nu}\}_\bfnn$.

\begin{defn}\label{def:rearrangment}
	Let $\nu \geq 1$ and using the same notations of Proposition \ref{prop:d-block-symbol_rearrangment}, define
	$$
	R_\bfgg = \left\{ p(\bfxx)\sum_{k=1}^\nu \lambda_k\left(\bfff_k(\bftheta)\right) \, : \, (\bfxx,\bftheta) \in \overline{\Omega} \times [-\pi,\pi]^d\right\}.
	$$
	
Let $\tilde{g} :  [0,1]\to [\min R_\bfgg , \max R_\bfgg ]$ be such that 
\begin{subequations}
		\begin{equation}\label{eq:rearrangment}
		\tilde{g}(x) = \inf\left\{ t \in [\min R_\bfgg , \max R_\bfgg ]\,:\, \phi(t)\geq x\right\}  
		\end{equation}
		where 
		\begin{equation}\label{eq:rearrangment2}
		\phi : [\min R_\bfgg ,\max R_\bfgg ] \to [0,1], \qquad \phi(t) := \frac{1}{\mu_m(\hat{D})}\mu_m \left(\left\{(\bfxx,\bftheta)\in \hat{D} \, : \, p(\bfxx)\sum_{k=1}^\nu \lambda_k\left(\bfff_k(\bftheta)\right) \leq t  \right\}\right).
		\end{equation}
\end{subequations}
		We call $\tilde{g}$ the \textnormal{monotone rearrangement} of $\bfgg$. 
\end{defn}
Clearly, $\tilde{g}$ is well-defined, univariate, monotone strictly increasing and right-continuous. Within our assumptions on $\bfgg$, it is then easily possible to extend \cite[Theorem 3.4]{DBFS93} for this multi-variate matrix-valued case, and it holds that
\begin{enumerate}[(i)]
	\item $\tilde{g}(0)= \min R_\bfgg$, $\tilde{g}(1)= \max R_\bfgg$;
	\item $\lim_{\bfnn\rightarrow \infty}\Sigma_{\lambda}(F,X_{\bfnn,\nu})=\int_0^1 F(\tilde{g}(x))d\mu_1\left(x\right).$
\end{enumerate}

We have the following results, see \cite[Section 3.4]{Bianchi}. The statements and proofs are exactly the same, we just adjusted them to fit with the multi-index notation we are adopting here.

\begin{thm}\label{thm:discrete_Weyl_law}
		Let $\tilde{g}: [0,1] \to [\min R_\bfgg,\max R_\bfgg]$ be the monotone rearrangement of a spectral symbol $\bfgg$ of the matrix sequence $\left\{ X_{\bfnn,\nu} \right\}$ and let $\left\{\lambda_k^{(\bfnn)} \right\}_{k=1}^{d_\bfnn}$ be the collection of eigenvalues of the matrices $X_{\bfnn,\nu}$, sorted in non-decreasing order. Let $\tilde{g}$ be piecewise Lipschitz continuous. Then 
		\begin{equation}\label{eq:discrete_Weyl_law1}
		\lim_{\bfnn \to \infty}\frac{\left| \left\{ k=1,\ldots,d_\bfnn \, : \, \lambda_k^{(\bfnn)}  \leq t \right\} \right|}{d_\bfnn} = \tilde{g}^{-1}(t).
		\end{equation}
		In particular, let $k=k(\bfnn)$ be such that $k(\bfnn)/d_\bfnn \sim x$ as $\bfnn\to \infty$ for a fixed $x \in [0,1]$ and let $\lambda_{k(\bfnn)}^{(\bfnn)}\in \tilde{g}\left([0,1]\right)=R_\bfgg$ definitely. Then
		\begin{align}
		&\lambda_{k(\bfnn)}^{(\bfnn)} \sim \tilde{g}\left(\frac{k(\bfnn)}{d_\bfnn}\right) \qquad \mbox{as } \bfnn\to \infty,\label{eq:discrete_Weyl_law2_1}\\
		&\left(\frac{k(\bfnn)}{d_\bfnn}, \lambda_{k(\bfnn)}^{(\bfnn)}  \right) \to \left(x, \tilde{g}(x)\right) \qquad \mbox{as } \bfnn\to \infty. \label{eq:discrete_Weyl_law2_2}
		\end{align}
\end{thm}

Notice that the relations \eqref{eq:discrete_Weyl_law2_1}-\eqref{eq:discrete_Weyl_law2_2} tell us that $\lambda_{k(\bfnn)}^{(\bfnn)}$ converges to $\tilde{g}(x)$, but they do not say anything about the rate of convergence, which can be slow.

\begin{cor}\label{cor:weak_clustering}
	In the same hypothesis of Theorem \ref{thm:discrete_Weyl_law}, it holds that
	$$
	\lim_{n \to \infty}\frac{\left| \left\{ k=1,\ldots,d_\bfnn \, : \, \lambda_k^{(\bfnn)} \notin R_\bfgg \right\} \right|}{d_\bfnn} = 0,
	$$
	that is, the number of possible outliers is $o(d_n)$.
\end{cor}

\begin{cor}\label{cor:differentiability}
		In the same hypothesis of Theorem \ref{thm:discrete_Weyl_law}, let $\tau: [\min R_\bfgg,\max R_\bfgg] \to \R$ be a differentiable real function and let $\{k(\bfnn)\}$ be a sequence of integers such that
		\begin{enumerate}[(i)]
			\item $\frac{k(\bfnn)}{d_\bfnn} \to x_0 \in [0,1]$;
			\item $\lambda_{k(\bfnn)}^{(\bfnn)}> \lambda_{k(\bfnn)-1}^{(\bfnn)}\in [\min R_\bfgg,\max R_\bfgg]$ definitely for $\bfnn \to \infty$.
		\end{enumerate} 
		Then
		$$
		d_\bfnn\left[\tau\left(\lambda_{k(\bfnn)}^{(\bfnn)} \right) - \tau\left(\lambda_{k(\bfnn)-1}^{(\bfnn)} \right)\right] \to \partial_x\left(\tau(\tilde{g}(x))\right)_{|x=x_0} \quad \mbox{a.e.} \qquad \mbox{as } \bfnn \to \infty.
		$$
\end{cor}

\begin{rem}\label{rem:rearrangment}
It may often happen that $\tilde{g}$ does not have an analytical expression or it is not feasible to calculate, therefore it is often needed an approximation. The simplest and easiest way to obtain it is by mean of sorting in non-decreasing order a uniform sampling of the original symbol function $p(\bfxx)f(\bftheta)$, in the case of real-valued symbol, or of sorting in non-decreasing order uniforms samplings  of $ p(\bfxx) \lambda_k\left(\bfff(\bftheta)\right)$ for $k=1,\ldots,\nu$, in the case of a matrix-valued symbol, by Remark \ref{rem:appendix_1}. See as a references \cite[Section 3]{GSERSH18} and \cite[Remark 2]{GMS18}. These approximations converge to $\tilde{g}$ as the mesh-refinement goes to zero, see \cite{T86}.
\end{rem}


\end{document}